\documentclass{article}
\usepackage{pstricks}
\usepackage[english]{babel}
\usepackage{amsmath,enumerate, amsfonts, amssymb,amsthm, xypic}
\DeclareMathOperator{\R}{\mathbf R}
\DeclareMathOperator{\C}{\mathbf C}
\DeclareMathOperator{\Q}{\mathbf Q}
\DeclareMathOperator{\N}{\mathbf N}
\DeclareMathOperator{\Z}{\mathbf Z}

\DeclareMathOperator{\st}{st}
\DeclareMathOperator{\lk}{lk}

\DeclareMathOperator{\Int}{Int}
\DeclareMathOperator{\id}{id}

\DeclareMathOperator{\Ima}{Im}
\DeclareMathOperator{\graph}{graph}
\DeclareMathOperator{\Dom}{Dom}
\DeclareMathOperator{\dis}{dis}

\begin{document}
\title{O-MINIMAL HAUPTVERMUTUNG FOR POLYHEDRA II}
\author{MASAHIRO SHIOTA}
\maketitle
\begin{abstract}
Hilbert initiated the finitary standpoint in foundations of mathematics.
From this standpoint, we allow only a finite number of repetitions of elementary operations when we construct objects and morphisms.
Let $S_n(X)$ denote the family of subsets of $\R^n$ obtained in this way, $n\in\N$, when we start from a subset $X$ of $\R^m$.
Then we assume that any element of $S_1(X)$ has only a finite number of connected components because for any closed subset $Y$ of $\R^n$ there exists $X$, which does not satisfy this assumption, such that $Y\in S_n(X)$ and the family $\{S_n(X):n\in\N\}$ is not interesting.
We call $X$ {\it tame} if the assumption is satisfied, and define a {\it tame morphism} in the same way.\par
In this paper we will show that a tame $C^0$ manifold is tamely homeomorphic to the interior of a compact PL (=piecewise linear) manifolds possibly with boundary and such a PL manifold possibly with boundary is unique up to PL homeomorphisms in the sense that if $M_1$ and $M_2$ are such PL manifolds possibly with boundary then $M_1$ and $M_2$ are PL homeomorphic (Theorem 1).
We modify this to Theorem 2 so that argument of model theory works, and we prove it.
We also consider the $C^r$ case, $0<r$ (Theorems 1$'$ and 2$'$).
\end{abstract}
\footnotetext{2010 {\it Mathematics Subject Classification}, 03C64, 57Q05, 57Q25.\\
\hspace{5.5mm}{\it Key words and phrases}. PL manifolds; compactification; Hauptvermutung; o-minimal structures.}
\section{Introduction}
Let $\N$ denote nonnegative integers.
For a subset $X$ of a Euclidean space $\R^m$, consider a family $S$ of $S_n(X),\ n\in\N$, such that each $S_n(X)$ consists of subsets of $\R^n$,\\ 
(0) $X\in S_m(X)$,\\
(i) {\it $S_n(X)$ is a Boolean algebra of subsets of $\R^n$,}\par\noindent
(ii) {\it if $Y\in S_n(X)$, then $\R\times Y$ and $Y\times\R$ are elements of $S_{n+1}(X)$,}\par\noindent
(iii) {\it the sets $\{(x,y,z)\in\R^3:x+y=z\}$ and $\{(x,y,z)\in\R^3:x y=z\}$ are elements of $S_3(X)$,}\par\noindent 
(iv) {\it the set $\{(x,y)\in\R^2:x<y\}$ is an element of $S_2(X)$,} and\par\noindent
(v) {\it if $Y\in S_{n+1}(X)$, then the image of $Y$ under the projection of $\R^{n+1}$ onto the first $n$ coordinates is an element of $S_n(X)$}.\par
We naturally give a partial order to the family of all $S$'s, and let $\{S_n(X)\}_{n\in\N}$ denote the smallest element in the family for simplicity of notation.
Thus starting from $X$ and the sets in (iii) and (iv) we obtain $\cup_{n\in\N}S_n(X)$ by finite repetitions of operations in (i), (ii) and (v).\par
We call $X$ {\it tame} if\\
(vi) {\it any element of $S_1(X)$ has only a finite number of connected components.}\vskip1mm\noindent
We call a map $f:X\to Y$ between tame sets {\it tame} if $\graph f$ is tame, and say that $X$ and $Y$ are {\it tamely homeomorphic} if there is a tame homeomorphism $f$.
We easily see that if $f$ is a tame homeomorphism then $f^{-1}$ is tame.
A {\it tame $C^0$ manifold} (over $\R$) is both a tame set and a $C^0$ manifold having an atlas $\{(U_\alpha,\psi_\alpha)\}_{\alpha\in A}$ such that $\psi_\alpha:U_\alpha\to\R^n$ are tame.
Examples of a tame set and a tame $C^0$ manifold are an algebraic complex (or real) variety in some $\C^n\,(=\R^{2n})$ (or $\R^n$) and a compact real analytic variety in some $\R^n$, and an example of a non-tame set is the graph of the function $y=\sin x$ because the intersection of the graph and the line $y=0$ is $\pi\Z\times\{0\}$.
A compact Euclidean polyhedron (i.e., a finite union of simplices in some $\R^n)$ is tame, but a Euclidean polyhedron is not necessarily tame.
For example, $\N$ in $\R$ is not tame.\par
The term ``\,tame\," is derived from two facts.
One is that we allow only a finite number of repetitions of the above operations (i), (ii) and (v).
Hence a wild manifold is never tame in our sense because we need a countable number of repetitions of pasting of coordinate neighborhoods to construct a wild manifold.
The other is that elements of $S_n(X)$ keep the well-known tame properties of semialgebraic sets.
(Here a semialgebraic set is defined to be a finite union of sets $\{x\in\R^n:f_i(x)\dag_i0,\,i\in I\}$, where $\dag_i$ is either = or $>$, $I$ is a finite set and $f_i$ are polynomial functions on $\R^n$ with coefficients in $\R$.)
For example, we can prove that a knot of any dimension is a tame knot (i.e., equivalent to a polyhedral knot) if it is tame in our sense by Triangulation theorem of definable sets explained below.\par
In this paper we will prove the following:\vskip2mm\noindent
{\bf Theorem 1}. {\it A compact tame topological manifold $M$ admits a unique PL manifold structure, i.e., $M$ is tamely homeomorphic to a compact PL manifold $M_1$ and $M_1$ is unique up to PL homeomorphisms in the sense that if $M_2$ is another PL manifold tamely homeomorphic to $M$ then $M_1$ and $M_2$ are PL homeomorphic.
A noncompact tame topological manifold is tamely homeomorphic to the interior of a compact PL manifold $M_1$ with boundary, and $M_1$ is unique up to PL homeomorphisms.}\vskip2mm
A manifold is called a {\it manifold with boundary} if it has boundary, and a {\it manifold} is a manifold without boundary.
We apply a theory of o-minimal structures, which is a notion of model theory, to prove this theorem because we can study topological properties of tame sets and maps by the theory.
Model theorists: Knight, Pillay and Steinhorn introduced o-minimal structures (\cite{P-S} and \cite{K-P-S}), and o-minimal structures have been extensively studied ever since.
Let $R$ denote an ordered field.
A {\it definable open interval} of $R$ is a subset of $R$ of the form $\{x\in R:a<x<b\}$ for some $a,b\in R\cup\{\pm\infty\}$.
Hence when $R$ is the rational number field $\Q$, we do not call the set $\{x\in\Q:-\pi<x<\pi\}$ an interval of $\Q$ because the end points are not elements of $\Q$.
We cannot define such a set by a formula of first-order logic.
An {\it o-minimal structure} over $R$ is a sequence $\{S_n\}_{n\in\N}$ such that for each $n\in\N$,\\
(i) {\it $S_n$ is a Boolean algebra of subsets of $R^n$,}\par\noindent
(ii) {\it if $X\in S_n$, then $R\times X$ and $X\times R$ are elements of $S_{n+1}$,}\par\noindent
(iii) {\it the sets $\{(x,y,z)\in R^3:x+y=z\}$ and $\{(x,y,z)\in R^3:x y=z\}$ are elements of $S_3$,}\par\noindent 
(iv) {\it the set $\{(x,y)\in R^2:x<y\}$ is an element of $S_2$,}\par\noindent
(v) {\it if $X\in S_{n+1}$, then the image of $X$ under the projection of $R^{n+1}$ onto the first $n$ coordinates is an element of $S_n$,} and\par\noindent
(vi) {\it an element of $S_1$ is a finite union of points and definable open intervals.}\vskip1mm
The simplest example of an o-minimal structure is the families of {\it semialgebraic} sets in $R^n,\ n\in\N$.
An element of $S_n$ is called {\it definable}, and a map between definable sets is called {\it definable} if its graph is definable.
When we need to clarify the o-minimal structure, we say that a set or a map is $\{S_n\}$-{\it definable}.
We endow a topology on $R$ and then $R^n$ as in the case of real numbers so that an open set in $R$ is a union of definable open intervals, and we say that two definable sets are {\it definably homeomorphic} if there is a definable homeomorphism between them.
We call $R^n$ a {\it Euclidean space} if no confusion is possible.
We naturally define a $C^0$ manifold, a polyhedron and a PL manifold over $R$.
However, we are interested only in definable ones.
The term ``\,a definable set\," means that the set is defined in some $R^n$ by a formula of first-order logic.
For example, a definable open interval means an open interval defined by its end points.
We easily see that $R$ is {\it definably connected}, i.e., $R$ does not contain a proper definable open and closed subset.\par
A {\it definable} $C^r$ {\it manifold} (over $R$) is both a definable subset of some $R^n$ and a $C^r$ manifold having an atlas $\{(U_\alpha,\psi_\alpha)\}_{\alpha\in A}$ of definable, i.e., $U_\alpha$ and $\psi_\alpha$ are definable, $0\le r<\infty$.
Here a definable $C^r$ manifold, $0<r<\infty$, is well defined for the following reason:
Differentiability of an $R$-valued function on $R$ at a point $x_0$ is described by a formula as follows: for some $c\in R$ and any $\epsilon>0$ there exists $\delta>0$ such that
$$|f(x)-f(x_0)-c(x-x_0)|\le\epsilon|x-x_0|\ \ \text{for any }x\in R\text{ with }|x-x_0|<\delta.$$
Hence by elementary argument of logic we see that given a definable function $g$ on $R^n$, the subset of $R^n$ where $g$ is of class $C^r$ is definable, $0\le r<\infty$, and the first derivatives of $g$ are definable if $g$ is of class $C^1$ everywhere.
In the same way, the following sets and function are shown to be definable: the closure $\overline X$ of a definable set $X$ in its ambient space $R^n$, the function on $R^n$ measuring the distance from $X$ and the set of critical points (values) of a definable $C^r$ map between definable $C^r$ manifolds.
It is also easy to see that the composite of two definable $C^r$ maps between definable open subsets of $R^n$ is a definable $C^r$ map.
Thus a definable $C^r$ manifold is well defined for positive $r$.\par
A {\it definable polyhedron} and a {\it definable PL manifold} are a Euclidean polyhedron and a Euclidean PL manifold contained and definable in some $R^n$.
Note that a definable set is always contained in some $R^n$ and we regard two definable sets as the same set if and only if their ambient spaces $R^n$ are the same and they coincide as the subsets.
From now on, a simplicial complex, a polyhedron and a PL manifold are always Euclidean ones.\par
Fundamental topological properties of definable sets are explained in \cite{S2} and \cite{v}.
In this series we further study topology of definable sets and definable $C^0$ maps, and one of the main theorems is the following:\vskip2mm\noindent
{\bf Theorem 2}. {\it A definable $C^0$ manifold over $R$ is definably homeomorphic to the interior of a unique compact PL manifold possibly with boundary.}\vskip2mm
Here we call a definable set in $R^n$ {\it compact} if it is bounded and closed in $R^n$, a polyhedron in $R^n$ a {\it compact polyhedron} if it is a finite union of simplices in $R^n$ and a map between compact polyhedra {\it PL} if its graph is not only a polyhedron but also a compact polyhedron.
In \cite{S3} we have seen the following facts: there exists a polyhedron closed and bounded in $R^n$ which cannot be a finite union of simplices (p.\,167); a compact definable set or a compact polyhedron is not compact in the usual sense unless $R$ is the real number field $\R$ (p.\,171); a polyhedron in $R^n$ is a compact polyhedron if and only if it is a compact definable set (p.\,167); the simplicial approximation theorem does not necessarily hold over $R$ (p.\,166); the simplicial homotopy theorem holds instead (Lemma 3.1 in \cite{S3}).\par
We state Theorem 2 over general $R$ but not $\R$.
One reason is that we construct a theory of manifolds and maps after starting from axioms as few as possible and we do not want to apply special properties of real numbers as axiom because we need an infinite number of repetitions of testing whether such a special property is satisfied or not, e.g.~compactness as shown below.
Indeed, if we adopt the special properties as ones of the axioms then most argument of model theory, particularly, the proof of this paper, do not work.
(Fewer axioms give more results.)
The second reason is that we expect applications of our theory to the other domains.
It is E.~Artin who first used general $R$ when he solved Hilbert's seventeenth problem.
There are many such applications to problems of real algebraic geometry (see \cite{B-C-R}).
We also see an application to topology in \cite{F-S}, where the field of convergent Puiseux series is used to study Milnor fibers.\par
We prove Theorem 2 by a machinery whose input is a definable $C^0$ manifold, the operations are only the operations which appear in (i), (ii) and (v), and the output is the couple of a PL manifold possibly with boundary and a definable homeomorphism between the $C^0$ manifold and the interior of the PL manifold possibly with boundary.
We call such a proof by a machine {\it constructive} (see p.\,166 in \cite{S3} for the details).
We cannot constructively show many properties of real numbers.
An example is compactness of a closed bounded set in $\R^n$ in the usual sense.
To be precise, given an open covering of a closed bounded subset of $\R^n$, we have no constructive methods to choose a finite subcovering because the machine, if exists, depends on the covering.
Another example is the Archimedean property of real numbers.
Given positive real numbers $x$ and $y$, there is no constructive methods to find a natural number $n$ such that $x<n y$ for the same reason.
Our proof in this paper is always constructive.
If a statement is constructively proved over some ordered field with some o-minimal structure, then in many cases it holds over any ordered field with any o-minimal structure.
This is the case for Theorem 2.\par
In a machinery we admit only finite repetition of operations.
This causes some special phenomena of definable manifolds and definable homeomorphisms as in Theorem 2, which never occur for general manifolds and homeomorphisms.
Let us remember four kinds of manifolds: first, a topological manifold which does not admit a PL manifold structure (see p.\,194 in \cite{K-S}), secondly, two PL manifolds which are homeomorphic but not PL homeomorphic (a counterexample to the Hauptvermutung for manifolds) (ibid.), thirdly, a noncompact $C^0$ (or PL) manifold which cannot be (PL) homeomorphic to the interior of a compact $C^0$ (PL) manifold with boundary (a noncompactifiable manifold) (e.g.~$\R^2-\Z^2$), fourthly, two compact $C^0$ (PL) manifolds with boundary which are not (PL) homeomorphic but whose interiors are (PL) homeomorphic (e.g.~a disk and a compact contractible manifold which is not a disk \cite{M}).
In o-minimal topology of definable manifolds and definable $C^0$ maps from Hilbert's standpoint, there are no such ``\,wild\," phenomena.
This is what we would argue in the series and, particularly, in this paper.\par
Tame sets are similar but different from definable sets in o-minimal structures over $\R$.
A definable set in an o-minimal structure over $\R$ is clearly tame, and for a tame set $X$, the family $\{S_n(X)\}_n$ is an o-minimal structure over $\R$.
Hence a set in a Euclidean space is tame if and only if the set is definable in some o-minimal structure over $\R$, and, moreover, a map $f:X\to Y$ between subsets of Euclidean spaces is tame if and only if $X,\ Y$ and $\graph f$ are definable in some one o-minimal structure over $\R$.
However, for two tame sets $X$ and $Y$ in Euclidean spaces there does not necessarily exist an o-minimal structure over $\R$ where $X$ and $Y$ are both definable \cite{R-S-W}.
Hence it is not trivial that a tame $C^0$ manifold is a definable $C^0$ manifold in some o-minimal structure over $\R$, although it is easy to prove.
(Remember that a polyhedral $C^0$ manifold in some Euclidean space is not necessarily a PL manifold, e.g.~the double suspension of the Mazer homology 3-sphere and Freedman's $E_8$-manifold.)
We will prove this.
Then Theorem 1 follows from Theorem 2.\par
Theorem 2 in the compact case is an easy consequence of Triangulation theorem of definable sets and Uniqueness theorem of definable triangulation.
The first theorem was shown in \cite{S3} and is a refinement of the well-known one, and the second is Theorem 1.1 in \cite{S3}.
We explain both the theorems in the next section.
In the noncompact case, Uniqueness theorem is false.
We will choose triangulations of noncompact definable sets, which we call {\it standard} triangulations, so that uniqueness holds (Theorem 4).
For the proof of Theorem 4 we use Triangulation theorem of definable $C^0$ functions (Theorem 3.2 in \cite{S3}), which was the key lemma of the proof of Uniqueness theorem of definable triangulation.\vskip1mm
The $C^r$ cases of Theorems 1 and 2 are the following, although we stray off the subject of the series:\vskip2mm\noindent
{\bf Theorem 1$'$}. {\it A tame $C^r$ manifold is tamely $C^r$ diffeomorphic to a nonsingular algebraic variety defined by polynomials with coefficients in $\R_{\rm alg}$=real algebraic numbers, uniquely up to Nash (i.e., semialgebraic and $C^\infty$) diffeomorphisms, $0<r<\infty$.
In the noncompact case, it is tamely $C^r$ diffeomorphic to the interior of a unique compact Nash manifold with boundary.}\vskip2mm\noindent
{\bf Theorem 2$'$}. {\it A definable $C^r$ manifold over $R$ is definably $C^r$ diffeomorphic to a nonsingular algebraic variety defined by polynomials with coefficients in $\R_{\rm alg}$ uniquely up to Nash diffeomorphisms, $0<r<\infty$.
In the noncompact case, it is also definably $C^r$ diffeomorphic to the interior of a unique compact Nash manifold with boundary.}\vskip2mm
The following are known: Any compact $C^r$ manifold over $\R$ is $C^r$ diffeomorphic to a nonsingular algebraic variety (Nash-Tognoli); any Nash manifold over $R$ is Nash diffeomorphic to a nonsingular algebraic variety \cite{C-S1}; in the semialgebraic case, Theorem 2$'$ hold for $r=\infty$ by \cite{C-S1} and \cite{C-S2}; a Nash manifold and a Nash map between Nash manifolds over $\R$ is of class $C^\omega$.
Note that we do not know whether Theorem 2$'$ holds in the case $r=\infty$ because for such $r$ we need infinite procedures to construct both of the manifolds with boundary and the diffeomorphisms and that Theorem 1$'$ quickly follows from Theorem 2$'$ as Theorem 1 does from Theorem 2.\vskip1mm
We compare two o-minimal structures over two ordered fields.
Let $R$ and $R'$ be ordered fields such that $R$ is contained in $R'$, and $\{S_n\}$ and $\{S'_n\}$ be o-minimal structures over $R$ and $R'$, respectively, such that for any $X\in S_n$ there exists $X'\in S'_n$ such that $X'\cap R^n=X$.
Moreover, we assume the following two natural conditions:
There are correspondences $S_n\to S'_n$ for all $n$ which carries $R$ and the sets in the axioms (iii) and (iv) for $R$ to $R'$ and the sets for $R'$, respectively, and commutes with the operations in (i), (ii) and (v).
For $X\in S_n$, let $X^{R'}$ denote the corresponding element in $S'_n$.
The other condition is that for $X\in S_n$ and a $\{S_n\}$-definable $C^0$ function $f$ on $X$, $(\graph f)^{R'}$ is the graph of some $\{S'_n\}$-definable $C^0$ function, say, $f^{R'}$ on $X^{R'}$.\par
An example, which satisfies the conditions, is the semialgebraic structure.
For a semialgebraic set $X$ over $R$, we let $X^{R'}$ be the semialgebraic set over $R'$ defined by the same polynomial functions as for $X$.
Then the conditions are clearly satisfied, and, moreover, the correspondence $\cup_{n\in\N}S_n\ni X\to X^{R'}\in\cup_{n\in\N}S'_n$ is unique.
Another example is the case where $R$ is dense in $R'$.
We define $X^{R'}$ for a $\{S_n\}$-definable set $X$ as follows:
For simplicity of notation we assume that $X$ is bounded in $R^n$.
Then by Triangulation theorem of definable sets we have a stratification $\{X_i\}$ of $\overline X$ into $\{S_n\}$-definable $C^0$ manifolds such that $X$ is a union of some $X_i$ and the frontier condition is satisfied.
We set $X^{R'}=\cup_{X_i\subset X}({\rm c l}^{R'}(X_i)-\cup_{X_{i'}\subset\overline{X_i}-X_i}{\rm c l}^{R'}(X_{i'}))$, where the symbol ${\rm c l}^{R'}$ stands for the closure in $R^{\prime n}$.
Obviously, this case also satisfies the conditions because $X^{R'}$ does not depend on the choice of the stratification.
Moreover, the correspondence $X\to X^{R'}$ is unique because by the conditions, $(\overline{X_i})^{R'}\supset{\rm c l}^{R'}(X_i)\supset X^{R'}_i$, $(\overline{X_i})^{R'}$ is the disjoint union of $X^{R'}_i$ and $\cup_{X'_{i'}\subset\overline{X_i}-X_i}X^{R'}_{i'}$ and hence $X^{R'}_i={\rm c l}^{R'}(X_i)-\cup_{X_{i'}\subset\overline{X_i}-X_i}{\rm c l}^{R'}(X_{i'})$.
In particular, if $R=R'$ then $X^{R'}=X$.\par
As a corollary of Uniqueness theorem of definable triangulation and Theorems 2$'$ and 4, we obtain the following, which says that there is no topological difference between $X$ and $X^{R'}$:\vskip2mm\noindent
{\bf Corollary}. {\it (i) Any $\{S_n\}$-definable sets $X_1$ and $X_2$ are $\{S_n\}$-definably homeomorphic if and only if $X^{R'}_1$ and $X^{R'}_2$ are $\{S'_n\}$-definably homeomorphic.
In particularly, in the case $R=R'$, $X_1$ and $X_2$ are $\{S_n\}$-definably homeomorphic if and only if they are $\{S'_n\}$-definably homeomorphic.\par
(ii) A $\{S_n\}$-definable set $X$ is a $\{S_n\}$-definable $C^r$ manifold if and only if $X^{R'}$ is a $\{S'_n\}$-definable $C^r$ manifold, $0\le r\le\infty$.
The $X$ is a Nash manifold if and only if so is $X^{R'}$.\par
(iii) Any $\{S_n\}$-definable $C^r$ manifolds $X_1$ and $X_2$ are $\{S_n\}$-definably $C^r$ diffeomorphic if and only if $X^{R'}_1$ and $X^{R'}_2$ are $\{S'_n\}$-definably $C^r$ diffeomorphic, $0<r<\infty$.
If $R=R'$, $X_1$ and $X_2$ are $\{S_n\}$-definably $C^r$ diffeomorphic if and only if they are $\{S'_n\}$-definably $C^r$ diffeomorphic.}\vskip2mm
We assume Axiom (iii) to use the field structure of $\R$ and $R$.
It is natural to regard $\R$ and $R$ as linear spaces and replace (iii) by the following:\vskip1mm\noindent
(iii)$'$ {\it the sets $\{(x,y,z)\in\R^3:x+y=z\}$ and $\{(x,y)\in\R^2:x=c y\},\,c\in\R$, (or $\{(x,y,z)\in R^3:x+y=z\}$ and $\{(x,y)\in R^2:x=c y\},\,c\in R$) are elements of $S_3(X)$ and $S_2(X)$ (or $S_3$ and $S_2$) respectively.}\vskip1mm\noindent
In the forthcoming paper entitled, O-minimal Hauptvermutung III, we will extend Theorems 1, 2, 1$'$ and 2$'$ by using results in \cite{S3} and weakening the axiom (iii) to (iii)$'$.
There we need some fundamental argument on o-minimal structures which are partially explained in Chapter V in \cite{S2}.
\section{Definitions and facts quoted from \cite{S3}}
Let $r$ be any positive integer.
In this section we make a list of definitions and facts which are given in \cite{S3} and used in this paper.
As we already mentioned, properties of PL topology over $\R$ or semialgebraic sets are not necessarily satisfied for general $R$ or a general o-minimal structure.
We give this list to clarify the satisfied properties.
We can define most terminology of PL topology over $R$ also, e.g.~a simplex and a simplicial complex.
In this list we only give definitions of those terms which causes no confusion.
Hence, from now on, admitting the list only we regard $R$ as $\R$, and we do not need to distinguish $R$ from $\R$.\vskip1mm\noindent
{\it Definition}. A {\it full} subcomplex $L$ of a simplicial complex $K$ is a subcomplex such that each $\sigma\in K$ whose vertices are all in $L$ is a simplex in $L$.\vskip1mm\noindent
{\it Definition}. A {\it cell} $\sigma$ in $R^n$ is a ``\,bounded\," subset of $R^n$ of the form $\{x\in R^n:f_i(x)\dag_i0,\,i\in I\}$, where $\dag_i$ is now either = or $\ge$, $I$ is a finite set and $f_i$ are linear functions with coefficients in $R$, which is equivalent to say that a cell is a convex hull of a finite set of points by (2) in \cite{S3}, Sect.\ 2.
The interior and boundary of $\sigma$ are denoted by $\Int\sigma$ and $\partial\sigma$ respectively.
We naturally define a {\it cell complex}.\vskip1mm\noindent
{\it Definition}. For a union $X$ of cells in a cell complex $L$, set $L|_X=\{\sigma\in L:\sigma\subset X\}$.\vskip1mm\noindent
{\it Definition}. A {\it derived} subdivision of a cell complex $L$ is a simplicial subdivision $L'$ of $L$ defined by induction on dimension of cells such that for each $\sigma\in L$, the restriction $L'|_\sigma$ is the family of one point $v_\sigma$ in $\Int\sigma$, the simplices in $L'_{\partial\sigma}$ and of the cones with vertex $v_\sigma$ and bases simplices in $L'|_{\partial\sigma}$.\vskip1mm\noindent
{\it Definition}. Let $X\supset Y$ be compact polyhedra.
A {\it collar} of $Y$ in $X$ is a PL embedding $\phi:Y\times[0,\,1]\to X$ such that $\phi(\cdot,0)=\id$ and $\phi(Y\times[0,\,1))$ is an open neighborhood of $Y$ in $X$.\vskip1mm\noindent
{\it Definition}. A compact polyhedral neighborhood $U$ of $Y$ in $X$ is called a {\it regular} neighborhood of $Y$ in $X$ if there exist compact polyhedra $X_1\supset U_1\supset Y_1$, a PL homeomorphism $h:(X;U,Y)\to(X_1;U_1,Y_1)$ and a simplicial decomposition $K_1$ of $X_1$ such that $K_1|_{Y_1}$ is a full subcomplex of $K_1$ and $U_1=\cup_{x\in Y_1}|\st(x,K'_1)|$.
Here a PL homeomorphism $h:(X;U,Y)\to(X_1;U_1,Y_1)$ is defined to be a PL homeomorphism from $X$ to $X_1$ carrying $U$ and $Y$ to $U_1$ and $Y_1$, respectively, and the notations $K'_1,\ K''_1,\ \st(x,K'_1)$ and $|K_1|$ denote a derived subdivision of $K_1$, a second derived subdivision of $K_1$, i.e., a derived subdivision of $K'_1$, the star of $x$ in $K'_1$ and the underlying polyhedron of $K_1$ respectively.\vskip1mm\noindent
{\it Definition}. Let $x*X$ denote the {\it cone} with vertex $x$ and base $X$, i.e., $x*X=\{t x+(1-t)x':t\in[0,\,1],\,x'\in X\}$, provided $t'x+(1-t')x'\not=t''x+(1-t'')x''$ for distinct $(t',x')$ and $(t'',x'')$ in $[0,\,1)\times X$.\vskip1mm\noindent
{\bf Fact 2.1} ((8) in \cite{S3}, Sect.\,2). {\it A PL map $g:X\to Y$ is extended to a PL map $g^*:x*X\to y*Y$ so that $g^*(x)=y$ and $g^*$ is linear on the segment with ends $x$ and any point of $X$ if $x*X$ and $y*Y$ exist, and $g^*$ is a PL homeomorphism if $g$ is so} (the Alexander trick).\vskip1mm\noindent
{\bf Fact 2.2} ((9.1) in \cite{S3}, Sect.\,2). {\it Given a finite simplicial complex $K$, a full subcomplex $L$ and two derived subdivisions $K_1$ and $K_2$ of $K$ such that $K_1|_{|L|}=K_2|_{|L|}$, there is an isomorphism $\tau:K_1\to K_2$ such that $\tau=\id$ on $|L|\cup|K^0|$ and }
$$\tau(\cup_{x\in|L|}|\st(x,K_1)|)=\cup_{x\in|L|}|\st(x,K_2)|.$$\vskip1mm\noindent
{\bf Fact 2.3} ((9.2) in \cite{S3}, Sect.\,2). {\it Given $K$ and $L$ as above, a full subcomplex $M$ of $L$ and a subdivision $K_1$ of $K$, there are derived subdivisions $K'$ of $K$ and $K'_1$ of $K_1$ such that }
$$\cup_{x\in|L|-|M|}|\st(x,K')|=\cup_{x\in|L|-|M|}|\st(x,K'_1)|.$$\vskip1mm\noindent
{\bf Fact 2.4} ((9.4) in \cite{S3}, Sect.\,2). {\it For a compact PL manifold $X$ with boundary $\partial X$, a regular neighborhood of $\partial X$ in $X$ is a collar of $\partial X$ in $X$}.\vskip1mm\noindent
{\bf Fact 2.5} ((13) in \cite{S3}, Sect.\,2). {\it Let $\phi$ and $\psi$ be nonnegative PL functions on a compact polyhedron $X$ such that $\phi^{-1}(0)=\psi^{-1}(0)$, and let $K$ be a simplicial decomposition of $X$.
Then there exists a PL isotopy $\tau_t,\ 0\le t\le1$, of $X$ preserving $K$ such that $\tau_t=\id$ on the zero set of $\phi$ and $\phi\circ\tau_1=\psi$ on some neighborhood of the zero set of $\phi$ in $X$}.
Here a PL isotopy $\tau_t,\,0\le t\le1$, of $X$ is defined to be a PL homotopy from $X$ into itself such that $\tau_0=\id$ and $\tau_t$ is a homeomorphism of $X$ for each $t$, and the sentence that $\tau_t$ is preserving $K$ means that $\tau_t(\sigma)=\sigma$ for any $t$ and $\sigma\in K$.\vskip1mm\noindent
{\it Definition}. An {\it open simplex (cell)} is the interior of a simplex (cell), which is never used in PL topology.\vskip1mm\noindent
{\it Definition}. A {\it semilinear} set is a semialgebraic set defined by linear functions in place of polynomial functions.
A {\it semilinear} map is a map between semilinear sets whose graph is semilinear.\vskip1mm\noindent
{\it Definition}. Let $X$ and $Y$ be compact definable sets in $R^n$.
Two definable $C^0$ maps $f,\,g:X\to Y$ are {\it definably isotopic} or $f$ is {\it definably isotopic} to $g$ if there exists a definable $C^0$ map $F:X\times[0,\,1]\to Y$ such that $F(\cdot,0)=f(\cdot)$, $F(\cdot,1)=g(\cdot)$ and $F(\cdot,t)$ is an embedding for each $t\in[0,\,1]$.
We write $F(\cdot,t)$ as $f_t(\cdot)$, $0\le t\le 1,$ and call it a {\it definable isotopy} of $f_0$ to $f_1$.\vskip1mm\noindent
{\it Definition}. A {\it definable isotopy} of $X$ (to $f$) is a definable $C^0$ map $F:X\times[0,\,1]\to X$ such that $F(\cdot,0)=\id$ and $F(\cdot,t)$ is a homeomorphism of $X$ for each $t\in[0,\,1]$ ($F(\cdot,1)=f(\cdot)$).
Note that a definable isotopy $F:X\times[0,\,1]\to X$ of the identity map does not mean a definable isotopy of $X$.
In the former case, $F(\cdot,t)$ is not necessarily a homeomorphism of $X$; it is an embedding.\vskip1mm\noindent
{\it Definition}. In order to distinguish them we call a definable $C^0$ map $F:X\times[0,\,1]\to Y$ a {\it definable isotopy through homeomorphisms} if for each $t$, the map $F(\cdot,t)$ is a definable homeomorphism from $X$ to $Y$, and call definable $C^0$ maps $f,g:X\to Y${\it definably isotopic through homeomorphisms} if there exists a definable isotopy of $f$ to $g$ through homeomorphisms.
Note that a definable isotopy $f_t:X\to Y,\ 0\le t\le1$, is one through homeomorphisms if $Y$ is the underlying polyhedron of a finite cell complex $L$, $f_t$ is preserving $L$ and $f_0$ is a homeomorphism, which is shown by reductio ad absurdum as follows.:\par
Assume that this is false, which is weaker than that there is a simplex $\sigma$ over $R$ and a definable homotopy $h_t:\partial\sigma\to\partial\sigma,\ 0\le t\le1$, such that $h_0=\id$ and $h_1$ is a constant map.
Then the simplicial homotopy theorem (Lemma 3.1 in \cite{S3}) says that we can assume that the map $H:\partial\sigma\times[0,\,1]\ni(x,t)\to h_t(x)\in\partial\sigma$ is PL.
Moreover, moving the values of $H$ we can suppose that $H$ is a simplicial map between simplicial complexes whose vertices are real algebraic numbers.
Thus we reduce the problem to the case $R=\R_{\rm alg}$.
In this case there does not exist such $H$.\vskip1mm\noindent
{\it Definition}. We call a point $s$ of $R^n$ with $|s|=1$ a {\it singular direction} for a definable subset $X$ of $R^n$ if the set $X\cap(a+R s)$ has interior points in the line $a+R s$ for some point $a$ in $X$, where $a+R s$ denotes the set $\{a+b s:b\in R\}$.\vskip1mm\noindent
{\it Definition and fact}. The (definable) {\it Alexander trick} is the following extension of a map:
Let $X$ and $Y$ be compact definable sets in $R^n$, let $x$ and $y$ points in $R^n$ such that there exist the cones with vertices $x$ and $y$ and bases $X$ and $Y$, respectively, and let $g:X\to Y$ be a definable $C^0$ map.
Then $g$ is extended to a definable map $g^*:x*X\to y*Y$ so that $g^*(x)=y$ and $g^*$ is linear on each segment with ends $x$ and a point of $X$.
The $g^*$ is a homeomorphism if $g$ is a homeomorphism.\vskip1mm\noindent
{\it Fact and definition} (II.1.8 in \cite{S2}). {\it A definable set $E$ in $R^n$ is a finite disjoint union of definable $C^r$ manifolds.}
We call $\{E_1,...,E_k\}$ a {\it definable $C^r$ stratification} of $E$ and simply write it as $\{E_i\}_i$ or $\{E_i\}$.
We naturally define the dimension of a definable manifold, and call $\max_i\dim E_i$ the {dimension} of $E$.\vskip1mm\noindent
{\it Fact} (II.1.14 in \cite{S2}). For simplicity of notation and without loss of generality we always assume that a definable $C^r$ stratification $\{E_i\}_i$ satisfies the {\it frontier condition} (i.e., if $E_i\cap(\overline{E_{i'}}-E_{i'})\not=\emptyset$ then $E_i\subset\overline{E_{i'}}$) and that each $E_i$ is definably connected, unless otherwise specified.\vskip1mm\noindent
{\it Definition}. The stratification $\{E_i\}_i$ is {\it compatible} with a finite family of definable sets $\{X_j\}_j$ in $R^n$ if each $E\cap X_j$ is the union of some $E_i$'s.\vskip1mm\noindent
{\it Fact and definition} (II.1.17 in \cite{S2}). {\it Given a definable $C^0$ map $g:X\to Y$ between definable sets, there exist definable $C^r$ stratifications $\{X_i\}_i$ of $X$ and $\{Y_j\}_j$ of $Y$ such that for each $i$, the restriction $g|_{X_i}$ is a surjective $C^r$ submersion onto some $Y_j$.}
We write $g:\{X_i\}_i\to\{Y_j\}_j$ and call it a {\it definable} $C^r$ {\it stratification} of $g$.
{\it Moreover, given finite families of definable sets $\{A_\nu\}_\nu$ of $X$ and $\{A'_{\nu'}\}_{\nu'}$ of $Y$, we can choose a definable $C^r$ stratification $g:\{X_i\}_i\to\{Y_j\}_j$ so that $\{X_i\}_i$ and $\{Y_j\}_j$ are compatible with $\{A_\nu\}_\nu$ and $\{A'_{\nu'}\}_{\nu'}$ respectively.}
We say that $g:\{X_i\}_i\to\{Y_j\}_j$ is {\it compatible} with $(\{A_\nu\}_\nu,\{A'_{\nu'}\}_{\nu'})$.\vskip2mm\noindent
{\bf Triangulation theorem of definable sets} (p.\,186 in \cite{S3}). {\it Given a finite simplicial complex $K$ in $R^n$ and a finite number of compact definable sets $\{X_i\}_i$ in $|K|$, there exists a definable homeomorphism $\tau$ of $|K|$ preserving $K$ such that the $\{\tau^{-1}(X_i)\}_i$ are polyhedra, and, moreover, there exists a definable isotopy $\tau_t,\ 0\le t\le 1$, of $|K|$ preserving $K$ such that $\tau_1=\tau$.}\vskip1mm\noindent
{\it Definition}. We call $\tau:\{\tau^{-1}(X_i)\}_i\to\{X_i\}_i$ a {\it definable triangulation} of $\{X_i\}_i$.\vskip2mm\noindent
{\bf Triangulation theorem of definable $C^0$ functions} (Theorem 3.2 in \cite{S3}).\\
{\it $(1)$ Let $f$ be a definable continuous $R$-valued function defined on a compact polyhedron $X$ in $R^n$ and let $P$ be a cellular decomposition of $X$.
Then there exists a definable homeomorphism $\pi$ of $X$ preserving $P$ such that $f\circ\pi$ is PL.\\
$(2)$ Such a $\pi$ is the finishing homeomorphism of some definable isotopy of $X$ preserving $P$.\\
$(3)$ Moreover, $\pi$ is unique in the following sense:
Let $\pi'$ be a definable homeomorphism of $X$ preserving $P$ such that $f\circ\pi'$ is PL.
Then there exists a definable isotopy $\omega_t,\ 0\le t\le 1$, of $\pi^{-1}\circ\pi'$ preserving $P$ such that $\omega_1$ is PL and $f\circ\pi\circ\omega_t=f\circ\pi'$ for $t\in[0,\,1]$.}\vskip2mm\noindent
{\bf Uniqueness theorem of definable triangulation} (Theorem 1.1 in \cite{S3}). {\it If two families of compact polyhedra $(X_1;X_{1,1},\allowbreak...,X_{1,k})$ and $(X_2;X_{2,1},...,X_{2,k})$ are definably homeomorphic then they are PL homeomorphic.}\vskip2mm
Here the notation $(X_1;X_{1,1},...,X_{1,k})$ indicates that $X_{1,i}$ are subsets of $X_1$, we say that $(X_1;X_{1,1},...,\allowbreak X_{1,k})$ and $(X_2;X_{2,1},...,X_{2,k})$ are definably (or PL) homeomorphic if there is a definable (or PL) homeomorphism from $X_1$ to $X_2$ carrying each $X_{1,i}$ to $X_{2,i}$, and we write $(X_1;X_{1,1},...,X_{1,k})$ as $(X_1;X_{1,i})$.\vskip2mm\noindent
{\bf Supplement} (p.\,227 in \cite{S3}). {\it Any definable homeomorphism $(X_1;X_{1,i})\to(X_2;X_{2,i})$ is definably isotopic to a PL homeomorphism through homeomorphisms.}\vskip2mm
In \cite{S3} the condition that the isotopy is one through homeomorphisms is missing.
However, all the arguments in \cite{S3} proceed under this condition.\par
Let us follow the proof of Triangulation theorem of definable $C^0$ functions (1) in \cite{S3}.
Let $p:R^n\times R\to R$ and $q_1:R^n\times R\to R^n$ be the projections and let $R^n\times R\stackrel{p_1}{\longrightarrow}R^{n-1}\times R\longrightarrow\cdots\stackrel{p_n}{\longrightarrow}R$ be the projections ignoring the respective first factors.
Set $A=\graph f$ and $A_t=\{x\in R^n:(x,t)\in A\}$ for each $t\in R$.
For simplicity of notation, we assume that $A$ is a compact definable subset of $R^n\times R$ of local dimension $n$ or $n+1$ everywhere.
There exists a definable $C^r$ stratification $\{A_j\}_j$ of $A$ compatible with $\{\sigma\times R:\sigma\in P\}$ such that the $\{p_1(A_j)\}_j,...,\{p_n\circ\cdots\circ p_1(A_j)\}_j$ are definable $C^r$ stratifications of $p_1(A),...,p_n\circ\cdots\circ p_1(A)$, respectively, and the $p_1|_A:\{A_j\}_j\to\{p_1(A_j)\}_j,...,\,p_n|_{p_{n-1}\circ\cdots\circ p_1(A)}:\{p_{n-1}\circ\cdots\circ p_1(A_j)\}_j\to\{p_n\circ\cdots\circ p_1(A_j)\}_j$ are definable $C^r$ stratifications of $p_1|_A:A\to p_1(A),...,\,p_n|_{p_{n-1}\circ\cdots\circ p_1(A)}:p_{n-1}\circ\cdots\circ p_1(A)\to p_n\circ\cdots\circ p_1(A)$ respectively.
We simply write $\{A_j\}_j\stackrel{p_1}{\longrightarrow}\{p_1(A_j)\}_j\stackrel{p_2}{\longrightarrow}\cdots\stackrel{p_n}{\longrightarrow}\{p_n\circ\cdots\circ p_1(A_j)\}_j$ and $A\stackrel{p_1}{\longrightarrow}p_1(A)\stackrel{p_2}{\longrightarrow}\cdots\stackrel{p_n}{\longrightarrow}p_n\circ\cdots\circ p_1(A)$.) 
Then, since $\dim p^{-1}_k(a)=1\ (k=1,...,n)$ for each $a\in p_k\circ\cdots\circ p_1(A)$, each $p_{k-1}\circ\cdots\circ p_1(A_j)$ is described by definable $C^r$ functions on $p_k\circ\cdots\circ p_1(A_j)$ as follows:
There exist definable $C^r$ functions $\phi_{k,j}$ and $\psi_{k,j}$ on $p_k\circ\cdots\circ p_1(A_j)$ such that $p_{k-1}\circ\cdots\circ p_1(A_j)$ is of the form
\begin{align*}\{(x_k,x',t)\in R\times(p_k\circ\cdots\circ p_1(A_j)):\phi_{k,j}(x',t)<x_k<\psi_{k,j}(x',t)\}\\
\text{or}\quad\{(x_k,x',t)\in R\times(p_k\circ\cdots\circ p_1(A_j)):x_k=\phi_{k,j}(x',t)\},\end{align*}
where $x'\in R^{n-k}$.
Then we need two conditions.
One is the following:\vskip1mm\noindent
{\bf Condition 2.6} ((**) on p.\,197 in \cite{S3}). {\it Each of $\phi_{k,j}$ and $\psi_{k,j}$ is extended to a definable $C^0$ function on $p_k\circ\cdots\circ p_1(\overline{A_j})$.}\par
Let $R^n\stackrel{r_1}{\longrightarrow}\cdots\stackrel{r_{n-1}}{\longrightarrow}R$ denote the projections ignoring the respective first factors, and let $\{A_{k,j}\}_j$ and $P_k$ denote $\{p_{k-1}\circ\cdots\circ p_1(A_j)\}_j$ and $\{r_{k-1}\circ\cdots\circ r_1(\sigma):\sigma\in P\}$ respectively.
The the other condition is the following:\vskip1mm\noindent
{\bf Condition 2.7} ((***) on p.\,201 in \cite{S3}). {\it Given $\sigma\in P_k$ and $A_j$ such that $p_k\circ\cdots\circ p_1(A_j)\subset\Int\sigma\times R$ and $\dim p_k\circ\cdots\circ p_1(A_j)\le\dim\sigma$, the restrictions $q_k|_{p_k\circ\cdots\circ p_1(A_j)}$ and $q_k|_{\overline{p_k\circ\cdots\circ p_1(A_j)}}$ are a $C^r$ embedding into $R^{n+1-k}$ and an injection, respectively}.\vskip2mm\noindent
{\bf Lemma 2.8} (Lemma 3.4 in \cite{S3}). {\it Assume that condition 2.6 is already satisfied for any small linear perturbation of $p_1,...,p_{n-1}$.
Then we can choose such $p_1,...,p_{n-1}$ and $\{A_j\}_j$ so that conditions 2.6 and 2.7 are satisfied.}\vskip2mm
The actual statement of Lemma 3.4 in \cite{S3} is a little different from this.
We translate it so that we can apply it in the forthcoming argument.
As the proof becomes easier by the assumption in the above lemma, we do not prove it.\vskip2mm\noindent
{\bf Proposition 2.9} (Proposition 3.6 in \cite{S3}). {\it Triangulation theorem of definable $C^0$ functions (1) holds under conditions 2.6 and 2.7.}\vskip2mm
We will use not only this proposition but also its proof.
In the proof in \cite{S3} we explicitly constructed a definable homeomorphism $\pi$ such that $f\circ\pi$ is PL.
We need that method of construction in the proof of Lemma 7 below on triangulations of definable $C^0$ maps into $R^2$.
\section{Theorems 1 and 2 in the compact case}
{\it Proof of Theorem 2 for a compact definable $C^0$ manifold}. Let $M$ be a compact definable $C^0$ manifold over $R$.
By Triangulation theorem of definable sets, $M$ is definably homeomorphic to a compact polyhedron.
Hence we can assume that $M$ is a compact polyhedron.
Moreover, by Uniqueness theorem of definable triangulation, if $M$ is definably homeomorphic to another compact polyhedron then they are PL homeomorphic.
Therefore, it remains to prove that $M$ is a PL manifold, i.e., the following statement:\par
{\it For each $x\in M$ there exists a polyhedral neighborhood $U_0$ of $x$ in $M$ such that $(U_0;x)$ is PL homeomorphic to $(\sigma_0;0)$, where $\sigma_0$ is a simplex such that $0\in\Int\sigma_0$.}\par
By definition of a definable $C^0$ manifold there exists a compact definable neighborhood $U$ of $x$ in $M$ such that $(U;x)$ is definably homeomorphic to $(\sigma_0;0)$.
Apply Triangulation theorem of definable sets to $U$ and a simplicial decomposition of $M$ such that $x$ is a vertex.
Then there is a definable homeomorphism $\tau$ of $M$ such that $\tau(x)=x$ and $\tau(U)$ is a polyhedron.
Since $\tau(U)$ is a compact polyhedral neighborhood of $x$ in $M$, replacing $U$ with $\tau(U)$ we assume that $U$ is a compact polyhedral neighborhood from the beginning.
Then by Uniqueness theorem of definable triangulation, $(U;x)$ is PL homeomorphic to $(\sigma_0;0)$.\qed\vskip2mm\noindent
{\it Proof of Theorem 1 for a compact tame $C^0$ manifold}. Let $M$ be a compact tame $C^0$ manifold over $\R$.
Let $S$ denote the o-minimal structure $\{S_n(M)\}_n$ over $\R$.
Then by Triangulation theorem of definable sets, $M$ is $S$-definably homeomorphic to a compact polyhedron $X$.
Here $M$ and $X$ are definably homeomorphic in any o-minimal structure containing $S$.
Thus, if $M$ is tamely homeomorphic to some compact polyhedron $X'$ then $M$ is definably homeomorphic to $X'$ in some o-minimal structure $S'$, and $M$ and $X$ are $S'$-definably homeomorphic since $S'\supset S$.
Therefore, $X$ and $X'$ are PL homeomorphic by Uniqueness theorem of definable triangulation.
Namely, $X$ is uniquely determined up to PL homeomorphisms.\par
It remains to see that $X$ is a PL manifold.
By the above argument it suffices to show that each point $x$ of $M$ has a neighborhood $U$ such that $(U;x)$ is tamely homeomorphic to $(\sigma_0;0)$, where $\sigma_0$ is a simplex such that $0\in\Int\sigma$.
This is obvious by the definition of a tame $C^0$ manifold.\qed
\section{Theorems 1 and 2 in the noncompact case}
Let $M\subset R^n$ be a noncompact definable $C^0$ manifold in Theorem 2.
Here we can assume that $M$ is bounded in $R^n$ since $R^n$ is semialgebraically homeomorphic to $\{x\in R^n:|x|<1\}$.
Apply Triangulation theorem of definable sets to $(\overline M;\overline M-M)$.
Then $M$ is definably homeomorphic to a definable polyhedron $X$ such that $\overline X$ is a compact polyhedron.
Moreover, we see that $X$ is a PL manifold by Uniqueness theorem of definable triangulation, which we call informally a {\it triangulation} of $M$.
Hence we need to find a triangulation $X$ such that $\overline X$ is a compact PL manifold with boundary.\par
For this, let us introduce unique triangulations of general noncompact definable sets.
Indeed, uniqueness of triangulations of general compact definable sets was the key of the proof in the compact case, although all the arguments in this paper become simple if the sets are locally closed in their ambient spaces, i.e., the sets are open in their closures.
The notion of semilinear is convenient for this.
However, the category of all semilinear sets and semilinear $C^0$ maps is too large.
There are two typical examples of semilinear sets $X_1$ and $X_2$ which are definably homeomorphic but not semilinearly homeomorphic.
First, $X_1=R$ and $X_2=(0,\,1)$.
Indeed, any semilinear function $f$ on $R$ is constant on $(c,\,\infty)$ for some $c\in R$ or $f(x)$ converges to an infinity as $x\to\infty$.
Secondly, $X_1=\sigma-\{0\}$ and $X_2=\sigma-\sigma/2$, where $\sigma$ is a 2-simplex such that $0\in\Int\sigma$.
We want to treat a natural family of semilinear sets where two semilinear sets are semilinearly homeomorphic if they are definably homeomorphic.
We avoid $X_1$ of the examples and consider only the next family which contains $X_2$ and any compact polyhedron.
We also avoid a semilinear set $X_1$ and choose $X_2$ if $X_1$ and $X_2$ are definably homeomorphic, $\overline{X_1}$ is not a PL manifold with boundary and $\overline{X_2}$ is a PL manifold with boundary, e.g.~$X_1=\partial\sigma-\{a\}$ and $X_2=(0,\,1)$ for a 2-simplex $\sigma$ and $a\in\partial\sigma$.
We precisely define the family below.\par
We call a semilinear set $X$ in $R^n$ {\it standard} if $X$ is bounded in $R^n$ and there is a {\it collar} of $\overline X-X$ in $\overline X$, i.e., there exists a semilinear $C^0$ embedding $\phi:(\overline X-X)\times[0,\,1]\to\overline X$ such that $\phi(\cdot,0)=\id$ and $\phi((\overline X-X)\times[0,\,1))$ is an open neighborhood of $\overline X-X$ in $\overline X$ (Fig.\ 1 (a)).
We also call the image of $\phi$ a {\it collar} when we know $\phi$ in the context.
(In PL topology a collar is defined on a closed subpolyhedron.) 
Note that $\phi((\overline X-X)\times(0,\,1])\subset X$, the image of $\phi$ is not a triangle in Fig.\ 1 (a) and a noncompact PL manifold $X$ in $R^n$ is a standard semilinear set if and only if $\overline X$ is a compact PL manifold with boundary.
The natural family is defined to be that of standard semilinear sets.
From now on, for a cell complex $K$, let $K'$ and $K''$ always denote the barycentric and second barycentric subdivisions of $K$, respectively, in this paper (Fig.\ 1 (b)).
The next remark shows how standard semilinear sets are unified.\vskip2mm
\begin{figure}[h]
\begin{pspicture}(-2.7,0.3)(12.7,2.6)
\psline(0.86,0.98)(2.2,0.5)
\psline(2.2,0.5)(2.2,2.5)
\psline(2.2,2.5)(0.8,2.5)
\psline(0.8,2.5)(0.8,2)
\psline[linestyle=dashed, dash=2pt 1pt](0.8,2)(0.8,1.07)
\pscircle(0.8,1){0.07}
\psdot(0.8,2)
\rput(1.5,1.6){$X$}
\psline(2.95,1.93)(2.95,1)
\pscircle(2.95,2){0.07}
\psdot(2.95,1)
\rput(2.95,2.4){$\overline X\!-\!X$}
\psline(3.6,1)(4.3,0.8)
\psline(4.3,0.8)(4.3,1.93)
\psline[linestyle=dashed, dash=2pt 1pt](4.23,2)(3.67,2)
\psline(3.6,1.93)(3.6,1)
\pscircle(3.6,2){0.07}
\pscircle(4.3,2){0.07}
\rput(2.5,0.2){(a)}
\rput(3.9,0.5){$\Ima\phi$}
\pscurve{->}(3.8,0.65)(3.9,1.1)(4,1.4)
\pspolygon(5.5,0.8)(6.5,0.8)(6,2)
\pspolygon(6.8,0.8)(7.8,0.8)(7.3,2)
\psline(6.8,0.8)(7.3,1.2)
\psline(7.3,0.8)(7.3,1.2)
\psline(7.8,0.8)(7.3,1.2)
\psline(7.55,1.4)(7.3,1.2)
\psline(7.3,2)(7.3,1.2)
\psline(7.05,1.4)(7.3,1.2)
\pspolygon(8.1,0.8)(9.1,0.8)(8.6,2)
\psline(8.1,0.8)(8.6,1.2)
\psline(8.6,0.8)(8.6,1.2)
\psline(8.43,0.93)(8.1,0.8)
\psline(8.43,0.93)(8.35,0.8)
\psline(8.43,0.93)(8.35,1)
\psline(8.43,0.93)(8.6,1.2)
\psline(8.43,0.93)(8.6,1)
\psline(8.43,0.93)(8.6,0.8)
\psline(8.43,0.93)(8.6,1)
\psline(9.1,0.8)(8.6,1.2)
\psline(8.77,0.93)(8.6,0.8)
\psline(8.77,0.93)(8.6,1.2)
\psline(8.77,0.93)(9.1,0.8)
\psline(8.77,0.93)(8.6,1)
\psline(8.77,0.93)(8.85,1)
\psline(8.77,0.93)(8.85,0.8)
\psline(8.85,1.4)(8.6,1.2)
\psline(8.85,1.13)(9.1,0.8)
\psline(8.85,1.13)(8.6,1.2)
\psline(8.85,1.13)(8.85,1.4)
\psline(8.85,1.13)(8.85,1)
\psline(8.85,1.13)(8.725,1.3)
\psline(8.85,1.13)(8.975,1.1)
\psline(8.6,2)(8.6,1.2)
\psline(8.68,1.53)(8.85,1.4)
\psline(8.68,1.53)(8.6,1.2)
\psline(8.68,1.53)(8.6,2)
\psline(8.68,1.53)(8.725,1.3)
\psline(8.68,1.53)(8.6,1.6)
\psline(8.68,1.53)(8.725,1.7)
\psline(8.35,1.4)(8.6,1.2)
\psline(8.53,1.53)(8.6,2)
\psline(8.53,1.53)(8.6,1.2)
\psline(8.53,1.53)(8.35,1.4)
\psline(8.53,1.53)(8.6,1.6)
\psline(8.53,1.53)(8.475,1.3)
\psline(8.53,1.53)(8.475,1.7)
\psline(8.35,1.13)(8.35,1.4)
\psline(8.35,1.13)(8.6,1.2)
\psline(8.35,1.13)(8.1,0.8)
\psline(8.35,1.13)(8.475,1.3)
\psline(8.35,1.13)(8.35,1)
\psline(8.35,1.13)(8.225,1.1)
\rput(6,2.4){$K$}
\rput(7.3,2.4){$K'$}
\rput(8.6,2.4){$K''$}
\rput(7.3,0.2){(b)}
\end{pspicture}
\caption{Collar and subdivisions}
\end{figure}
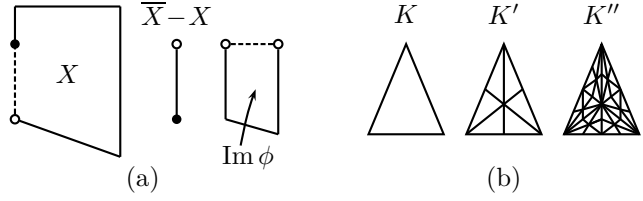\noindent
{\bf Remark 3.}
{\it Let $X$ be a standard semilinear set in $R^n$.\par\noindent
{\rm (i)} A semilinear $C^0$ function on $X$ is extendable to a PL function on $\overline X$.\par\noindent
{\rm (ii)} Set $X(0)=X,\ X(j+1)=\overline{X(j)}-X(j)\ (j=0,1,...).$
Then $\overline{X(j)}$ is the disjoint union of $X(j)$ and $X(j+1)$, $X(j)\supset X(j+2)$, $X(j)-X(j+2)$ is the set of inner points of $X(j)$ in $\overline{X(j)}$, $X(j)$ is not necessarily standard (e.g.\ $X$ is the union of an open 2-simplex and a point of its boundary), and $X(2j)$ is standard.\par\noindent
{\rm (iii)} It is known that there are two compact PL manifolds with boundary which are not PL homeomorphic but whose interiors are PL homeomorphic \cite{M}.
However, two compact PL manifolds with boundary are PL homeomorphic if their interiors are semilinearly homeomorphic.
Moreover, a semilinear homeomorphism $f:X\to Y$ between standard semilinear sets is extended to a PL homeomorphism $\overline f:\overline X\to\overline Y$.
The last statement is not the case unless $X$ and $Y$ are standard.
A counter-example is $X=[0,\,1)\cup(1,\,2]$ and $Y=[0,\,1)\cup(2,\,3]$.}\vskip2mm\noindent
{\it Proof of Remark 3}. (i) Let $g$ be a semilinear $C^0$ function on $X$ and let $\phi:(\overline X-X)\times[0,\,1]\to\overline X$ be a collar of $\overline X-X$ in $\overline X$.
Then it suffices to see that the semilinear $C^0$ function $g\circ\phi|_{(\overline X-X)\times(0,\,1]}$ is extendable to a PL function on $(\overline{\overline X-X})\times[0,\,1]$.
Let $K$ be a simplicial decomposition of $(\overline{\overline X-X})\times[0,\,1]$ such that $(\overline X-X)\times(0,\,1]$ is the union of some open simplices in $K$ and $g\circ\phi$ is linear on each open simplex in $K$ contained in $(\overline X-X)\times(0,\,1]$.
This is possible because a semilinear set is a finite union of open simplices.
For such an open simplex $\Int\sigma$, $g\circ\phi|_{\Int\sigma}$ is extended to a linear function on $\sigma$.
We need to show uniqueness of these extensions.
For this it suffices to see the following statement:\vskip1mm\noindent
{\it Let $x_0\in\overline X-X$ and let $l_1$ and $l_2$ be distinct segments in $(\overline X-X)\times[0,\,1]$ with one end $(x_0,0)$ such that $l_1-\{(x_0,0)\}$ and $l_2-\{(x_0,0)\}$ are contained in $(\overline X-X)\times(0,\,1]$ and $l_1$ and $l_2$ are contained in some simplices in $K$ respectively.
Then $\lim_{l_1\ni(x,t)\to(x_0,0)}g\circ\phi(x,t)=\lim_{l_2\ni(x,t)\to(x_0,0)}g\circ\phi(x,t)$.}\par
This is easy to prove.
Let $\sigma_i\ (i=1,2)$ denote the 2-simplices with two 1-faces $l_i$ and $\{x_0\}\times[0,\,1]$ if $l_i\not\subset\{x_0\}\times[0,\,1]$, and set $\sigma_i=l_i$ otherwise.
Then $\sigma_i-\{(x_0,0)\}\subset(\overline X-X)\times(0,\,1]$ and we can reduce the problem to the case where $\sigma_1$ and $\sigma_2$ are contained in some simplices in $K$.
Hence we see as above that $g\circ\phi|_{\sigma_i-\{(x_0,0)\}}$ are extended to PL functions on $\sigma_i$.
Thus
$$\lim_{l_1\ni(x,t)\to(x_0,0)}g\circ\phi(x,t)=\lim_{t\to 0}g\circ\phi(x_0,t)=\lim_{l_2\ni(x,t)\to(x_0,0)}g\circ\phi(x,t).$$\par
(ii) Only the last statement is not obvious.
By induction it suffices to show that $X(2)$ is standard.
Let $\phi$ be a collar as above.
We will see that $\phi|_{X(3)\times[0,\,t_0]}$ is a collar of $X(3)$ in $\overline{X(2)}$ for some $t_0\in(0,\,1]$, i.e., $\phi^{-1}(\overline{X(2)})\cap(X(1)\times[0,\,t_0])=X(3)\times[0,\,t_0]$.
The last equality is equivalent to the following statement:\vskip1mm\noindent
{\it Let $x\in X(1)$.
Then $\phi(\{x\}\times(0,\,t_0])\subset X-X(2)$ if and only if $x\in X(1)-X(3)$.}\par
 Note that the image of a compact polyhedron under a semilinear $C^0$ map is compact.
Assume that $x\in X(1)-X(3)$.
Then $x$ is an inner point of $X(1)$ in $\overline{X(1)}$, $X(1)\cap U$ is compact for some compact polyhedral neighborhood $U$ of $x$ in $R^n$, $\phi((X(1)\cap U)\times[0,\,1])$ is a compact neighborhood of $\phi(x,t)$ in $X$ for each $t\in(0,\,1)$ and hence $\phi(x,t)\not\in X(2)$.
Conversely assume that $\phi(x,t)\in X-X(2)$ for some $t\in(0,\,1)$.
Then there is a closed polyhedral neighborhood $U$ of $(x,t)$ in $X(1)\times[0,\,1]$ such that $\phi(U)$ is compact.
Such a $U$ is compact because $\phi^{-1}|_{\phi(U)}$ is a semilinear $C^0$ map.
Hence $x\in X(1)-X(3)$.
Thus the statement is proved.\par
(iii) We prove that $f:X\to Y$ in (iii) is extended to PL $\overline f:\overline X\to\overline Y$.
By (i), $f$ and $f^{-1}$ are extended to PL maps $\overline f:\overline X\to\overline Y$ and $\overline{f^{-1}}:\overline Y\to\overline X$, respectively, and $\overline{f^{-1}}\circ\overline f:\overline X\to\overline X$ and $\overline f\circ\overline{f^{-1}}:\overline Y\to\overline Y$ are the $C^0$ extensions of the identity maps $f^{-1}\circ f:X\to X$ and $f\circ f^{-1}:Y\to Y$ respectively.
Hence $\overline f:\overline X\to\overline Y$ is a PL homeomorphism.\qed\vskip2mm
We need to consider multiple semilinear sets at once.
We call a family of semilinear sets $(X;X_i)$ {\it standard} if $\{X_i\}$ is finite, $X$ is standard and we can choose the semilinear $C^0$ embedding $\phi:(\overline X-X)\times[0,\,1]\to\overline X$ so that $\phi^{-1}(X_i)=X'_i\times(0,\,1]$ for each $i$ and some semilinear subset $X'_i$ of $\overline X-X$.
We call $\phi$ or $\Ima\phi$ a {\it collar} of $\overline X-X$ in $(\overline X;X,X_i)$.
Remark 3, (ii) is generalized so that $(X(2j);X(2j)\cap X_i)$ is standard for each $j$.
Note that for noncompact, semilinear and bounded $X$ and a simplicial decomposition $K$ of $\overline X$, $(X;\sigma\cap X:\sigma\in K)$ is never standard if $X$ is of dimension $>1$ locally at $\overline X-X$.
This means that the usual argument of PL topology is not sufficient for our problem.
Now we state uniqueness triangulations of noncompact definable sets.\vskip2mm\noindent
{\bf Theorem 4} (Unique standard triangulation theorem of definable sets).
{\it A finite family of definable sets $(X;X_i)$ is definably homeomorphic to a unique standard family of semilinear sets when $X$ is not compact.}\vskip2mm
Theorems 1 and 2 in the noncompact case follow quickly from Theorem 4 as in the compact case.\vskip2mm\noindent
{\it Proof of Theorem 2 for a noncompact definable $C^0$ manifold under Theorem 4}. Let $M$ be a noncompact definable $C^0$ manifold over $R$.
By Theorem 4, $M$ is definably homeomorphic to a standard semilinear set $X$.
Then $X$ is a PL manifold for the same reason as in the proof of Theorem 2 in the compact case, and $\overline X$ is a compact PL manifold with boundary as we already noted.
It remains to see uniqueness of $\overline X$.
Let $M$ be definably homeomorphic to the interior of another compact PL manifold $M_1$ with boundary.
Then $\Int M_1$ is standard.
Hence by uniqueness in Theorem 4, $\overline X$ and $M_1$ are PL homeomorphic.\qed\vskip2mm\noindent
{\it Proof of Theorem 1 for a noncompact tame $C^0$ manifold under Theorem 4}. We proceed as in the compact case.
Let $M$ be a noncompact tame $C^0$ manifold over $\R$.
Let $S$ denote the o-minimal structure $\{S_n(M)\}$.
Then by Theorem 4, $M$ is $S$-definably homeomorphic to a standard semilinear set $X$, and $M$ and $X$ are definably homeomorphic in any o-minimal structure containing $S$.
Thus, if $M$ is tamely homeomorphic to some standard semilinear set $X'$ then $M$ is definably homeomorphic to $X'$ in some o-minimal structure $S'$, and $M$ and $X$ are $S'$-definably homeomorphic.
Therefore, $X$ and $X'$ are semilinearly homeomorphic by Theorem 4.
It follows that $\overline X$ and $\overline{X'}$ are PL homeomorphic.
Hence it remains to see that $\overline X$ is a PL manifold with boundary, i.e., $X$ is a PL manifold.
However, we have already shown this in the proof of the compact case.\qed
\section{Lemmas for Theorem 4}
In this and the next sections we are engaged in proving Theorem 4.
We will use the following lemmas in its proof:\vskip2mm\noindent
{\bf Lemma 5.} {\it Let $K$ be a cell complex generated by one cell $\sigma_0$, let $\sigma_1$ be a (possibly empty) proper face of $\sigma_0$, and let $\sigma_2,...,\sigma_k$ be some proper faces of $\sigma_1$.
Set $\delta=\cup_{j=1}^k\Int\sigma_j$.
Then a definable isotopy of $\partial\sigma_0-\delta$ preserving $\{\sigma-\delta:\sigma\in K|_{\partial\sigma_0}\}$ is extendable to a definable isotopy of $\sigma_0-\delta$.}\vskip2mm\noindent
{\it Proof of Lemma 5}.
We proceed by induction on the dimension of $\delta$.
Let $\partial\alpha_t,\ 0\le t\le 1$, denote the definable isotopy of $\partial\sigma_0-\delta$.
Assume that $0\in\Int\sigma_0$.
Then it is natural to define an extension $\alpha_t$ by $\alpha_t(sy)=s\partial\alpha_t(y)$ for $(y,s,t)\in(\partial\sigma_0-\delta)\times[0,\,1]^2$.
However, such an $\alpha_t$ is not defined on the interior of $\sigma_0$ because $\{s y:s\in[0,\,1],y\in\partial\sigma_0-\delta\}\not=\sigma_0$.
We need to modify the definition.\par
First we reduce the problem to the case where $\partial\alpha_t=\id$ on $\partial\sigma_1-\delta$ (not $\partial\sigma_0-\delta$).
There is a cellular subdivision $K_1$ of $K$ such that $K|_{\partial\sigma_1}=K_1|_{\partial\sigma_1}$ and for each $\sigma$ in $K_1$ not contained in $\partial\sigma_1$, $(\sigma;\sigma\cap\sigma_j:j=2,...,k)$ satisfies the conditions on $(\sigma_0;\sigma_j:j=1,...,k)$ in Lemma 5.
Indeed, choose one point in each open cell in $K$ not contained in $\partial\sigma_1$, and define a cell complex $K_1$ such that $K|_{\partial\sigma_1}=K_1|_{\partial\sigma_1}$ and $K_1^0$ consists of the points and $K^0\cap\partial\sigma_1$ in the same way as we have defined a derived subdivision of a cell complex.
Then for each $\sigma$ in $K_1$ not contained in $\partial\sigma_1$, $(\sigma;\sigma\cap\sigma_j:j=2,...,k)$ satisfies the conditions on $(\sigma_0;\sigma_j:j=1,...,k)$.
Let us consider $\partial\alpha_t|_{\partial\sigma_1-\delta}$ and $K_1$.
Then by using the induction hypothesis, we can extend $\partial\alpha_t|_{\partial\sigma_1-\delta}$ to a definable isotopy $\beta_t,\ 0\le t\le1$, of $\sigma_0-(\delta-\Int\sigma_1)$ preserving $\{\sigma-(\delta-\Int\sigma_1):\sigma\in K_1\}$.
Set $\partial\gamma_t=\beta_t^{-1}\circ\partial\alpha_t,\ 0\le t\le1$, which is a definable isotopy of $\partial\sigma_0-\delta$ preserving $\{\sigma-\delta:\sigma\in K|_{\partial\sigma_0}\}$ and whose restriction to $\partial\sigma_1-\delta$ is the identity map.
If $\partial\gamma_t,\ 0\le t\le1$, is extended to a definable isotopy $\gamma_t,\ 0\le t\le1$, of $\sigma_0-\delta$, then $\beta_t\circ\gamma_t,\ 0\le t\le1$, is a definable isotopy of $\sigma_0-\delta$ and is an extension of $\partial\alpha_t,\ 0\le t\le1$.
Hence we can assume that $\partial\alpha_t=\id$ on $\partial\sigma_1-\delta$ from the beginning.\par
Let $\xi:\partial\sigma_0\times[0,\,1]\to\sigma_0$ be a semialgebraic $C^0$ map such that $\xi(x,s)=x$ for $(x,s)\in(\partial\sigma_0\times\{0\})\cup(\sigma_1\times[0,\,1])$ and $\xi|_{(\partial\sigma_0-\sigma_1)\times[0,\,1]}$ is a $C^0$ embedding.
For example, let $\xi$ be the composite of a map $\partial\sigma_0\times[0,\,1]\ni(x,s)\to(x,s\rho(x))\in\partial\sigma_0\times[0,\,1]$ and a collar of $\partial\sigma_0$ in $\sigma_0$, where $\rho:\partial\sigma_0\to[0,\,1]$ is a PL map with zero set $\sigma_1$.
Let us define isotopies $\tilde\alpha_t$, $0\le t\le1$, of $(\partial\sigma_0-\delta)\times[0,\,1]$ to be $\tilde\alpha_t(x,s)=(\partial\alpha_{t(1-s)}(x),s)$, and $\alpha_t,\ 0\le t\le1$, of $\Ima\xi-\delta$ by 
$$\alpha_t(x)=\xi\circ\tilde\alpha_t\circ\xi^{-1}(x)\quad\text{for}\ x\in\Ima\xi.$$
Then $\tilde\alpha_t(x,0)=(\partial\alpha_t(x),0)$, $\tilde\alpha_t=\id$ on $(\partial\sigma_0-\delta)\times\{1\}$, $\alpha_t,\ 0\le t\le1$, is a well-defined definable isotopy of $\Ima\xi-\delta$ since $\partial\alpha_t=\id$ on $\partial\sigma_1-\delta$, and $\alpha_t=\id$ on $\xi(\partial\sigma_0\times\{1\})-\delta$---the boundary of $\Ima\xi-\delta$ in $\sigma_0-\delta$.
Hence we can extend $\alpha_t,\ 0\le t\le1$, to a definable isotopy of $\sigma_0-\delta$ by setting $\alpha_t=\id$ outside of the image of $\xi$.
Thus Lemma 5 is proved.\qed\vskip2mm\noindent
{\bf Lemma 6} (Lifts of definable homeomorphisms).
{\it Let $\gamma:L\to M$ be a cellular map between finite cell complexes and $D$ a subset of $|M|$.
Let $\delta_M$ be a definable homeomorphism of $|M|$ preserving $M$ such that $\delta_M=\id$ on $D$.
Then there exists a definable homeomorphism $\delta_L$ of $|L|$ preserving $L$ such that $\delta_L=\id$ on $\gamma^{-1}(D)$ and $\delta_M\circ\gamma=\gamma\circ\delta_L$.}\vskip2mm\noindent
{\it Proof of Lemma 6}.
For each $\sigma\in L$ we will construct a definable homeomorphism $\delta_\sigma$ of $\sigma$ such that $\delta_M\circ\gamma=\gamma\circ\delta_\sigma$ on $\sigma$ and $\delta_\sigma=\delta_{\sigma'}$ on $\sigma\cap\sigma'$ for other $\sigma'\in L$.
If $\gamma|_\sigma$ is injective, then $\delta_\sigma$ is uniquely determined.
In the other case, we define $\delta_\sigma$ by a cross-section as follows:
Let $c_\sigma:\gamma(\sigma)\to\sigma$ be a PL cross-section of $\gamma|_\sigma$ such that $c_\sigma(x)\in\Int(\gamma|_\sigma)^{-1}(x)$ for each $x\in\gamma(\sigma)$ and $c_\sigma=c_{\sigma'}$ on $\gamma(\sigma')$ for each face $\sigma'$ of $\sigma$ with $\sigma'=(\gamma|_\sigma)^{-1}(\gamma(\sigma'))$, which we can construct by double induction on $\dim\sigma-\dim\gamma(\sigma)$ and $\dim\sigma$ by the Alexander trick (Fact 2.1).
Note $(\gamma|_\sigma)^{-1}(x)$ is a cell of the form $c_\sigma(x)*\partial((\gamma|_\sigma)^{-1}(x))$.
By induction on the dimension of $\sigma$ we assume that $\delta_\sigma$ is already constructed on the boundary of $\sigma$.
Let $\delta_{\partial\sigma}$ denote the homeomorphism of the boundary of $\sigma$.
We extend naturally the restriction of $\delta_{\partial\sigma}$ to $\cup_{x\in\gamma(\sigma)}\partial((\gamma|_\sigma)^{-1}(x))$ (not $\delta_{\partial\sigma}$ itself) to $\delta_\sigma$ by 
$$\delta_\sigma\big(ty+(1-t)c_\sigma(x)\big)=t\delta_{\partial\sigma}(y)+(1-t)c_\sigma\circ\delta_M(x)\quad\text{for}\ (y,t)\in\partial((\gamma|_\sigma)^{-1}(x))\times[0,\,1].$$
Then $\delta_\sigma$ is a definable homeomorphism of $\sigma$ preserving $L|_\sigma$, $\gamma\circ\delta_\sigma=\delta_M\circ\gamma$ on $\sigma$, moreover, $\delta_\sigma=\delta_{\partial\sigma}$ on $\partial\sigma-\cup_{x\in\gamma(\sigma)}\partial((\gamma|_\sigma)^{-1}(x))$ because of the method of construction of $\delta_\sigma$, hence $\delta_\sigma$ is an extension of $\delta_{\partial\sigma}$ and $\delta_\sigma=\id$ on $\sigma\cap\gamma^{-1}(D)$.
Thus Lemma 6 is proved.\qed\vskip2mm\noindent
{\bf Lemma 7} (Weak local triangulations of definable $C^0$ maps).
{\it Let $K$ be a finite simplicial complex in $R^n$ with underlying polyhedron $X$ and let $H=(H_1,H_2):X\to R^2$ be a definable $C^0$ map such that $H_1\ge 0,\ H_2\ge 0$ and the zero sets of $H_1$ and $H_2$ are the underlying polyhedra of some subcomplexes of $K$.
Then there exist a definable isotopy $\zeta_t,\ 0\le t\le 1$, of $X$ preserving $K$ and a definable neighborhood $N$ of $(0,\,\epsilon]\times\{0\}$ in $(0,\,\epsilon]\times R$ for some positive $\epsilon\in R$ such that the $H\circ\zeta_1|_{H^{-1}(\overline N)}$ is extendable to a PL map $\tilde H=(\tilde H_1,\tilde H_2):X\longrightarrow R^2$ satisfying $\tilde H_1\ge0,\ \tilde H_2\ge0,\ \tilde H_1^{-1}(0)=H_1^{-1}(0)$ and $\tilde H_2^{-1}(0)=H_2^{-1}(0)$.\par
If $H_1$ is PL from the beginning, then we can choose $\zeta_t$ so that $H_1\circ\zeta_1=H_1$ on $H^{-1}(\overline N)$.}\vskip2mm
Note that a ``\,proper\," local triangulation of $H$ is impossible in general, e.g.\ the blowing-up $[0,\,1]^2\ni(x_1,x_2)\allowbreak\to(x_1,x_1x_2)\in R^2$ does not admit a triangulation locally at each point of $\{0\}\times[0,\,1]$.
Hence the neighborhood $N\cap\{x\}\times R$ of $(x,0)$ in $\{x\}\times R$ for $x\in (0,\,\epsilon]$ converges to $\{(0,0)\}$ as $x\to0$ in the sense that $\dis(N\cap\{x\}\times R,0\times0)\to0$.\vskip2mm\noindent
{\it Proof of Lemma 7}.
If $\Ima H\cap N$ is of dimension $\le1$ for some $N$ in Lemma 7 (i.e., $\Ima H\cap N$ is the empty set or $(0,\,\epsilon]\times\{0\}$ for some smaller $N$), then Lemma 7 is obvious or follows from Triangulation theorem of definable $C^0$ functions (1).
Hence we assume that $\Ima H\cap N$ is of dimension $2$ for any $N$.
Let us consider the former half of Lemma 7.
We only construct a definable homeomorphism $\zeta$ of $X$ preserving $K$ such that the $H\circ\zeta|_{H^{-1}(\overline N)}$ is extendable to a PL map because the other condition that $\zeta$ is the finishing homeomorphism of some isotopy can be clearly satisfied by the usual argument.
We use the same argument as in the proof of Triangulation theorem of definable $C^0$ functions (1).
Set $A=\graph H$, and let $R^n\times R^2\stackrel{p_1}{\longrightarrow}\cdots\stackrel{p_n}{\longrightarrow}R^2$ denote the projections ignoring the respective first factors.
If there is a sequence of definable $C^1$ stratifications $\{A''_{j''}\}\stackrel{p_1}{\longrightarrow}\cdots\stackrel{p_n}{\longrightarrow}\{p_n\circ\cdots\circ p_1(A''_{j''})\}$ such that $\{A''_{j''}\}$ is a stratification of $A$ compatible with $\{A\cap(\sigma\times R^2):\sigma\in K\}$, conditions 2.6 and 2.7 are satisfied and $\{p_n\circ\cdots\circ p_1(A''_{j''})\}$ is the set of the open simplices in some simplicial complex in $R^2$, then by Proposition 2.9 and its proof there exists a definable homeomorphism $\tau$ of $X$ preserving $K$ such that $H\circ\tau$ is PL.
We will obtain such a sequence of definable $C^1$ stratifications by shrinking $A$ to $A\cap(R^n\times\overline N)$.\par
Let $p:R^n\times R^2\to R^2$ denote the projection and $p_{n+1}:R^2\to R$ the projection to the first (not last) factor.
Set $B_{(u,v)}=\{x\in R^n:(x,u,v)\in B\}$ for $(u,v)\in R^2$ and for any subset $B$ of $R^n\times R^2$.
Let $p|_A:\{A_j\}_j\to\{ p(A_j)\}_j$ be a definable $C^1$ stratification of $p|_A$ compatible with $(\{A\cap(\sigma\times R^2):\sigma\in K\},\emptyset)$.
Let $A'$ denote the union of the $A_j$ such that $A_{j(u,v)}$ is of dimension smaller than $n$ for each $(u,v)\in R^2$.
Set $S^{n-1}=\{\lambda\in R^n:|\lambda|=1\}$, let $T_{(u,v)}\ (\subset S^{n-1})$ denote the closure of the set of singular direction for $A'_{(u,v)}$, and set $T=\{(\lambda,u,v)\in S^{n-1}\times R^2:\lambda\in T_{(u,v)}\}$.
Then $T$ is definable, and $T_{(u,v)}$ is a definable set of dimension smaller than $n-1$ for each $(u,v)$.
If $\cup_{(u,v)\in R^2}T_{(u,v)}\not=S^{n-1}$, we can show a triangulation of $H$.
However, this assumption is not necessarily the case.
We will choose $N$ such that $\cup_{(u,v)\in\overline N}T_{(u,v)}\not=S^{n-1}$, i.e., there is a nonsingular direction for $A'_{(u,v)}$ for any $(u,v)\in\overline N$.\par
We will find a definable closed subset, say $V_{(0,0)}$, of $S^{n-1}$ of dimension smaller than $n-1$ any definable neighborhood of which contains $\cup_{(u,v)\in\overline N}T_{(u,v)}$ for some $N$.
As $V_{(0,0)}$ we will choose $\cap_{u_0>0}\cup_{0\le u\le u_0}\overline{\cup_{0\le v\le h(u)}T_{(u,v)}}$ for some small nonnegative definable $C^0$ function $h$ on $R$ with zero set $\{0\}$.
The $V_{(0,0)}$ should contains $T_{(0,0)}$.
If $(\overline T)_{(0,0)}$ is of dimension smaller than $n-1$, then it satisfies the requirement.
However, it may be of dimension $n-1$ since the parameter of $T_{(u,v)}$ is of two variables.
(For example, the set $C=\{(x,u,v)\in [0,\,1]^3:x u=v\}$ has the property\,: $C=\overline{\cup_{(u,v)\not=(0,0)}C_{(u,v)}\times\{(u,v)\}},\ \dim C_{(u,v)}=0$ for $(u,v)\not=(0,0)$ with $C_{(u,v)}\not=\emptyset$, and $\dim C_{(0,0)}=1$.
If $C$ is of one parameter, $C=\cup_{u\in R}C_u\times\{u\}$, $C_0\not=\emptyset$ and $C=\overline{\cup_{u\not=0}C_u\times\{u\}}$, then $\dim C_0=\dim C_u$ and $C_u\not=\emptyset$ for some $u\not=0$ arbitrarily near 0.) 
Hence we choose another set smaller than $(\overline T)_{(0,0)}$ by reducing the problem to the case of one variable.
Set 
$$U=\cup_{u\in[0,\,1]}(T\cap(S^{n-1}\times\{u\}\times R)\overline),$$
where the notation $(\ \overline)$ denote the closure of a set $(\ )$.
Then $T\subset U$, $U$ is a definable subset of $S^{n-1}\times R^2$, and $U_{(u,v)}$ is closed and of dimension smaller than $n-1$ for each $(u,v)$.
However, $U_{(0,0)}$ does not necessarily satisfies the requirement.
We need to enlarge $U_{(0,0)}$ this time.
Next we set 
$$V=\cup_{v\in[0,\,1]}(U\cap(S^{n-1}\times R\times\{v\})\overline).$$
Then $U\subset V$, $V$ is a definable subset of $S^{n-1}\times R^2$, $V\cap(S^{n-1}\times R\times\{0\})$ is closed, $V_{(u,v)}$ is closed and of dimension smaller than $n-1$ for each $(u,v)$, and, moreover, $V_{(0,0)}$ fulfills the requirement as follows:\par
Let $Q$ be any definable neighborhood of $V_{(0,0)}$ in $S^{n-1}$.
We need to see that the set $\{(u,v)\in(0,\,\epsilon]\times R:Q\supset V_{(u,v)}\}$ is a neighborhood of $(0,\,\epsilon]\times\{0\}$ in $(0,\,\epsilon]\times R$ for some $\epsilon\in(0,\,1]$.
Let $N$ denote this set.
There exists $\epsilon\in(0,\,1]$ such that $Q\times[0,\,\epsilon]\times\{0\}$ is a neighborhood of $V\cap(S^{n-1}\times[0,\,\epsilon]\times\{0\})$ in $S^{n-1}\times[0,\,\epsilon]\times\{0\}$ since $V\cap(S^{n-1}\times[0,\,\epsilon]\times\{0\})$ is closed.
Hence $N\cap(\{u\}\times R)$ is a neighborhood of $(u,0)$ in $\{u\}\times[0,\,\epsilon]$ for each $u\in[0,\,\epsilon]$.
Then, since $N$ is definable, there is a positive definable function $h$ on $(0,\,\epsilon]$ such that $N\supset\{(u,v)\in(0,\,\epsilon]\times R:0\le v\le h(u)\}$.
Shrink $\epsilon$ so that $h$ is of class $C^0$ on $(0,\,\epsilon]$.
Then $N$ is a neighborhood of $(0,\,\epsilon]\times\{0\}$ in $(0,\,\epsilon]\times R$.\par
Thus we obtain $V_{(0,0)}$ and $N$ such that the $\cup_{(u,v)\in\overline N}T_{(u,v)}$ is contained in a small definable neighborhood of $V_{(0,0)}$.
Then $\cup_{(u,v)\in\overline N}T_{(u,v)}\not= S^{n-1}$.
Hence there exists a nonsingular direction for $A'_{(u,v)}$ for any $(u,v)\in\overline N$.
After changing linearly the coordinate system of $R^n$ we can assume that the element $(1,0,...,0)$ of $R^n$ or any element of $R^n$ near $(1,0,...,0)$ is not a singular direction for $A'_{(u,v)}$ for any $(u,v)\in\overline N$.
Let $\{N_k\}$ be a definable $C^1$ stratification of $\overline N$, and let $\{A'_{j'}\}\stackrel{p_1}{\longrightarrow}\{p_1(A'_{j'})\}$ be a definable $C^1$ stratification of $p_1|_{A\cap(R^n\times\overline N)}$ compatible with $(\{A_j\},\emptyset)$ such that $\{A'_{j'}\cap(R^n\times N_k)-A'\}=\{A_j\cap(R^n\times N_k)-A'\}$.
Then condition 2.6 is satisfied for $\{A'_{j'}\}\overset{p_1}{\longrightarrow}\{p_1(A'_{j'})\}$.\par
Repeat the same argument for $\{p_1(A'_{j'})\}$ and $R^{n-1}\times R^2\stackrel{p_2}{\longrightarrow}\cdots\stackrel{p_n}{\longrightarrow}R^2$.
Then shrinking $N$ we obtain a sequence of definable $C^1$ stratifications $\{A''_{j''}\}\stackrel{p_1}{\longrightarrow}\cdots\stackrel{p_n}{\longrightarrow}\{p_n\circ\cdots\circ p_1(A''_{j''})\}$ such that $\{A''_{j''}\}$ is a stratification of $A\cap(R^n\times\overline N)$ compatible with $\{A\cap(\sigma\times\overline N):\sigma\in K\}$ and condition 2.6 is satisfied.
Moreover, this holds after any small perturbation of $p_1,...,p_{n-1}$.
Here we can assume by Lemma 2.8 that condition 2.7 is also satisfied.
However, the remaining condition that $\{p_n\circ\cdots\circ p_1(A''_{j''})\}$ is the set of the open simplices in some simplicial complex in $R^2$ cannot be satisfied because of $\cup_{j''}p_n\circ\cdots\circ p_1(A''_{j''})=\overline N\cap\Ima H$.
We solve this problem as follows:\par
Shrinking $N$ we can assume that $\overline N\subset\Ima H$, $\overline N=\{(u,v)\in[0,\,\epsilon]\times R:0\le v\le h(u)\}$ and $\{p_n\circ\cdots\circ p_1(A''_{j''})\}$ is the natural stratification of $\overline N$, i.e., $\{p_n\circ\cdots\circ p_1(A''_{j''})\}$ consists of $(0,0),\ (\epsilon,0),\ (\epsilon,h(\epsilon)),\ \allowbreak\graph h|_{(0,\,\epsilon)},\ (0,\,\epsilon)\times\{0\},\,\{\epsilon\}\times(0,\,h(\epsilon))$ and $\{(u,v)\in(0,\,\epsilon)\times R:0<v<h(u)\}$, where $h$ is a nonnegative definable $C^1$ function on $[0,\,\epsilon]$ with zero set $\{0\}$.
Let $\sigma_N$ denote the 2-simplex in $R^2$ with vertices $(0,0),\,(\epsilon,0)$ and $(\epsilon,h(\epsilon))$, and let $\beta:\overline N\to\sigma_N$ be a definable homeomorphism fixing $[0,\,\epsilon]\times\{0\}$ and carrying each segment $l$ in $\overline N$ parallel to $\{0\}\times R$ to the segment $l'$ in $\sigma_N$ such that $l$ and $l'$ are extended to the same line.
Set $\tilde\beta(x,u,v)=(x,\beta(u,v))$ for $(x,u,v)\in R^n\times\overline N$.
Let us consider $\beta\circ H|_{H^{-1}(\overline N)},\,\tilde\beta(A\cap(R^n\times\overline N))$ and its definable $C^1$ stratification $\{\tilde\beta(A''_{j''})\}$ in place of $H|_{H^{-1}(\overline N)},\,A\cap(R^n\times\overline N)$ and $\{A''_{j''}\}$.
Then $\tilde\beta(A\cap(R^n\times\overline N))$ is the graph of $\beta\circ H|_{H^{-1}(\overline N)}$, and the sequence $\{\tilde\beta(A''_{j''})\}\stackrel{p_1}{\longrightarrow}\cdots\stackrel{p_n}{\longrightarrow}\{p_n\circ\cdots\circ p_1\circ\tilde\beta(A''_{j''})\}$ satisfies the remaining condition since $\{p_n\circ\cdots\circ p_1\circ\tilde\beta(A''_{j''})\}$ is the set of open simplices in the simplicial complex generated by $\sigma_N$.\par
There is another problem when we apply Proposition 2.9.
Namely, a triangulation of $\beta\circ H|_{H^{-1}(\overline N)}$ is not necessarily extended to $X$ through a definable homeomorphism preserving $K$.
To take this into account, the proof of Proposition 2.9 says that there exist a compact subpolyhedron $X_1$ of $X$ and a definable homeomorphism $\tau:X_1\to H^{-1}(\overline N)$ such that $\tau(\sigma\cap X_1)\subset\sigma$ for each $\sigma\in K$, $\beta\circ H\circ\tau$ is PL and $\tau$ is extended to a definable homeomorphism $\tau$ of $X$ preserving $K$.
Shrink $N$ to another $N_1$ so that $H^{-1}(\overline{N_1})\subset X_1$, which is possible because $\tau$ preserves the inverse images of $(0,0)$ and $R\times\{0\}$ under $H$ and hence $X_1$ is a neighborhood of $H^{-1}((0,\,\epsilon]\times\{0\})$ for shrunk $\epsilon$.
Then $\beta\circ H\circ\tau|_{H^{-1}(\overline{N_1})}$ is extendable to a PL map $\tilde H=(\tilde H_1,\tilde H_2):X\to R^2$.\par
It remains to find a definable homeomorphism $\zeta$ of $X$ preserving $K$ such that $\beta\circ H\circ\tau=H\circ\zeta$ on $H^{-1}(\overline{N_2})$ for some smaller $N_2$.
Let $\beta\circ H\circ\tau|_{X_1}:L\to M$ be a simplicial decomposition of $\beta\circ H\circ\tau|_{X_1}:X_1\to\sigma_N$ such that each simplex in $L$ is contained in some simplex in $K$.
Let $\sigma_M$ be a small simplex in $R^2$ such that one of the vertices is the origin, another lies on $R\times\{0\}$ and $\sigma_M$ and $\beta(\sigma_M)$ are contained in a simplex in the barycentric subdivision $M'$.
Then we can replace $N_1$ by $\beta^{-1}(\sigma_M)-\{0\}$, and by the Alexander trick we can extend $\beta|_{\sigma_M}$ to a definable homeomorphism $\beta_M$ of $\sigma_N$ preserving $M$.
Hence existence of $\zeta$ follows from Lemma 6.
Moreover, the conditions $\tilde H_1\ge0,\ \tilde H_2\ge0,\ \tilde H_1^{-1}(0)=H_1^{-1}(0)$ and $\tilde H_2^{-1}(0)=H_2^{-1}(0)$ are obviously satisfied.
Thus the former half of Lemma 7 is proved.\par
Let us consider the latter half of Lemma 7.
Apply Fact 2.5 to PL functions $H_1$ and $\tilde H_1$.
Then there is a PL isotopy $\zeta'_t,\ 0\le t\le1$, of $X$ preserving $K$ such that $\tilde H_1\circ\zeta'_1=H_1$ on $H_1^{-1}([0,\epsilon'])$ for some sufficiently small $\epsilon'>0$.
Replace $\epsilon$ with $\epsilon'$, $N$ with $N\cap[0,\,\epsilon']\times R$, $\zeta_t$ with $\zeta_t\circ\zeta'_t$ and $\tilde H$ with $\tilde H\circ\zeta'_1$.
Then the latter half follows, and we complete the proof of Lemma 7.\qed
\section{Proof of Theorem 4}
{\it Proof of the existence in Theorem 4}.
By Triangulation theorem of definable sets we can assume that $X$ and $X_i$ are the unions of some open simplices in some finite simplicial complex $K$ in $R^n$ with $|K|=\overline X$.
Let $v_\sigma$ denote the barycenter of a simplex $\sigma$.
For each $\sigma\in K$ we regard $\Int\sigma$ as a slit in $\overline X$ and enlarge it to its neighborhood in $\overline X$ by a semialgebraic (not semilinear) $C^0$ embedding $\tau_\sigma:\overline X-\Int\sigma\to\overline X-\Int\sigma$ as follows (Figure 2):
We set $\tau_\sigma=\id$ outside $\Int|\st(v_\sigma,K')|$.
In the case where $\sigma$ is a vertex, we define $\tau_\sigma$ on $\Int|\st(v_\sigma,K')|-\Int\sigma$ so that the $\overline{X-\Ima\tau_\sigma}$ is $\sigma*|\lk(\sigma,K'')|$ and $\tau_\sigma$ linearly carries each open segment $l$ joining $\sigma$ with a point in $|\lk(\sigma,K')|$ into $l$.
Here $A*B$ denotes the join of compact polyhedra $A$ and $B$, i.e. the smallest polyhedron containing $A$ and $B$ under the condition that any two distinct segments, with ends in $A$ and $B$, do not meet except in $A\cup B$, and $\lk(\sigma,K'')$ denotes the link of $\sigma$ in $K''$.
Note 
$$\tau_\sigma((1-t)\sigma+t x)=(1-t)\sigma/2+(1+t)x/2\quad\text{for}\ (x,t)\in|\lk(\sigma,K')|\times(0,\,1].$$
In the general case of $\sigma$ we naturally extend $\tau_\sigma$ to $\Int|\st(v_\sigma,K')|-\Int\sigma$ so that $\overline{X-\Ima\tau_\sigma}=\cup_{\sigma_1\in K'',\Int\sigma_1\subset\Int\sigma}|\st(\sigma_1,K'')|$ and $\tau_\sigma$ linearly carries each open segment $l$ joining a point in $\Int\sigma$ with a point in $|\lk(v_\sigma,K')|$ into $l$.
Then 
$$\Ima\tau_{\sigma_1}\circ\tau_{\sigma_2}=\Ima\tau_{\sigma_1}\cap\Ima\tau_{\sigma_2}\quad\text{for}\ \sigma_1\not=\sigma_2\in K\ \text{with}\ \dim\sigma_1\le\dim\sigma_2.$$
\begin{figure}[h]
\begin{pspicture}(3.7,0.7)(12.7,2.4)
\pspolygon(9.9,1.5)(11.4,0.7)(12.9,1.5)(11.4,2.3)
\psline(11.4,0.7)(11.4,2.3)
\pspolygon(10.9,1.5)(11.4,0.7)(11.9,1.5)(11.4,2.3)
\pspolygon(11.15,1.5)(11.4,0.7)(11.65,1.5)(11.4,2.3)
\psdot(11.4,1.5)
\rput(10.66,1.3){$\sigma$}
\pscurve{->}(10.8,1.25)(11.1,1.22)(11.36,1.2)
\rput(10.65,1.68){$v_\sigma$}
\pscurve{->}(10.83,1.67)(11.1,1.6)(11.33,1.52)
\rput(10,1.15){$v_l$}
\pscurve{->}(10.17,1.2)(10.5,1.39)(10.87,1.5)
\rput(10.7,0.71){$v_1$}
\pscurve{->}(10.85,0.7)(11.1,0.68)(11.35,0.69)
\rput(10.05,2.3){$|\!\st(\sigma,K)|$}
\rput(8.6,1.6){$|\!\lk(\sigma,K)|$}
\pscurve{->}(9.35,1.6)(9.6,1.65)(9.85,1.55)
\pscurve{->}(10.07,2.1)(10.2,1.9)(10.4,1.52)
\pscurve{->}(10.8,2.3)(11.4,2.4)(11.8,2.2)(11.9,1.9)
\pscurve{->}(10.83,1.67)(11.1,1.6)(11.3,1.52)
\rput(13.37,1){$|\!\st(v_\sigma,K')|$}
\pscurve{->}(12.5,1.08)(12.1,1.27)(11.75,1.5)
\pscurve{->}(12.5,1)(12.1,1.1)(11.04,1.39)
\rput(13.36,2){$\overline X\!-\!\Ima\tau_\sigma$}
\pscurve{->}(12.53,1.95)(11.9,1.8)(11.25,1.7)
\pscurve{->}(12.53,1.86)(12.1,1.75)(11.48,1.59)
\end{pspicture}
\caption{Enlarged slit $\overline X-\Ima\tau_\sigma$}\end{figure}
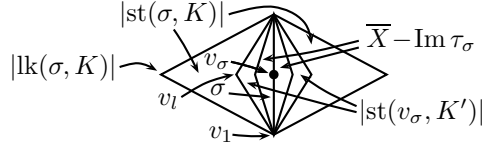\vskip2mm\noindent
Fig.\ 2 describes the image of $\tau_\sigma$, where $K$ is generated by two 2-simplices with common 1-face $\sigma$.\par
Let $\{\sigma_1,...,\sigma_k\}=\{\sigma\in K:X\cap\Int\sigma=\emptyset\}$, order the simplices so that $\dim\sigma_1\le\cdots\le\dim\sigma_k$, and set $\tau=\tau_{\sigma_1}\circ\cdots\circ\tau_{\sigma_k}$.
Then $\tau$ is a semialgebraic $C^0$ embedding of $(X;X_i)$ into $(X;X_i)$, $\overline{\Ima\tau}\subset X$, 
$$\overline X-\overline{\Ima\tau}=\overline X-\cap_{j=1}^k\overline{\Ima\tau_{\sigma_j}}=\{x\in\Int|\st(v,K'')|:v\in K^{\prime\prime0}-X\}, $$
hence $\overline X-\overline{\Ima\tau}$ is a neighborhood of $\overline X-X$ in $\overline X$, $(\Ima\tau;\tau(X_i))$ is a family of semilinear sets and is semialgebraically homeomorphic to $(X;X_i)$, and it is, moreover, standard for the following reason:\par
Set $V=\overline{\Ima\tau}-\Ima\tau$.
We will define a collar of $V$ in $\overline{\Ima\tau}$ by using alterations of $\tau$.
Let the above $\tau_\sigma$ be rewritten as $\tau_{\sigma 1/2}$ and define $C^0$ embeddings $\tau_{\sigma t_0}:\overline X-\Int\sigma\to\overline X-\Int\sigma$ for each $t_0\in(0,\,1]$ so that each open segment $l$ joining a point $x$ in $\Int\sigma$ with a point $y$ in $|\lk(v_\sigma,K')|$, the restriction of $\tau_{\sigma t_0}$ to $l$ is linear and its image is the open segment joining $t_0x+(1-t_0)y$ with $y$.
Let us define $\tau_{t_0}$ to be $\tau_{\sigma_1 t_0}\circ\cdots\circ\tau_{\sigma_k t_0}$ for $t_0\in(0,\,1]$.
Then $\tau_{t_0}$ are semialgebraic $C^0$ embeddings of $(X;X_i)$ into $(X;X_i)$, $\overline{\Ima\tau_{t_0}}\subset\Ima\tau_{t_0'}$ for $t_0<t_0'$, and $\tau_1=\id$.
We define a candidate of a collar $\pi:V\times[0,\,1]\to\overline{\Ima\tau_{1/2}}$ by $\pi(x,t)=\tau_{1-t/2}(x)$.
Then $\pi$ is a semialgebraic $C^0$ embedding, $\pi(\cdot,0)=\id$, $\pi(V\times[0,\,1))\ (=\overline{\Ima\tau_{1/2}}-\overline{\Ima\tau_{1/4}})$ is an open semilinear neighborhood of $V$ in $\overline{\Ima\tau_{1/2}}$, and $\pi^{-1}(\tau_{1/2}(X_i))=(V\cap X_i)\times(0,\,1]$.
Thus $\pi$ satisfies the conditions on a collar except that $\pi$ is of semilinear class.
We need to modify $\pi$ to a semilinear embedding, say, $\xi:V\times[0,\,1]\to\overline{\Ima\tau_{1/2}}$.
It is easy in the case $\overline V=V$ and will become apparent.\par
Consider the other case.
There are two problems.
One is that for any convergent directed family of points $V\ni x_\alpha\to x\in\{x\in\sigma:\sigma\in K'',\,\sigma\subset\overline V-V\}$, the $\pi(\{x_\alpha\}\times[0,\,1])$ converges to the point $\{x\}$ but if $\xi$ exists then $\xi(\{x_\alpha\}\times[0,\,1])$ does not necessarily converge to a point as in Fig.\ 1 (a).
We solve this problem by shrinking $V\times[0,\,1]$---the domain of definition of $\xi$---to the next set $\Xi$ so that the $\xi(\Xi\cap(\{x_\alpha\}\times[0,\,1]))$ converges to a point for the time being.\par
Note that $V$ is the union of some open simplices in $K''$ and
$$\pi^{-1}(\sigma\cap\overline{\Ima\tau_{1/2}})=(V\cap\sigma)\times[0,\,1]\quad\text{for}\ \sigma\in K'\ \text{with}\ V\cap\sigma\not=\emptyset.$$
We will modify $\pi$ to $\xi$ keeping the last equality because the equality $\xi^{-1}(\tau_{1/2}(X_i))=(V\cap X_i)\times[0,\,1]$ follows from it and hence $\xi$ is a collar of $V$ in $(\overline{\Ima\tau_{1/2}};\Ima\tau_{1/2},\tau_{1/2}(X_i))$ if $\xi$ is a collar of $V$ in $\overline{\Ima\tau_{1/2}}$.
Let $\theta:K''|_{\overline V}\to\{0,1,[0,\,1]\}$ be the simplicial map such that $\theta=0$ on $K^{\prime\prime0}\cap(\overline V-V)$ and $\theta=1$ at the other vertices, and set 
$$\Xi=\{(x,t)\in V\times[0,\,1]:0\le t\le\theta(x)\}.$$
Note that for $x\in V$, $\theta(x)=0$ if and only if $\pi(\{x\}\times[0,\,1])$ is a point.
We shrink $V\times[0,\,1]$ to $\Xi$.
Then the map of natural modification $\Xi\ni(x,t)\to\pi(x,t/\theta(x))\in\Ima\pi$ is well defined and extendable to a semialgebraic $C^0$ embedding $\overline\Xi\to\overline{\Ima\pi}$.
However, the map is not necessarily semilinear even if $\overline V=V$, i.e., $\theta\equiv1$.
This is the second problem.
We can modify the map by the following statement, which is easily proved by the Alexander trick:\vskip1mm\noindent
{\it Given a finite cell complex $L$ and a definable $C^0$ embedding $\rho:|L|\to R^n$ such that $\rho(\sigma)$ is polyhedral for each $\sigma\in L$, there exists a PL embedding $\rho':|L|\to R^n$ such that $\rho'(\sigma)=\rho(\sigma)$ for $\sigma\in L$ and $\rho'=\rho$ on $\sigma$ if $\rho|_\sigma$ is PL.}\par
Hence there exists a semilinear homeomorphism $\xi:\Xi\to\Ima\pi$ such that $\Xi(\cdot,0)=\id$ on $V$ and $\xi(\Xi\cap((\sigma\cap V)\times[0,\,1]))=\pi((\sigma\cap V)\times[0,\,1])$ for $\sigma\in K'$ since $\overline{\pi((\sigma\cap V)\times[0,\,1])}$ are polyhedral.
Thus the first and second problems are solved at once.\par
It remains to enlarge $\Xi$ to $V\times [0,\,1]$.
Clearly there exists a semilinear homeomorphism $\xi':V\times[0,\,1]\to\Xi$ of the form $\xi'(x,t)=(\xi''(x,t),t)$ for $(x,t)\in V\times[0,\,1]$ such that 
$$\xi'((V\cap\sigma)\times[0,\,1])=\Xi\cap((V\cap\sigma)\times[0,\,1])\quad\text{for}\ \sigma\in K'.$$
Hence $\xi\circ\xi':V\times[0,\,1]\to\overline{\Ima\tau_{1/2}}$ is the required collar of $V$ in $\overline{\Ima\tau_{1/2}}$, and the existence is proved.\qed\vskip2mm\noindent
{\bf Remark 8.}
{\it The map $\tau^{-1}:\Ima\tau\to X$ in the above proof is extendable to a semialgebraic $C^0$ map $\overline{\tau^{-1}}:\overline{\Ima\tau}\to\overline X$ such that $\overline{\tau^{-1}}(\sigma\cap\overline{\Ima\tau})=\sigma$ for $\sigma\in K'$ with $\Int\sigma\subset X$.}\vskip2mm
For a technical reason, we will prove the uniqueness in Theorem 4 in the following general form:\vskip1mm\noindent
{\bf Statement 9.} {\it A definable homeomorphism between standard families of semilinear sets is definably isotopic to a semilinear homeomorphism through homeomorphisms.}\vskip1mm
We prepare for proving statement 9.
Let $(X;X_i)$ be a finite family of semilinear sets bounded in $R^n$ such that $X$ is not compact, and $K$ be as in the above proof.
Set $W=\overline X-X$ and $U_K=\{x\in|\st(v,K'')|:v\in K^{\prime\prime0}\cap W\}$, i.e., the closure of $X-\Ima\tau$ in the above proof, (note that $\Ima\tau=(X-U_K)\cup\{x\in\sigma:\sigma\in K'',\,\sigma\subset\overline W-W\}$), let $L$ be another simplicial decomposition of $\overline X$ having the same property as $K$, let $L_1$ and $L_2$ be derived subdivisions of $L$ and $L_1$, respectively, and define $U_L$ by $L_2$ likewise $U_K$ by $K''$.
Then the $(U_K;U_K\cap X,U_K\cap X_i)$ and $(U_L;U_L\cap X,U_L\cap X_i)$ are PL homeomorphic for the following reason:
This is Theorem 3.24 in \cite{R-S} and a small generalization of Fact 2.4.
We briefly repeat the proof in \cite{R-S}.
We can assume that $L$ is a subdivision of $K$ by replacing $L$ with a simplicial subdivision of the cell complex $\{\sigma_K\cap\sigma_L:\sigma_K\in K,\,\sigma_L\in L\}$.
Let $\tilde K_1$ and $\tilde K_2$ be any derived subdivisions of $K$ and $\tilde K_1$, respectively, such that $\tilde K_1|_{\overline W}=K'|_{\overline W}$ and $\tilde K_2|_{\overline W}=K''|_{\overline W}$, let $\tilde L_1$ and $\tilde L_2$ be any derived subdivisions of $L$ and $\tilde L_1$, respectively, such that $\tilde L_1|_{\overline W}=L_1|_{\overline W}$ and $\tilde L_2|_{\overline W}=L_2|_{\overline W}$, and define $\tilde U_K$ and $\tilde U_L$ by $\tilde K_2$ and $\tilde L_2$ likewise $U_K$ by $K''$.
Then by Fact 2.2, $(U_K;U_K\cap X,U_K\cap X_i)$ and $(U_L;U_L\cap X,U_L\cap X_i)$ are PL homeomorphic to $(\tilde U_K;\tilde U_K\cap X,\tilde U_K\cap X_i)$ and $(\tilde U_L;\tilde U_L\cap X,\tilde U_L\cap X_i)$, respectively, and by Fact 2.3 there exist $\tilde K_1,\,\tilde L_1,\,\tilde K_2$ and $\tilde L_2$ such that 
$$(\tilde U_K;\tilde U_K\cap X,\tilde U_K\cap X_i)=(\tilde U_L;\tilde U_L\cap X,\tilde U_L\cap X_i).$$
Hence $(U_K;U_K\cap X,U_K\cap X_i)$ and $(U_L;U_L\cap X,U_L\cap X_i)$ are PL homeomorphic.
Note also that the homeomorphism is extendable to a PL homeomorphism of $(\overline X;X,X_i)$.\par
We call this property on $(U_K;U_K\cap X,U_K\cap X_i)$ and $(U_L;U_L\cap X,U_L\cap X_i)$ {\it the property of invariance}.
Using this we will show the following remark, by which we can assume that any standard semilinear family is of the form $(\Ima\tau;\tau(X_i))$ when we prove statement 9:\vskip2mm\noindent
{\bf Remark 10.}
{\it Let $(X;X_i),\,K,\,W,\,U_K$ and $\tau$ be the same as above.
Assume that $(X;X_i)$ is standard.
Then $(X;X_i)$ is semilinearly homeomorphic to $(\Ima\tau;\tau(X_i))$.}\vskip2mm\noindent
{\it Proof of Remark 10}.
Let $\phi:W\times[0,\,1]\to\overline X$ be a collar of $W$ in $(\overline X;X,X_i)$ and let $\overline\phi:\overline W\times[0,\,1]\to\overline X$ be the PL extension of $\phi$ (Remark 3, (i)).
We will reduce the problem to the one on $W\times[0,\,1]$ through $\phi$.
Let $Y$ denote the boundary of $\Ima\phi$ in $\overline X$.
Then $\overline\phi^{-1}(Y)\cap(W\times\{0\})=\emptyset$ since $\phi(W\times[0,\,1))$ is a neighborhood of $W$ in $\overline X$.
Let $h$ be a PL function on $\overline W$ such that $h>0$ on $W$, $h<1$ and the set $Z=\{(x,t)\in W\times[0,\,1]:t\le h(x)\}$ does not intersect with $\overline\phi^{-1}(Y)$.
By the property of invariance we can subdivide $K$ and assume that there is a simplicial decomposition $M$ of $\overline W\times[0,\,1]$ such that $\overline\phi:M\to K$ is simplicial and $Z$ is the union of some open simplices in $M$.
Let us define $U_M$ by $W\times\{0\}$ in $\overline W\times[0,\,1]$ and the second barycentric subdivision $M''$.
Then $\phi(U_M\cap(W\times(0,\,1]))=U_K\cap X$ and $\phi^{-1}(\Ima\tau)=W\times[0,\,1]-U_M$.
Hence it suffices to find a semilinear homeomorphism $\rho$ from $W\times(0,\,1]$ to $W\times[0,\,1]-U_M$ such that $\rho(\phi^{-1}(X_i))\subset\phi^{-1}(X_i)$ and $\rho=\id$ outside of $Z$.
We define $\rho$ as follows:\par
Note that $U_M$ is contained in the set $\{(x,t)\in\overline W\times[0,\,1]:t\le h(x)/2\}$.
Let $\overline\rho$ be the PL homeomorphism from $(\overline W\times\{0\})\cup\graph h$ to $\partial U_M\cup\graph h$ of the form $\overline\rho(x,t)=(x,\overline\rho_1(x,t))$ such that $\overline\rho=\id$ on $\graph h$, where $\partial U_M$ denotes the boundary of $U_M$ in $\overline W\times[0,\,1]$.
By the Alexander trick we can extend $\overline\rho$ to a PL homeomorphism $\overline\rho:\overline W\times[0,\,1]\to\overline W\times[0,\,1]-\Int U_M$ of the same form such that $\overline\rho=\id$ on $\{(x,t):t\ge h(x)\}$.
Then the restriction $\rho=\overline\rho|_{W\times(0,\,1]}$ satisfies the conditions.\qed\vskip2mm
We continue to prepare for proving statement 9.
We will replace $U_K$ by a set according to the proof of statement 9 as follows:
Let us define a semialgebraic (not semilinear) $C^0$ function $f_K$ on $\overline X-\{x\in\sigma:\sigma\in K'',\,\sigma\subset\overline W-W\}$ by 
$$f_K=\left \{
\begin{array}{l}\!
0\quad\text{on}\ K^{\prime\prime0}\cap W\\
\!1\quad\text{on}\ K^{\prime\prime0}-\overline W,
\end{array}
\right.
\quad
f_K(\sum_{j=1}^{l_3}t_j v_j)=\sum_{j=l_1+1}^{l_2}t_j/\sum_{j=1}^{l_2}t_j $$
for $t_1,...,t_{l_3}\in[0,\,1]$ with $\sum_{j=1}^{l_3} t_j=1$ and $\sum_{j=1}^{l_2}t_j>0$ and for a simplex in $K''$ whose vertices $v_1,...,v_{l_1}$ are in $W$, $v_{l_1+1},...,v_{l_2}$ are outside of $\overline W$ and $v_{l_2+1},...,v_{l_3}$ are in $\overline W-W$.
Then $\overline X-\Dom f_K\subset\overline W-W\subset X$ (the notation $\Dom$ denotes the domain of definition); $f_K$ is linear on the simplex $\{\sum_{j=1}^{l_3}t_j v_j:t_j\in[0,\,1],\sum_{j=1}^{l_3}t_j=1,\sum_{j=1}^{l_2}t_j=c\}$ for each $c\in(0,\,1]$ and the above $v_j$; its image is $[0,\,1]$; $f_K$ is PL if and only if $\overline W=W$ (i.e., $X$ is locally closed in $R^n$); $\overline{\cup_{t\in[0,\,1)}f^{-1}_K(t)}=U_K$; the map $f_K:\{\sigma\in K'':\sigma\cap(\overline W-W)=\emptyset\}\to\{[0,\,1],\{0\},\{1\}\}$ is simplicial; $(\overline X;U_K,X,X_i)$ is PL homeomorphic to $(\overline X;\overline{f^{-1}_K([0,t])},X,X_i)$ for any $t\in(0,\,1)$ by the above argument of the property of invariance because there exists a derived subdivision $L$ of $K'$ such that $\overline{f^{-1}_K([0,\,t])}$ is defined by $L$ likewise $U_K$ by $K''$. Hence we can replace $U_K$ by $\overline{f^{-1}_K([0,\,1/2])}$.\vskip2mm
\begin{figure}[h]
\begin{pspicture}(-3,-2.3)(12.5,2.1)
\pspolygon(1,0.5)(3,0.7)(1.5,2)
\psline[linestyle=dashed, dash=2pt 1pt](1.98,1.05)(1.5,2)
\psline(1,0.5)(1.93,0.96)
\psline(2.07,0.98)(3,0.7)
\pscircle(2,1){0.07}
\rput(2,2.2){$X$}
\rput(2.1,0.3){\footnotesize$v_{l_1+1},...,v_{l_2}$}
\rput(3,1.7){\footnotesize$v_{l_1+1},...,v_{l_2}$}
\pscurve{-}(1.67,0.4)(1.66,0.55)(1.66,0.65)
\pscurve[linestyle=dashed, dash=1pt 0.7pt]{->}(1.66,0.65)(1.67,0.8)(1.74,1)
\pscurve{-}(1.81,0.4)(1.8,0.55)(1.8,0.65)
\pscurve[linestyle=dashed, dash=1pt 0.7pt]{->}(1.8,0.65)(1.82,0.8)(1.84,1)
\psline{->}(2,0.4)(2,0.6)
\pscurve{->}(2.18,0.4)(2.13,0.65)(2,0.79)
\pscurve{->}(2.7,1.58)(2.56,1.47)(2.29,1.35)
\pscurve{->}(2.8,1.58)(2.5,1.3)(2.14,1.17)
\pscurve{->}(2.9,1.58)(2.7,1.2)(2.5,0.85)
\pscurve{->}(3,1.58)(2.95,1)(3,0.73)
\rput(0.83,1.26){$W$}
\pscurve{->}(1.05,1.23)(1.25,1.3)(1.45,1.3)
\psline[linestyle=dashed, dash=2pt 1pt](4,0.5)(5,1)
\psline(5,1)(4.52,1.97)
\psline[linestyle=dashed, dash=2pt 1pt](4.47,1.97)(4,0.5)
\pscircle(4.5,2){0.07}
\rput(5,2.2){$W$}
\rput(5.5,1.7){\footnotesize$v_{1},...,v_{l_1}$}
\pscurve{->}(5.1,1.6)(5,1.55)(4.77,1.5)
\pscurve{->}(5.2,1.6)(5,1.36)(4.55,1.19)
\pscurve{->}(5.3,1.6)(5.2,1.3)(5.01,1.04)
\psline(6.3,0.5)(7.27,0.98)
\psline(6.3,0.5)(6.9,2)
\pscircle(7.3,1){0.07}
\rput(7,2.2){$\overline W-W$}
\rput(7.8,1.5){\footnotesize$v_{l_2+1},...,v_{l_3}$}
\psline{->}(7.15,1.65)(6.91,1.97)
\pscurve{->}(7,1.45)(6.8,1.36)(6.63,1.25)
\pscurve{->}(7.06,1.37)(6.6,0.85)(6.36,0.56)
\pscurve{->}(7.15,1.33)(7,1.1)(6.82,0.77)
\psline(1.8,-2)(2.55,-1.63)
\psline(1.8,-2)(2.3,-0.5)
\rput(2.4,-0.3){$\overline X-\Dom f_K$}
\psline(4,-2)(4.68,-1.7)
\psline(4,-2)(4.5,-0.5)
\pspolygon(4.68,-1.7)(5.2,-1.52)(4.97,-1.3)(4.83,-0.9)(4.55,-0.54)(5.05,-0.9)(5.35,-1.3)(5.75,-1.62)
\rput(6.3,-2.05){\footnotesize$f^{-1}_K(1/4)\cup(\overline X-\Dom f_K)$}
\rput(6.4,-0.4){\footnotesize$f^{-1}_K(1/2)\cup(\overline X-\Dom f_K)$}
\pscurve{->}(5.5,-0.65)(5.3,-0.8)(4.97,-1.1)
\pscurve{->}(4.6,-2)(4.4,-1.8)(4.5,-1.4)
\end{pspicture}
\caption{Inverse image of $f_K$}
\end{figure}
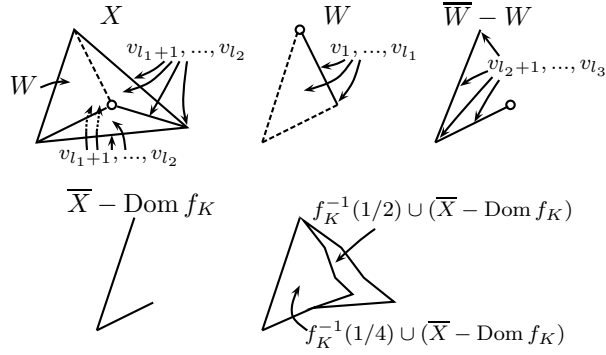\noindent
In Fig.\ 3, $K$ is generated by a 3-simplex, $W$ is the union of one open 2-simplex, one open 1-simplex and one 0-simplex in $K$, and $v_j$ point out the vertices in $K'$ only.\par
There is another reason why we introduce $f_K$.
For two definably homeomorphic standard families of semilinear sets we will define two $f_K$'s as above and prove that the families are semilinearly homeomorphic by comparing two $f_K$'s.
This idea is similar to the idea of applying Triangulation theorem of definable $C^0$ functions in the proof of Uniqueness theorem of definable triangulation.\vskip2mm
We have completed preparation of the proof in the case where the domain is locally closed.
The proof in the non-locally-closed case is long but the idea is the same as in the locally closed one.
We separate the proof and first prove the locally closed one to show the idea.\vskip2mm\noindent
{\it Proof of statement 9 in the locally closed case}.
Let $\eta:(Y;Y_i)\to(Z;Z_i)$ be a definable homeomorphism between standard families of semilinear sets in $R^n$.
Let $K_Y$ be a simplicial decomposition of $\overline Y$ such that $Y$ and $Y_i$ are the unions of some open simplices in $K_Y$.
Set $W_Y=\overline Y-Y$, and define a semialgebraic $C^0$ function $f_Y$ by $K_Y$ likewise $f_K$ by $K$.
Let $K_Z,\,W_Z$ and $f_Z$ be given for $(Z;Z_i)$ in the same way.
Assume that $Y$ is noncompact and locally closed.
We want to compare $f_Y$ and $f_Z$.
For this we need the condition that $\eta$ is extendable to a definable (not necessarily injective) $C^0$ map $\overline\eta:\overline Y\to\overline Z$, which is possible for the following reason:\par
Set $X=\graph\eta$ and $X_i=\graph\eta|_{Y_i}$, and let $p_Y:(X;X_i)\to(Y;Y_i)$ and $p_Z:(X;X_i)\to(Z;Z_i)$ denote the projections.
Then $p_Y$ and $p_Z$ are extendable to definable $C^0$ maps $\overline{p_Y}:\overline X\to\overline Y$ and $\overline{p_Z}:\overline X\to\overline Z$ respectively.
By Triangulation theorem of definable sets we can regard $\overline X$ as the underlying polyhedron of a finite simplicial complex $K$ such that $X$ and $X_i$ are the unions of some open simplices in $K$.
Let us define $\tau:(X;X_i)\to(X;X_i)$ as above.
Then the $(\Ima\tau;\tau(X_i))$ is a standard family of semilinear sets, and $\tau^{-1}|_{\Ima\tau}:\Ima\tau\to X$ is extendable to a semialgebraic $C^0$ map from $\overline{\Ima\tau}$ into $\overline X$ (Remark 8).
Hence $p_Y\circ(\tau^{-1}|_{\Ima\tau}):(\Ima\tau;\tau(X_i))\to(Y;Y_i)$ and $p_Z\circ(\tau^{-1}|_{\Ima\tau}):(\Ima\tau;\tau(X_i))\to(Z;Z_i)$ are definable homeomorphisms between standard families of semilinear sets and are extendable to definable $C^0$ maps $\overline{\Ima\tau}\longrightarrow\overline Y$ and $\overline{\Ima\tau}\longrightarrow\overline Z$ respectively.
Thus by replacing the $\eta:(Y;Y_i)\to(Z;Z_i)$ with the pair of $p_Y\circ(\tau^{-1}|_{\Ima\tau}):(\Ima\tau;\tau(X_i))\to(Y;Y_i)$ and $p_Z\circ(\tau^{-1}|_{\Ima\tau}):(\Ima\tau;\tau(X_i))\to(Z;Z_i)$, we can assume that $\eta$ is extendable to a definable $C^0$ map $\overline\eta:\overline Y\to\overline Z$.\par
Next we modify $K''_Y$ and $K''_Z$ to certain cell complexes so that $(Y;Y\cap\sigma:\sigma\in K''_Y)$ and $(Z;Z\cap\sigma:\sigma\in K''_Z)$ become standard since $(Y;Y\cap\sigma:\sigma\in K''_Y)$ and $(Z;Z\cap\sigma:\sigma\in K''_Z)$ are not standard if $Y$ is of dimension $>1$ locally at $\overline Y-Y$ as noted at the definition of a standard family.
We set
\begin{align*}\tilde Y=f^{-1}_Y((1/2,\,1]),\ \ \tilde Y_i=\tilde Y\cap Y_i,\ \ \tilde f_Y=2f_Y-1,\\
\tilde K_Y=\{\sigma\cap f^{-1}_Y([1/2,\,1]),\,\sigma\cap f^{-1}_Y(1/2):\sigma\in K''_Y\},\end{align*}
and owing to Remark 10, replace $(Y;Y_i),\,K''_Y$ and $f_Y$ by $(\tilde Y;\tilde Y_i),\,\tilde K_Y$ and $\tilde f_Y$ respectively.
Then we can assume the following:
$\tilde K_Y$ is a cellular (not simplicial if $Y$ is of dimension $>1$ locally at some point of $\overline Y-Y$) decomposition of $\overline Y$; $\tilde f_Y:\tilde K_Y\to\{0,1,[0,\,1]\}$ is a cellular map; for each cell $\sigma$ in $\tilde K_Y$ such that $\sigma\cap W_Y\not=\emptyset$ and $\sigma\not\subset W_Y$; $\sigma\cap\tilde f^{-1}_Y([0,\,1/2])$ is linearly isomorphic to $(\sigma\cap W_Y)\times[0,\,1]$; $\tilde K_Z$ satisfies the same conditions.
Hence $(Y;Y\cap\sigma:\sigma\in\tilde K_Y)$ and $(Z;Z\cap\sigma:\sigma\in\tilde K_Z)$ are standard.\par
There are then semialgebraic isotopies $\pi_{Y t}:\overline Y\to\overline Y$ and $\pi_{Z t}:\overline Z\to\overline Z$, $0\le t\le1$, of the identity maps preserving $\{\sigma\cap Y:\sigma\in\tilde K_Y\}$ and $\{\sigma\cap Y:\sigma\in\tilde K_Z\}$, respectively, such that $\Ima\pi_{Y t}=\tilde f^{-1}_Y([t/2,\,1])$ and $\Ima\pi_{Z t}=\tilde f^{-1}_Z([t/2,\,1])$.
Assume that $\tilde f_Y=\tilde f_Z\circ\overline\eta$ on $\tilde f^{-1}_Y([0,\,t_0])$ for some $t_0\in(0,\,1]$, which is possible as shown below.
Then there is a definable isotopy $\eta_t:(Y;Y_i)\to(Z;Z_i)$, $0\le t\le1$, of $\eta$ through homeomorphisms such that $\eta_1$ is extendable to a definable homeomorphism from $\overline Y$ onto $\overline Z$.
Therefore, statement 9 follows from Uniqueness theorem of definable triangulation and its supplement in the case where $Y$ and $Z$ are locally closed, and $\tilde f_Y=\tilde f_Z\circ\eta_1$ on $\tilde f_Y^{-1}([0,\,t_0])$.\par
It remains to reduce the problem to the case $\tilde f_Y=\tilde f_Z\circ\overline\eta$ on some neighborhood of $W_Y$ in $\overline Y$.
Set $g_Y=\tilde f_Z\circ\overline\eta$, which is a definable $C^0$ function on $\overline Y$ with zero set $W_Y\,(=\tilde f^{-1}_Y(0))$.
Apply Triangulation theorem of definable $C^0$ functions to $g_Y$.
Then there exists a definable isotopy $\xi_t,\ 0\le t\le1$, of $\overline Y$ preserving $\tilde K_Y$ such that the $g_Y\circ\xi_1$ is PL.
Next apply Fact 2.5 to $\tilde f_Y$ and $g_Y\circ\xi_1$.
Then we have a PL isotopy $\delta_t,\ 0\le t\le1$, of $\overline Y$ preserving $\tilde K_Y$ such that $g_Y\circ\xi_1\circ\delta_1=\tilde f_Y$ on some definable neighborhood of $W_Y$ in $\overline Y$.
Hence translating $\overline\eta$ to $\overline\eta\circ\xi_1\circ\delta_1$ through $\overline\eta\circ\xi_t\circ\delta_t,\ 0\le t\le1$, we can assume that $\tilde f_Y=\tilde f_Z\circ\overline\eta$ around $W_Y$.
Thus statement 9 is proved in the locally closed case.\qed\vskip2mm
When $Y$ is not locally closed, we proceed as follows:
First we show that the equality $\tilde f_Y=\tilde f_Z\circ\overline\eta$ can hold on a neighborhood of $W_Y$ in $\overline Y$.
However, the neighborhood is not of the form $\tilde f_Y^{-1}([0,\,t_0]),\ t_0\in(0,\,1)$.
In the last step of the proof we will move the neighborhood to some $\tilde f_Y^{-1}([0,\,t_0])$.\par
We prepare for the proof of statement 9 in the non-locally-closed case.
Assume that $Y$ is not locally closed.
Let $\eta:(Y;Y_i)\to(Z;Z_i),\,K_Y,\,\allowbreak W_Y,\,f_Y,\,K_Z,\,W_Z$ and $f_Z$ be the same as in the locally closed case.
Note that we can replace $(Y;Y_i)$ or $(Z;Z_i)$ by a standard family of semilinear sets through a semilinear homeomorphism and a definable homeomorphism if the definable homeomorphism is isotopic to a semilinear homeomorphism through homeomorphisms.
By such replacement we will reduce the problem to a simple case.
To be more precise, we will subdivide (or modify) $K''_Y$ and $K''_Z$ into $K_{jY}$ and $K_{jZ},\ j=1,2$, respectively, in sequence, we will define $f_{jY}$ and $f_{jZ}$ by $K_{j Y}$ and $K_{j Z}$, respectively, likewise $f_Y$ by $K''_Y$, and we will reduce the problem to the case where the following six conditions are satisfied in order:\vskip2mm\noindent
(a) (i) {\it For any $t\in[0,\,1)$ there exists a natural cellular isomorphism $\kappa_{Y t}$ from $K_{1Y}$ to the cell complex $\{\sigma\cap(f_{1Y}^{-1}(t)\overline),\sigma\cap(f^{-1}_{1Y}([t,\,1])\overline):\sigma\in K_{1Y}\}$ such that the map $\overline Y\times[0,\,1)\ni(y,t)\to\kappa_{Y t}(y)\in\overline Y$ is continuous;} (ii) {\it the restriction of $f_{1Y}$ to $\{\sigma\in K_{1Y}:\sigma\subset\Dom f_{1Y}\}$ is a cellular map into $\{0,1,[0,\,1]\}$;} (iii) {\it each cell in $K_{1Y}$ is described as $\sigma_1*\sigma_2$ for unique $\sigma_1,\sigma_2\in K_{1Y}$ such that $\sigma_1\subset\overline Y-\Dom f_{1Y}$ and $\sigma_2\subset\Dom f_{1Y}$; 
$$f_{1Y}(ty_1+(1-t)y_2)=f_{1Y}(y_2)\quad\text{for}\ (y_1,y_2,t)\in\sigma_1\times\sigma_2\times[0,\,1)\leqno{\rm(iv)}$$
for such $\sigma_1$ and $\sigma_2$;} (v) {\it the $\sigma_2\cap f^{-1}_{1Y}([0,\,1/2])$ is linearly isomorphic to $(\sigma_2\cap W_Y)\times[0,\,1]$ for each $\sigma_2\in K_{1Y}$ such that $\sigma_2\not\subset W_Y$, $\sigma_2\subset\Dom f_{1Y}$ and $\sigma_2\cap W_Y\not=\emptyset$;} (vi) {\it the $(Y;Y\cap\sigma:\sigma\in K_{1Y})$ is standard; (vii) the property of invariance holds for $K_{1Y}$ and any second derived subdivision of $K_Y$.
The $K_{1Z}$ also satisfies the same conditions, and $\{\Int\sigma\cap Z:\sigma\in K_{1Z}\}$ is compatible with $\{\eta(\sigma\cap Y):\sigma\in K_{1Y}\}$.}\vskip2mm\noindent
(b) {\it The $\eta$ is extended to a definable homeomorphism $\overline\eta:Y\cup(\overline{W_Y}-W_Y\overline)\to Z\cup(\overline{W_Z}-W_Z\overline)$ and $\overline\eta|_{(\overline{W_Y}-W_Y\overline)}$ is a PL homeomorphism onto $(\overline{W_Z}-W_Z\overline)$.}\vskip2mm\noindent
(c) {\it The $K_{1Y}$ and $K_{1Z}$ are subdivided into simplicial complexes $K_{2Y}$ and $K_{2Z}$, respectively, such that the $\overline\eta|_{(\overline W_Y-W_Y\overline)}:K_{2Y}|_{(\overline W_Y-W_Y\overline)}\allowbreak\to K_{2Z}|_{(\overline W_Z-W_Z\overline)}$ is simplicial and $\eta(\overline Y-\Dom f_{2Y})=\overline Z-\Dom f_{2Z}$.}\vskip2mm\noindent
(d) {\it Let $\phi_Y$ and $\phi_Z$ be any nonnegative PL functions on $\overline Y$ and $\overline Z$ with zero set $\overline Y-\Dom f_{2Y}$ and $\overline Z-\Dom f_{2Z}$, respectively, which always exists.
Then the functions $\phi_Y f_{2Y}$ and $(\phi_Z f_{2Z})\circ\eta$ on $Y\cap\Dom f_{2Y}$ are extendable to PL and definable $C^0$ functions on $\overline Y$ with zero set $\overline{W_Y}$, say, $\overline{\phi_Y f_{2Y}}$ and $\overline{(\phi_Z f_{2Z})\circ\eta}$, respectively, $\phi_Y=\phi_Z\circ\eta$ on $Y\cap G^{-1}([0,\,\epsilon]\times[0,\,1])\,(=Y\cap\phi^{-1}_Y([0,\,\epsilon]))$ for some $\epsilon\in(0,\,1)$, and there exists a definable neighborhood $N$ of $(0,\,\epsilon]\times\{0\}$ in $(0,\,\epsilon]\times[0,\,1]$ such that the $\overline{(\phi_Z f_{2Z})\circ\eta}|_{G^{-1}(\overline N)}$ is extendable to a nonnegative PL function $\tilde F_2$ on $\overline Y$ with zero set $\overline{W_Y}$, where $G=(G_1,G_2)=(\phi_Y,\overline{\phi_Y f_{2Y}})$.}\vskip2mm\noindent
(e) {\it It holds that $\tilde F_2=G_2$ on $G^{-1}(\overline N)$ and hence $\overline{\phi_Y\!f_{2Y}}=\overline{(\phi_Z f_{2Z})\circ\eta}$ there.}\vskip2mm\noindent
(f) {\it The $f_{2Z}\circ\eta$ coincides with $f_{2Y}$ on $Y\cap G^{-1}(\hat N)$ for some definable neighborhood $\hat N$ of $[\epsilon,\,1]\times\{0\}$ in $[\epsilon,\,1]\times[0,\,1]$.}\vskip2mm
Under (a) we will find a definable isotopy $\alpha_t,\ 0\le t\le1$, of $Y$ preserving $\{Y\cap\sigma:\sigma\in K_{1Y}\}$ by induction on the dimension of $Y$ such that $\eta\circ\alpha_1$ is semilinear.
(We do not directly show that $\eta$ is definably isotopic to a semilinear homeomorphism through homeomorphisms because the induction process becomes complicated for the direct construction.) 
To obtain (b) we introduce the induction hypothesis that under (a) there exists the required isotopy $\alpha_t,\,0\le t\le1$, of $Y$ for smaller dimensional $Y$.\vskip2mm\noindent
{\it Proof that (a) can be satisfied}.
Set
$$\tilde Y=(\overline Y-\Dom f_Y)\cup f^{-1}_Y((1/2,\,1]),\quad\tilde Y_i=\tilde Y\cap Y_i,$$
$$\tilde K_Y=\{\sigma\cap(\tilde Y\overline),\,\big(\sigma\cap f^{-1}_Y(1/2)\overline{\big)}:\sigma\in K''_Y\}\quad\text{and}\quad\tilde f_Y=2f_Y-1.$$
Then $\tilde K_Y$ is a cell complex; $((\tilde Y\overline);\tilde Y,\tilde Y_i)$ and $(\overline Y;Y,Y_i)$ are PL homeomorphic (Remarks 3, (iii) and 7); $\tilde f_Y$ is defined by $\tilde K_Y$ likewise $f_Y$ by $K''_Y$; $\tilde f_Y:\{\sigma\in\tilde K_Y:\sigma\subset\Dom\tilde f_Y\}\to\{0,\,1,[0,\,1]\}$ is cellular; each cell in $\tilde K_Y$ is of the form either $\sigma_1*(\sigma_2\cap f^{-1}_Y(1/2))$ or $\sigma_1*(\sigma_2\cap f^{-1}_Y([1/2,\,1]))$ for some $\sigma_1,\sigma_2\in K''_Y$ such that $\sigma_1\subset\overline Y-\Dom f_Y$ and $\sigma_2\subset\Dom f_Y$; 
\begin{align*}\tilde f_Y(t y_1+(1-t)y_2)=2f_Y(t y_1+(1-t)y_2)-1=2f_Y(y_2)-1=\tilde f_Y(y_2)\\
\qquad\qquad\text{for}\ (y_1,y_2,t)\in\sigma_1\times(\sigma_2\cap f_Y^{-1}(1/2))\times[0,\,1)\\
\text{or}\ (y_1,y_2,t)\in\sigma_1\times(\sigma_2\cap f_Y^{-1}([1/2,\,1]))\times[0,\,1))\end{align*}
for such $\sigma_1$ and $\sigma_2$; the correspondence $\sigma_3\cap f^{-1}_Y(1/2)\to\sigma_3\cap f^{-1}_Y(1/2+t/2)$ induces a cellular isomorphism $\kappa_{Y t}$ from $\tilde K_Y$ to $\{\sigma\cap(\tilde f^{-1}_Y(t)\overline),\sigma\cap(\tilde f^{-1}_Y([t,\,1])\overline):\sigma\in\tilde K_Y\}$ for each $t\in[0,\,1)$ and each 1-simplex $\sigma_3$ in $K''_Y$ with one vertex in $\Dom f_Y$ and the other in $\overline Y-\Dom f_Y$; the map $(\tilde Y\overline)\times[0,\,1)\ni(y,t)\to\tilde\kappa_{Y t}(y)\in(\tilde Y\overline)$ is continuous; $\sigma_2\cap\tilde f_Y^{-1}([0,\,1/2])$ is linearly isomorphic to $(\sigma_2\cap f^{-1}_Y(1/2))\times[0,\,1]$ for each $\sigma_2\in\tilde K_Y$ such that $\sigma_2\not\subset\overline{f_Y^{-1}(1/2)},\ \sigma_2\subset\Dom\tilde f_Y$ and $\sigma_2\cap f_Y^{-1}(1/2)\not=\emptyset$; $(\tilde Y;\tilde Y\cap\sigma:\sigma\in\tilde K_Y)$ is standard.
(These are not the case unless we replace $U_K$ by $\overline{f^{-1}_Y([0,\,1/2])}$.) 
Hence, we can replace $(Y;Y_i),\,K''_Y$ and $f_Y$ by $(\tilde Y;\tilde Y_i),\,\tilde K_Y$ and $\tilde f_Y$ respectively.
We use new notation $K_{1Y}$ and $f_{1Y}$ for $\tilde K_Y$ and $\tilde f_Y$ since $K_{1Y}$ is now a cellular decomposition of $\overline Y$.
In the same way we can modify $K''_Z$.
Then the conditions in (a) are satisfied except the one that $\{\Int\sigma\cap Z:\sigma\in K_{1Z}\}$ is compatible with $\{\eta(\sigma\cap Y):\sigma\in K_{1Y}\}$.
For this condition we modify the above construction and $\eta$ as follows:\par
There is a semialgebraic homeomorphism $\lambda_Y:(\tilde Y;\tilde Y\cap\sigma:\sigma\in\tilde K_Y)\to(Y;Y\cap\sigma:\sigma\in K''_Y)$ such that $\lambda_Y|_{(\tilde Y;\tilde Y_i)}:(\tilde Y;\tilde Y_i)\to(Y;Y_i)$ is definably isotopic to a semilinear homeomorphism $(\tilde Y;\tilde Y_i)\to(Y;Y_i)$ through homeomorphisms but not to $(\tilde Y;\tilde Y\cap\sigma:\sigma\in\tilde K_Y)\to(Y;Y\cap\sigma:\sigma\in K''_Y)$.
First we only replace $(Y;Y_i)$, $K''_Y$ and $f_Y$ by $(\tilde Y;\tilde Y_i)$, $K_{1Y}$ and $f_{1Y}$ and keep the notation $(Y;Y_i)$ and $\eta$ for $(\tilde Y;\tilde Y_i)$ and $\eta\circ\lambda_Y$.
Next by Triangulation theorem of definable sets there is a definable isotopy $\beta_t,\ 0\le t\le1$, of $\overline Z$ preserving $K_Z$ such that $\beta_1(\overline{\eta(Y\cap\sigma)})\ (\sigma\in K_{1Y}$) are polyhedra.
Thirdly, we replace $\eta$ by $\beta_1\circ\eta$.
Then by subdividing $K_Z$ we can assume that $\{\Int\sigma\cap Z:\sigma\in K_Z\}$ is compatible with $\{\eta(\sigma\cap Y):\sigma\in K_{1Y}\}$.
Fourthly, we define $\tilde K_Z$ and $(\tilde Z;\tilde Z_i)$ by this $K_Z$.
Finally, we replace $(Z;Z_i),\ K''_Z,\ f_Z$ and $\eta$ by $(\tilde Z;\tilde Z_i),\ K_{1Z},\ f_{1Z}$ and $\lambda_Z^{-1}\circ\eta$ and keep the notation $(Z;Z_i)$ and $\eta$, where $\lambda_Z$ is defined in the same way as $\lambda_Y$.
Then the remaining condition in (a) is satisfied.\qed\vskip2mm\noindent
{\it Proof that (b) can be satisfied}.
As we have mentioned, we use the induction hypothesis here.
Obviously $\dim\,(\overline{W_Y}-W_Y)<\dim Y$.
By Remark 3, (ii), $(\overline{W_Y}-W_Y;Y_i\cap\overline W_Y)$ is standard.
It is also easy to see that $(\overline{W_Y}-W_Y;Y_i\cap\overline{W_Y}),\,K_{1Y}|_{(\overline{W_Y}-W_Y\overline)},\,\eta|_{\overline{W_Y}-W_Y},\allowbreak\,(\overline{W_Z}-W_Z;Z_i\cap\overline{W_Z})$ and $K_{1Z}|_{(\overline{W_Z}-W_Z\overline)}$ satisfy the conditions on $(Y;Y_i),\,K_{1Y},\,\eta,\,(Z;Z_i)$ and $K_{1Z}$ in (a).
Hence by the induction hypothesis there exists a definable isotopy $\alpha_{Wt},\ 0\le t\le 1$, of $\overline{W_Y}-W_Y$ preserving $\{\sigma\cap\overline{W_Y}-W_Y:\sigma\in K_{1Y}\}$ such that $\eta\circ\alpha_{W1}$ is semilinear.
Then by Remark 3, (iii), $\eta\circ\alpha_{W1}$ is extended to a PL homeomorphism from $(\overline{W_Y}-W_Y\overline)$ onto $(\overline{W_Z}-W_Z\overline)$ since $\overline{W_Y}-W_Y$ and $\overline{W_Z}-W_Z$ are standard, and by Lemma 5 we can extend $\alpha_{Wt}$ to a definable isotopy $\tilde\alpha_{Wt},\ 0\le t\le 1$, of $Y$ preserving $\{Y\cap\sigma:\sigma\in K_{1Y}\}$.
(To be precise, by induction on $l\in\N$, we construct a definable isotopy $\alpha_{l t},\ 0\le t\le1$, of $Y\cap|K^l_{1Y}|$ preserving $\{Y\cap\sigma:\sigma\in K^l_{1Y}\}$ such that $\alpha_{l t}$ is an extension of $\alpha_{l-1t}$ and $\alpha_{l t}=\alpha_{W t}$ on $|K^l_{1Y}|\cap\overline W_Y-W_Y$.) 
Hence by replacing $\eta$ with $\eta\circ\tilde\alpha_{W1}$ we can assume (b), i.e., that the extension $\overline\eta|_{(\overline{W_Y}-W_Y\overline)}$ of $\eta|_{\overline{W_Y}-W_Y}$ is a PL homeomorphism onto $(\overline{W_Z}-W_Z\overline)$.\qed\vskip2mm\noindent
{\it Proof that (c) can be satisfied}.
Let $\hat K_{1Y}$ and $\hat K_{1Z}$ be simplicial subdivisions of $K_{1Y}$ and $K_{1Z}$, respectively, such that $\overline\eta|_{(\overline{W_Y}-W_Y\overline)}:\hat K_{1Y}|_{(\overline{W_Y}-W_Y\overline)}\to\hat K_{1Z}|_{(\overline{W_Z}-W_Z\overline)}$ is an isomorphism, and set $K_{2Y}=\hat K''_{1Y}$ and $K_{2Z}=\hat K''_{1Z}$.
Then $\eta(\overline Y-\Dom f_{2Y})=\overline Z-\Dom f_{2Z}$.
Hence replacing $K_{1Y},\,f_{1Y},\,K_{1Z}$ and $K_{1Z}$ with $K_{2Y},\,f_{2Y},\,K_{2Z}$ and $K_{2Z}$ we assume (c).\qed\vskip2mm
Now we define a simplicial map $\phi_Y:K_{2Y}\to\{0,1,[0,\,1]\}$ concretely by 
$$\phi_Y=
\left \{
\begin{array}{l}\!
0\quad\text{on}\ K^0_{2Y}\cap(\overline{W_Y}-W_Y)\\
\!1\quad\text{on}\ K^0_{2Y}-(\overline{W_Y}-W_Y)
\end{array}
\right.$$
and a simplicial map $\phi_Z:K_{2Z}\to\{0,1,[0,\,1]\}$ in the same way for simplicity of notation, although the next argument works for general $\phi_Y$ and $\phi_Z$.
Then $\phi_Y^{-1}(0)=\overline Y-\Dom f_{2Y}$, $\phi_Z^{-1}(0)=\overline Z-\Dom f_{2Z}$, and $\phi_Y\!f_{2Y}$, $\phi_Z f_{2Z}$ and $(\phi_Z f_{2Z})\circ\eta$ are extendable to a simplicial map $\overline{\phi_Y\!f_{2Y}}:K_{2Y}\to\{0,1,[0,\,1]\}$ with zero set $\overline{W_Y}$, a simplicial map $\overline{\phi_Z f_{2Z}}:K_{2Z}\to\{0,1,[0,\,1]\}$ with zero set $\overline{W_Z}$ and a definable $C^0$ function $\overline{(\phi_Z f_{2Z})\circ\eta}$ on $\overline Y$ with zero set $\overline{W_Y}$ respectively.
Note that $\phi_Z\circ\eta$ is not necessarily extendable to a definable $C^0$ function on $\overline Y$, and
$$\Ima(\phi_Y,\overline{\phi_Y\!f_{2Y}})=\Ima(\phi_Z,\overline{\phi_Z f_{2Z}})=\{(u,v)\in[0,\,1]^2:v\le u\}.$$
Set $F=(F_1,F_2)=(\phi_Z\circ\eta,\overline{(\phi_Z f_{2Z})\circ\eta})$ and $G=(G_1,G_2)=(\phi_Y,\overline{\phi_Y\!f_{2Y}})$.\vskip2mm\noindent
{\it Proof that (d) can be satisfied}.
The first condition in (d) is automatically satisfied as we have already mentioned.\par
For the second it suffices to find a definable isotopy $\mu_t,\ 0\le t\le1$, of $Y$ preserving $\{\sigma\cap Y:\sigma\in K_{2Y}\}$ such that $\phi_Y=\phi_Z\circ\eta\circ\mu_1$ on $Y\cap G^{-1}([0,\,\epsilon]\times[0,\,1])$ for some $\epsilon>0$.
First we reduce the problem to the case where $\phi_Z\circ\eta$ is extendable to a definable $C^0$ function on $\overline Y$, say, $\overline{\phi_Z\circ\eta}$.
Set $X=\graph\phi_Z\circ\eta\subset Y\times R$, and let $\phi_1$ and $\phi_2$ be the projection $\overline X\to R$ and the composite of the projection $\overline X\to\overline Y$ and $\phi_Y$ respectively.
Then only for construction of $\mu_t$ we can replace $\phi_Y:\overline Y\to R$ and $\phi_Z\circ\eta:Y\to R$ by $\phi_2$ and $\phi_1$ respectively.
Hence we can assume that $\phi_Z\circ\eta$ is extendable to a definable $C^0$ function on $\overline Y$.
By Triangulation theorem of definable $C^0$ functions we can regard $\phi_Y$ as a PL function on a compact polyhedron $\sigma\cap Y$ for each $\sigma$ in the original $K_{2Y}$ as the union of some open simplices in a new simplicial complex $K_{2Y}$.
Once more by Triangulation theorem of definable $C^0$ functions there is a definable isotopy $\nu_t,\ 0\le t\le1$, of $\overline Y$ preserving $K_{2Y}$ such that $\overline{\phi_Z\circ\eta}\circ\nu_1$ is PL.
Hence we assume that $\overline{\phi_Z\circ\eta}$ is PL from the beginning.
Then by Fact 2.5 we can suppose that $\phi_Y=\phi_Z\circ\eta$ on $Y\cap G^{-1}([0,\,\epsilon]\times[0,\,1])$ for some $\epsilon>0$.\par
The third condition is satisfied by Lemma 7, and (d) can be satisfied.\qed\vskip2mm\noindent
{\it Proof that (e) can be satisfied}.
We will find a PL homeomorphism $\delta$ of $\overline Y$ and a PL isotopy of $\overline Y$ whose finishing homeomorphism is $\delta$ such that $\tilde F_2\circ\delta=G_2$ on $G^{-1}(\overline N)$ (i.e., $\overline{(\phi_Z f_{2Z})\circ\eta}\circ\delta=\overline{\phi_Y\!f_{2Y}}$).
Then the equality $\phi_Y\circ\delta=\phi_Y$ should hold on $G^{-1}(\overline N)$.
If we consider $\tilde F_2$ and $G_2$ on $\phi^{-1}_Y(\epsilon')$ only for one $\epsilon'\in(0,\,1]$ and construct a PL homeomorphism $\delta_{\epsilon'}$ of $\phi^{-1}_Y(\epsilon')$ such that $\tilde F_2\circ\delta_{\epsilon'}=G_2$ on $G^{-1}(\overline N)\cap\phi^{-1}_Y(\epsilon')$, then the equality $\phi_Y\circ\delta_{\epsilon'}=\phi_Y$ is obvious.
We will choose $\epsilon'$ so close to 0 that such $\delta_{\epsilon'}$ is extendable to $\overline Y$ and the equality holds on $G^{-1}([0,\,\epsilon']\times[0,\,1])$.\par
Let $K_{3Y}$ be a simplicial subdivision of the $K_{2Y}$ in (c) such that $G|_{K_{3Y}}$ is a simplicial map onto some simplicial complex and the restriction of $\tilde F_2$ to each simplex in $K_{3Y}$ is linear (we cannot find $K_{3Y}$ such that both of $\tilde F_2|_{K_{3Y}}$ and $G|_{K_{3Y}}$ are simplicial maps onto some simplicial complexes, and a counterexample is given before Lemma 2.1 in \cite{S3}).
Let $\epsilon'\in R$ be positive and so small that $\phi_Y(K^0_{3Y})\cap(0,\,\epsilon']=\emptyset$, and let $K_{4Y}$ be the canonical simplicial subdivision of the cell complex $\{\sigma\cap\phi^{-1}_Y(\sigma_{\epsilon'}):\sigma\in K_{3Y},\,\sigma_{\epsilon'}\in\{0,\epsilon',1,[0,\,\epsilon'],[\epsilon',\,1]\}\}$.
Set $K_{4Y\epsilon'}=K_{4Y}|_{\phi^{-1}_Y(\epsilon')}$, and compare $\tilde F_2$ and $G_2$ on $\phi^{-1}_Y(\epsilon')$.
(Note that the $\{G(\sigma)\cap([0,\,\epsilon']\times[0,\,1]):\sigma\in K_{4Y}\}$ is a simplicial complex generated by the simplices $0*G(\sigma),\ \sigma\in K_{4Y\epsilon'}$.) 
Since both are PL and nonnegative and have the same zero set there exists a PL homeomorphism $\delta_{\epsilon'}$ (and an isotopy) of $\phi^{-1}_Y(\epsilon')$ preserving $K_{4Y\epsilon'}$ such that $\tilde F_2\circ\delta_{\epsilon'}=G_2$ on a neighborhood of $\phi^{-1}_Y(\epsilon')\cap\overline{W_Y}$ in $\phi^{-1}_Y(\epsilon')$, say, $G^{-1}(O)$ for a small segment $O$ with one end $(\epsilon',0)$ and the other in $\{\epsilon'\}\times[0,\,1]$.
We need to extend $\delta_{\epsilon'}$ to a PL homeomorphism $\delta$ of $\overline Y$ preserving $K_{4Y}$.
This is obvious on $\phi^{-1}_Y([\epsilon',\,1])$ by the Alexander trick.
We define $\delta$ on $\phi^{-1}_Y([0,\,\epsilon'])$ by 
\begin{align*}\delta(ry_1+(1-r)y_2)=ry_1+(1-r)\delta_{\epsilon'}(y_2)\quad\text{for}\ (y_1,y_2,r)\in\sigma_1\times\sigma_2\times[0,\,1],\\ \text{where}\ \sigma_1\in K_{4Y}|_{\phi^{-1}_Y(0)},\, \sigma_2\in K_{4Y\epsilon'}\ \text{with}\ \sigma_1*\sigma_2\in K_{4Y}.\end{align*}
(Any simplex in $K_{4Y}|_{\phi_Y^{-1}([0,\,\epsilon'])}$ is of form $\sigma_1*\sigma_2$ for such $\sigma_1$ and $\sigma_2$.) 
Then for the same $(y_1,y_2,r)$ with $G(y_2)\in O$, 
$$\tilde F_2\circ\delta(r y_1+(1-r)y_2)=\tilde F_2(ry_1+(1-r)\delta_{\epsilon'}(y_2))=r\tilde F_2(y_1)+(1-r)\tilde F_2\circ\delta_{\epsilon'}(y_2)$$
$$\quad\qquad\qquad\qquad\qquad=(1-r)\tilde F_2\circ\delta_{\epsilon'}(y_2)=(1-r)G_2(y_2)=G_2(r y_1+(1-r)y_2),$$
since $\tilde F_2$ and $G_2$ are linear on $\sigma_1*\sigma_2$.
Hence $\tilde F_2\circ\delta=G_2$ on $G^{-1}((0,0)*O)$.
Here the cone $(0,0)*O$ is a semilinear neighborhood of $(0,\,\epsilon']\times\{0\}$ in $[0,\,\epsilon']\times[0,\,1]$.
On the other hand, $\phi_Y\circ\delta(r y_1+(1-r)y_2)=(1-r)\epsilon'=\phi_Y(r y_1+(1-r)y_2)$ on $G^{-1}([0,\,\epsilon]\times[0,\,1])$.
Therefore, by shrinking $N$ to $((0,0)*O)\cap N$ and using the same notation $\epsilon$ and $N$, we can assume (e), i.e., that $\tilde F_2=G_2$ on $G^{-1}((0,0)*O)$ and hence $\tilde F_2=G_2$ on $G^{-1}(\overline N)$.\qed\vskip2mm\noindent
{\it Proof that (f) can be satisfied}.
Note that $f_{2Z}\circ\eta|_{Y\cap G^{-1}([\epsilon,\,1]\times[0,\,1])}$ is extended to a definable $C^0$ function on $G^{-1}([\epsilon,\,1]\times[0,\,1])$.
Let $g_Z$ denote the extension, and set $g_Y=f_{2Y}|_{G^{-1}([\epsilon,\,1]\times[0,\,1])}$.
Then $g_Z=g_Y$ on $G^{-1}(N\cap(\{\epsilon\}\times[0,\,1]))$, and it suffices to find a definable isotopy $\pi_t,\,0\le t\le1$, of $G^{-1}([\epsilon,\,1]\times[0,\,1])$ preserving $K_{4Y}|_{G^{-1}([\epsilon,\,1]\times[0,\,1])}$ such that $g_Z\circ\pi_1=g_Y$ on $G^{-1}([\epsilon,\,1]\times[0,\,\epsilon''])$ for some $\epsilon''>0$ because $\pi_t,\,0\le t\le1$, is extendable to a definable isotopy $\tilde\pi_t,\,0\le t\le1$, of $\overline Y$ preserving $K_{4Y}$ such that $\phi_Y\circ\tilde\pi_1=\phi_Y$ on $G^{-1}([0,\,\epsilon]\times[0,\,1])$ and $\overline{(\phi_Z f_{2Z})\circ\eta}\circ\tilde\pi_1=\overline{\phi_Y\!f_{2Y}}$ on $G^{-1}(\overline N)$ for the same reason as in the above proof.
Here only for construction of $\pi_Y$ we can assume that $g_Y$ and $g_Z$ are PL by Triangulation theorem of definable $C^0$ functions (1) and (2) and for the same reason as above.
Then existence of $\pi_t$ follows from Fact 2.5 because $g_Y$ and $g_Z$ satisfy the conditions in Fact 2.5.
Hence (f) can be satisfied.\qed\vskip2mm\noindent
{\it Proof of statement 9 in the non-locally-closed case}.
By the above argument we assume (a),...,(f).
We will find a definable isotopy $\alpha_t,\ 0\le t\le 1$, of $Y$ preserving $\{Y\cap\sigma:\sigma\in K_{1Y}\}$ such that $\eta\circ\alpha_1:(Y;Y_i)\to(Z;Z_i)$ is semilinear.\par
First of all, we see that $(\overline Y;Y,Y_i)$ and $(\overline Z;Z,Z_i)$ are definably homeomorphic.
Set $G_Z=(\phi_Z,\overline{\phi_Z f_{2Z}})$, which is the definable $C^0$ extension of $(\phi_Z,\phi_Z f_{2Z})$ to $\overline Z$, and let $L$ be the simplicial complex generated by the simplex $\{(u,v)\in[0,\,1]^2:v\le u\}$.
Let us remember that $G:K_{2Y}\to L$ and $G_Z:K_{2Z}\to L$ are simplicial.
Let $\psi_1$ and $\psi_2$ be definable nonnegative $C^0$ functions defined on $[0,\,1]$ with zero set $\{0\}$ such that $\psi_1(u)<u$ and $\psi_2(u)<u$ for $u\in(0,\,1]$, 
and set 
$$\Psi_j=\{(u,v)\in[0,\,1]^2:\psi_j(u)\le v\le u\},\quad j=1,2.$$
Note that $G^{-1}(\Psi_j)\subset Y$ and $G^{-1}_Z(\Psi_j)\subset Z$.
Later we define $\psi_1$ and $\psi_2$ more explicitly so that the $G^{-1}(\Psi_1)$ and $G^{-1}_Z(\Psi_1)$ are definably homeomorphic to $\overline Y$ and $\overline Z$, respectively, and $\eta|_{G(\Psi_2)}$ is a homeomorphism onto $G_Z^{-1}(\Psi_2)$.
Hence we need the following lemma, which we will prove later:\vskip2mm\noindent
{\bf Lemma 11.}
{\it There exist definable isotopies $\delta_{Y t},\ 0\le t\le 1$, of $\overline Y$ preserving $K_{2Y}$ and $\delta_{Z t},\ 0\le t\le 1$, of $\overline Z$ preserving $K_{2Z}$ such that $\delta_{Y1}(G^{-1}(\Psi_1))=G^{-1}(\Psi_2)$, $\delta_{Z1}(G_Z^{-1}(\Psi_1))=G_Z^{-1}(\Psi_2)$, $\delta_{Y t}=\id$ on a definable neighborhood of $W_Y$ in $\overline Y$ and $\delta_{Z t}=\id$ on a definable neighborhood of $W_Z$ in $\overline Z$.}\vskip2mm\noindent
{\it Continued proof of statement 9 in the non-locally-closed case.}
Assume that Lemma 11 is proved.
Choose $\psi_1$ and $\psi_2$ so that $\psi_1(u)=d u$ for some small positive $d\in R$, $\psi_2=\psi_1$ on $[\epsilon,\,1]$ and the set $\{(u,v)\in(0,\,1]^2:v\le\psi_2(u)\}$ is contained in $N\cup(\hat N\cap G_Z\circ\eta(Y\cap G^{-1}(\hat N)))$.
Then $G^{-1}(\Psi_1)=\overline{f_{2Y}([d,\,1])}$ and $G^{-1}_Z(\Psi_1)=\overline{f_{2Z}([d,\,1])}$ by definition of $G$ and $G_Z$; hence there are PL homeomorphisms $\kappa_Y:\overline Y\to G^{-1}(\Psi_1)$ and $\kappa_Z:\overline Z\to G^{-1}_Z(\Psi_1)$ by (a); moreover, $\eta(G^{-1}(\Psi_2))=G^{-1}_Z(\Psi_2)$ for the following reason:\par
Since the boundaries of $G^{-1}(\Psi_2)$ in $\overline Y$ and $G^{-1}_Z(\Psi_2)$ in $\overline Z$ are $\{y\in Y:\phi_Y\!f_{2Y}(y)=\psi_2\circ\phi_Y(y)\}$ and $\{z\in Z:\phi_Zf_{2Z}(z)=\psi_2\circ\phi_Z(z)\}$, respectively, it suffices to see 
$$\{y\in Y:\phi_Y f_{2Y}(y)=\psi_2\circ\phi_Y(y)\}=\{y\in Y:(\phi_Z f_{2Z})\circ\eta(y)=\psi_
2\circ\phi_Z\circ\eta(y)\}.\leqno{\text{(g)}}$$
We consider (g) separately on $\phi^{-1}_Y(0),\ \phi^{-1}_Y((0,\,\epsilon])$ and $\phi^{-1}_Y([\epsilon,\,1])$.
First the following is obvious:
\begin{align*}\{y\in\phi^{-1}_Y(0):\phi_Y f_{2Y}(y)=\psi_2\circ\phi_Y(y)\}=\phi^{-1}_Y(0)\\
=\{y\in\phi^{-1}_Y(0):(\phi_Z f_{2Z})\circ\eta(y)=\psi_2\circ\phi_Z\circ\eta(y)\}.\end{align*}
Secondly, we prove
\begin{align*}\{y\in\phi^{-1}_Y((0,\,\epsilon]):\phi_Y f_{2Y}(y)=\psi_2\circ\phi_Y(y)\}\\
=\{y\in\phi^{-1}_Y((0,\,\epsilon]):(\phi_Z f_{2Z})\circ\eta(y)=\psi_2\circ\phi_Z\circ\eta(y)\}.\end{align*}
Let $y$ be an element of the left-side set.
Then 
$$G(y)=(\phi_Y(y),\phi_Y\!f_{2Y}(y))=(\phi_Y(y),\psi_2\circ\phi_Y(y))\in N.$$
$$(\phi_Y(y),\phi_Y\!f_{2Y}(y))=(\phi_Z\circ\eta(y),(\phi_Z f_{2Z})\circ\eta(y))\quad\text{by (d) and (e).}\leqno{\text{Hence}}$$\vskip-7mm
$$(\phi_Z f_{2Z})\circ\eta(y)=\phi_Y\!f_{2Y}(y)=\psi_2\circ\phi_Y(y)=\psi_2\circ\phi_Z\circ\eta(y).\leqno{\text{Therefore}}$$
Namely, $y$ is an element of the right-side set.
In the same way we see that the right-side set is contained in the left-side.
Thus (g) holds on $\phi^{-1}_Y((0,\,\epsilon])$.\par
Finally, we see (g) on $\phi^{-1}_Y([\epsilon,\,1])$.
In the same way as above, by using (f) and the inclusions $\{(u,v)\in[\epsilon,\,1]\times(0,\,1]:v\le\psi_2(u)\})\subset\hat N$ and $G^{-1}_Z(\{(u,v)\in[\epsilon,\,1]\times(0,\,1]:v\le\psi_2(u)\})\subset\eta(Y\cap G^{-1}(\hat N))$, we obtain 
$$\{y\in\phi^{-1}_Y([\epsilon,\,1]):f_{2Y}(y)=d\}=\{y\in\phi^{-1}_Y([\epsilon,\,1]):f_{2Z}\circ\eta(y)=d\}.$$\vskip-7mm
$$\ \{y\in\phi^{-1}_Y([\epsilon,\,1]):\phi_Y\!f_{2Y}(y)=\psi_2\circ\phi_Y(y)\}\qquad\qquad\qquad\qquad\qquad\qquad\qquad\qquad\qquad\qquad\leqno{\text{Hence}}$$
$$=\{y\in\phi^{-1}_Y([\epsilon,\,1]):f_{2Y}(y)=d\}\quad\text{since}\ \psi_2(u)=d u\ \text{for}\ u\in[\epsilon,\,1]\qquad\qquad$$
$$=\{y\in\phi^{-1}_Y([\epsilon,\,1]):f_{2Z}\circ\eta(y)=d\}\qquad\qquad\qquad\qquad\qquad\qquad\qquad\qquad\quad$$
$$\qquad\qquad=\{y\in\phi^{-1}_Y([\epsilon,\,1]):(\phi_Z f_{2Z})\circ\eta(y)=\psi_2\circ\phi_Z\circ\eta(y)\}\quad\text{since}\ \psi_2(u)=d u\ \text{for}\ u\in[\epsilon,\,1].$$
Therefore, (g) holds on $\phi^{-1}_Y([\epsilon,\,1])$ and hence on $Y$.
Thus $\eta(G^{-1}(\Phi_2))=G^{-1}(\Psi_2)$.\par
Set $\omega_Y=\delta_{Y1}$ and $\omega_Z=\delta_{Z1}$.
We have obtained a sequence of definable homeomorphisms
$$\overline Y\stackrel{\kappa_Y}{\longrightarrow}G^{-1}(\Psi_1)\stackrel{\omega_Y}{\longrightarrow}G^{-1}(\Psi_2)\stackrel{\eta}{\longrightarrow}G^{-1}_Z(\Psi_2)\stackrel{\omega_Z^{-1}}{\longrightarrow}G^{-1}_Z(\Psi_1)\stackrel{\kappa_Z^{-1}}{\longrightarrow}\overline Z.$$
Moreover, the composite $\kappa^{-1}_Z\circ\cdots\circ\kappa_Y$ carries $\sigma\cap Y$ onto $\eta(\sigma\cap Y)$ for each $\sigma\in K_{1Y}$ by the condition in (a) that $\{\Int\sigma\cap Z:\sigma\in K_{1Z}\}$ is compatible with $\{\eta(\sigma\cap Y):\sigma\in K_{1Y}\}$.
Hence the composite is a definable homeomorphism from $(\overline Y;Y,Y_i)$ to $(\overline Z;Z,Z_i)$.
Consequently, $(\overline Y;Y,Y_j)$ and $(\overline Z;Z,Z_j)$ are definably homeomorphic.\par
Next we see that the restriction to $Y$ of the composite is the finishing homeomorphism of some isotopy of $\eta:Y\to Z$ through homeomorphisms.
For this we find isotopies whose finishing homeomorphisms are $\kappa_Y,\ \omega_Y,\ \kappa_Z$ and $\omega_Z$, and construct the required isotopy $\alpha_t,\ 0\le t\le1$, of $Y$ by the isotopies.
Let us define $\psi_{j t}$ and $\Psi_{j t}$ to be $t\psi_j$ and $\{(u,v)\in[0,\,1]^2:\psi_{j t}(u)\le v\le u\}$, respectively, for each $t\in[0,\,1]$ and $j=1,2$, and repeat the above argument for $\psi_{j t}$ and $\Psi_{j t}$ for $t\in(0,\,1]$.
Then by the above proof we have definable isotopies $\omega_{Y t},\ 0\le t\le 1$, of $\overline Y$ preserving $K_{2Y}$ and $\omega_{Z t},\ 0\le t\le 1$, of $\overline Z$ preserving $K_{2Z}$ such that $\omega_{Y1}=\omega_Y$, $\omega_{Z1}=\omega_Z$, $\omega_{Y t}(G^{-1}(\Psi_{1t}))=G^{-1}(\Psi_{2t})$ and $\omega_{Z t}(G^{-1}_Z(\Psi_{1t}))=G^{-1}_Z(\Psi_{2t})$ for $t\in[0,\,1]$.
On the other hand, by (a) there are definable isotopies $\kappa_{Y t}$ and $\kappa_{Z t},\ 0\le t\le1$, of the identity maps of $\overline Y$ and $\overline Z$ preserving $\{Y\cap\sigma:\sigma\in K_{1Y}\}$ and $\{Z\cap\sigma:\sigma\in K_{1Z}\}$, respectively, such that $\kappa_{Y1}=\kappa_Y$, $\kappa_{Z1}=\kappa_Z$, $\Ima\kappa_{Y t}=\overline{f_{2Y}([t d,\,1])}$ and $\Ima\kappa_{Z t}=\overline{f_{2Z}([t d,\,1])}$.
Set $\kappa'_{Y t}=\omega_{Y t}\circ\kappa_{Y t}$ and $\kappa'_{Z t}=\omega_{Z t}\circ\kappa_{Z t}$.
Then $\kappa'_{Y t}$ and $\kappa'_{Z t}$, $0\le t\le 1$, are definable isotopies of the identity maps of $\overline Y$ and $\overline Z$ preserving $\{Y\cap\sigma:\sigma\in K_{1Y}\}$ and $\{Z\cap\sigma:\sigma\in K_{1Z}\}$, respectively, such that $\Ima\kappa'_{Y t}=G^{-1}(\Psi_{2t})$ and $\Ima\kappa'_{Z t}=G^{-1}_Z(\Psi_{2t})$.
Note that $\Ima\eta\circ\kappa'_{Y t}=\Ima\kappa'_{Z t}$ for each $t\in(0,\,1]$.
Hence $\kappa^{\prime-1}_{Z t}\circ\eta\circ\kappa'_{Y t}|_Y:(Y;Y_i)\to(Z;Z_i)$, $0\le t\le1$, is a definable isotopy of $\eta$ to the map $\kappa_{Z1}^{\prime-1}\circ\eta\circ\kappa'_{Y1}|_Y\,(=\kappa^{-1}_Z\circ\omega^{-1}_Z\circ\eta\circ\omega_Y\circ\kappa_Y|_Y)$ through homeomorphisms and $\{\Int\sigma\cap Z:\sigma\in K_{1Z}\}$ is compatible with $\{\kappa^{\prime-1}_{Z t}\circ\eta\circ\kappa'_{Y t}(\sigma\cap Y):\sigma\in K_{1Y}\}$ by (a) for each $t$.\par
If $\kappa^{\prime-1}_{Z1}\circ\eta\circ\kappa'_{Y1}$ is PL, we define $\alpha_t$ to be $\eta^{-1}\kappa_{Z t}^{\prime-1}\circ\eta\circ\kappa'_{Y t}|_Y$.
Then $\alpha_t$ satisfies the conditions that $\alpha_t,\ 0\le t\le 1$, is a definable isotopy of $Y$ preserving $\{Y\cap\sigma:\sigma\in K_{1Y}\}$ and $\eta\circ\alpha_1$ is semi-linear.
Assume that $\kappa^{\prime-1}_{Z1}\circ\eta\circ\kappa'_{Y1}$ is not PL.
Then we modify it to a PL homeomorphism by using the property that $\kappa^{\prime-1}_{Z1}\circ\eta\circ\kappa'_{Y1}|_Y$ is the finishing homeomorphism of an isotopy through homeomorphisms.
By Uniqueness theorem of definable triangulation and its supplement, there is a definable isotopy $\mu_t:(\overline Y;Y\cap\sigma:\sigma\in K_{1Y})\to(\overline Z;\eta(Y\cap\sigma):\sigma\in K_{1Y}),\ 0\le t\le1$, of $\kappa^{\prime-1}_{Z1}\circ\eta\circ\kappa'_{Y1}$ through homeomorphisms such that $\mu_1$ is PL since $\kappa_{Z1}^{\prime-1}\circ\eta\circ\kappa'_{Y1}(Y\cap\sigma)=\eta(Y\cap\sigma)$ for $\sigma\in K_{1Y}$ and $\eta(Y\cap\sigma)$ are semilinear by (a).
We define $\alpha_t$ by
$$\alpha_t=\left \{
\begin{array}{l}
\eta^{-1}\circ\kappa^{\prime-1}_{Z2t}\circ\eta\circ\kappa'_{Y2t}|_Y\hspace{7mm}\text{for}\ t\in[0,\,1/2]\\
\eta^{-1}\circ\mu_{2t-1}|_Y\hspace{21mm}\text{for}\ t\in[1/2,\,1].
\end{array}
\right.$$
Then $\alpha_t$ satisfies the conditions for the same reason as above.
Thus statement 9 in the non-locally-closed case is proved.\qed\vskip2mm\noindent
{\it Proof of Lemma 11}.
We construct $\delta_{Yt}$ only.
Let $\psi_0$ and $\psi_3$ be functions having the same properties as $\psi_1$ and $\psi_2$ such that $\psi_0<\psi_j<\psi_3,\ j=1,2$, on $(0,\,1]$.
First define a definable isotopy $\delta_{L t},\ 0\le t\le 1$, of $|L|$ preserving $L$ so that for each $(u,v)\in|L|$, 
$$\delta_{L t}(u,v)=
\left \{
\begin{array}{l}
\!(u,v)\hspace{34mm}\text{for}\ v\in[0,\,\psi_0(u)]\cup[\psi_3(u),\,u]\\
\!\big(u,t\psi_2(u)+(1-t)\psi_1(u)\big)\hspace{3mm}\text{for}\ v\text{ equal to }\psi_1(u),
\end{array}
\right.$$
and for each fixed $u$, the function $\delta_{L t}(u,v)$ of the variable $v$ is linear on the segments $[\psi_0(u),\,\psi_1(u)]$ and $[\psi_1(u),\,\psi_3(u)]$.
Then $\delta_{L1}(\Psi_1)=\Psi_2$.
Next lift $\delta_{L t}$ to a definable isotopy $\delta_{Yt}$ of $\overline Y$ preserving $K_{2Y}$ by Lemma 6.
Then $\delta_{Y1}(G^{-1}(\Psi_1))=G^{-1}(\Psi_2)$ and $\delta_{Y t}=\id$ on $G^{-1}(\{v\le\psi_0(u)\})$---a definable neighborhood of $W_Y$ in $\overline Y$.
Thus Lemma 11 is proved.\qed
\section{Proofs of Theorems 1$'$ and 2$'$ and the corollary}
In this section we apply some results in \cite{S1} and \cite{S2}, which are stated only over $\R$ but hold over any $R$.
Assume that $0<r<\infty$ and if $R=\R$, $0<r\le\infty$.
A {\it definable $C^r$ function} on a definable set in $R^n$ is the restriction to the set of some definable $C^r$ function defined on some definable neighborhood of the set in $R^n$.
A {\it definable $C^r$ isotopy (homotopy)} $f_t:M\to N$, $0\le t\le1$, between definable $C^r$ manifolds is a definable isotopy such that the map $M\times[0,\,1]\ni(x,t)\to(f_t(x),t)\in N\times[0,\,1]$ is a $C^r$ embedding (respectively, map).
A definable $C^r$ {\it manifold with corners} over $R$ is a definable set in which each point has a neighborhood definably $C^r$ diffeomorphic to $R\times\cdots\times R\times[0,\,\infty)\times\cdots\times[0,\,\infty)$.
Hence a definable $C^r$ manifold with boundary is a definable $C^r$ manifold with corners.
A {\it definable $C^r$ cell} is a definable $C^r$ manifold with corners definably $C^r$ diffeomorphic to a cell.
Let $M$ or $M^m$ be a definable $C^r$ manifold possibly with corners of dimension $m$ over $R$.
We treat only a definable $C^r$ function on $M$ whose restriction to a proper face of $M$, if exists, is non-critical or constant.
Here a {\it face} of a definable $C^r$ manifold with corners is the closure of a definably connected component (i.e., a minimal definable open and closed subset) of a stratum of its canonical definable $C^r$ ($C^1$) stratification provided that the closure is a compact $C^r$ manifold possibly with corners (see p.\,189 in \cite{S2} for the definition of a canonical definable $C^1$ stratification).
For example, the faces of a simplex is ones in the usual sense.
Note that if $M$ is semialgebraic then its canonical definable $C^r$ stratification is semialgebraic and its definably connected component is semialgebraic.
Note also that $\partial M$ is not necessarily the union of its proper faces, e.g.~a 2-simplex which is smoothed at two points.
However, we treat such manifolds with corners only.
Moreover, we always assume that the closure of each stratum of the above stratification is a compact $C^r$ manifold possibly with corners.\par
A definable $C^r$ function on $M$ is called {\it Morse} if it is of the form $\sum_{i=1}^k x^2_i-\sum_{i=k+1}^m x^2_i+\rm{const}$ locally at each critical point in $\Int M$ for some definable local coordinate system $(x_1,\ldots,x_m)$.
Such a critical point or the function at the point are called {\it of type} $k$.
Given a definable $C^0$ function $f$ on $M$, a definable $C^0$ function $g$ on $M$ is called a $C^0$ $\phi$-{\it approximation} or simply a $C^0$ {\it approximation} of $f$ if $|f-g|<\phi$ for a small positive definable $C^0$ function $\phi$ on $M$.
Note that if $M$ is compact then a $C^0$ approximation $g$ of $f$ is defined by $|f-g|<c$ for a small $c>0\in R$.
We define a $C^1$ {\it approximation} of a definable $C^1$ map between definable $C^1$ manifolds by replacing the domain by a finite system of definable coordinate neighborhoods and the target space by its ambient Euclidean space and by also considering the first derivatives of the map.
When $M$ is compact and we choose $c$ in $\R$, we define a $C^0\,(C^1)$ {\it approximation in the $(C^1)\ \R$-topology}.
Note that given a definable $C^0$ or $C^1$ function $f$ on $\{x\in R^n:|x|\le1\}$ whose image is contained in $\{x\in R:|x|\le1\}$, the polynomial function theorem and the simplicial approximation theorem hold for $f$ in the $\R$-topology.\par
A special case of Theorem II.6.5 and Corollary II.6.6 in \cite{S2} is the next lemma in the case $\partial M=\emptyset$.
The case $\partial M\not=\emptyset$ is obviously proved in the same way.
Hence we do not prove it.\vskip2mm\noindent
{\bf Lemma 12} ($C^r$ trivialization of definable function).
{\it Let $f$ be a definable $C^r$ function on the above $M$, $0<r<\infty$.
Assume that $\Ima f$ is a closed, open, or half-open interval, say $I$, $f|_{\Int M}$ is non-critical and $f:M\to I$ is proper.
Then there exists a definable $C^r$ diffeomorphism $\pi:f^{-1}(\epsilon)\times I\to M$ for any $\epsilon\in I$ such that $f\circ\pi:f^{-1}(\epsilon)\times I\to I$ is the projection.
Moreover, if a definable $C^r$ diffeomorphism $\pi:(f^{-1}(\epsilon)\cap\partial M)\times I\to\cup_iF_i$ is given so that $f\circ\pi$ is the projection onto $I$ from the beginning, where $F_i$ are the proper faces of $M$ where $f$ is non-critical, then we can extend $\pi$ to the above global $\pi$.}\vskip2mm
Moreover, Theorem II.6.5 and Corollary II.6.6 in \cite{S2} states the same conclusion as above for a compact contractible definable $C^r$ manifold $I'$ with corners and a proper definable surjective $C^r$ submersion $f':M\to I'$.
In the proof of Theorem II.6.7 in \cite{S2} we proved and used the following:\vskip2mm\noindent
{\bf Lemma 13} (Morse approximation).
{\it Let $f$ be a definable $C^0$ function on the $M$ such that $f|_{\partial M}$ satisfies the conditions on definable $C^r$ functions, i.e., $f|_{\partial M}$ is of class $C^r$ and its restriction to each proper face of $M$ is non-critical or constant.
Then $f$ is $C^1$ approximated by a definable Morse $C^r$ function $g$ such that $g=f$ on $\partial M$ and $g(x_1)\not=g(x_2)$ for any distinct critical points $x_1$ and $x_2$ in $\Int M$ of $g$, $0<r<\infty$.}\vskip2mm
Let a function $\psi$ on $R^m$ be defined by $\psi(x)=\sum_{i=1}^k x^2_i-\sum_{i=k+1}^m x^2_i$ for $x=(x_1,\ldots,x_m)\in R^m$, where $0<k<m$.
Set $N=\{x\in R^m:|\psi(x)|\le1,|x|\le l\}$ for sufficiently large $l\in\N$.
Then $N$ is a compact Nash manifold with corners and there are unions of proper faces $N_{-1},\,N_0$ and $N_1$ of $N$ such that they are Nash manifolds with boundary, $\partial N_0=N_0\cap(N_{-1}\cup N_1)$, $\partial N=N_{-1}\cup N_0\cup N_1$, $\psi(N_{-1})=-1,\ \psi(N_1)=1$, $N_{-1}$ is Nash diffeomorphic to $D^k\times\partial D^{m-k}$ and $\psi|_{\Int N_0}$ is non-critical, where $D^i$ denotes the unit $i$-disk $\{x\in R^i:|x|\le1\}$.\vskip2mm\noindent
{\bf Remark 14} (Morse function).
{\it Assume that $M^m$ is a compact definable $C^r$ manifold with boundary, $0<r<\infty$, and let $g$ be a definable Morse $C^r$ function on $M$ such that there is only one critical point, say $x$, of $g|_{\Int M}$, $x$ is of type $k$ with $0<k<m$, $g(\partial M)=\{-2,2\}$ and $g(x)=0$.
Then we obtain $M$ from $g^{-1}([-2,\,-1])$ attaching an $(m-k)$-handle of Smale and smoothing the corners (see \cite{Sm}), which is precisely stated as follows:\par
By the definition of a definable $C^r$ manifold with boundary, there is a definable $C^r$ embedding $\tau:N\to M$ such that $g\circ\tau=\epsilon\psi$ for some $\epsilon>0\in R$.
Here we assume $\epsilon=1$ for simplicity of notation, which is clearly possible by the next argument.
By Lemma 12 we have a semialgebraic $C^r$ diffeomorphism $\pi_0:(N_0\cap\psi^{-1}(0))\times[-1,\,1]\to N_0$ such that $\psi\circ\pi_0$ is the projection onto $[-1,\,1]$.
Once more by Lemma 12, the diffeomorphism $\tau(N_0\cap\psi^{-1}(0))\times[-1,\,1]\ni(y,t)\to\tau\circ\pi_0(\tau^{-1}(y),t)\in\tau(N_0)$ is extended to a definable $C^r$ diffeomorphism $\pi:\big(g^{-1}(0)\cap(M-\pi(N)\overline)\big)\times[-1,\,1]\to(g^{-1}([-1,\,1])-\pi(N)\overline)$ so that $g\circ\pi$ is the projection onto $[-1,\,1]$.
A definable $C^r$ diffeomorphism $\pi_1:g^{-1}(1)\times[1,\,2]\to g^{-1}([1,\,2])$ also exists so that $g\circ\pi_1$ is the projection onto $[1,\,2]$.
We regard $\tau(N)\cup\pi_1(\tau(N_1)\times[1,\,2])$ as a handle.
Thus we obtain $M$ from $g^{-1}([-2,-1])$ attaching this handle and smoothing the corners.}
We call the union of $g^{-1}([-2,-1])$ and the handle a {\it manifold with reflex angles on the boundary}, which is not a manifold with corners but whose {\it corners} are naturally defined.\par
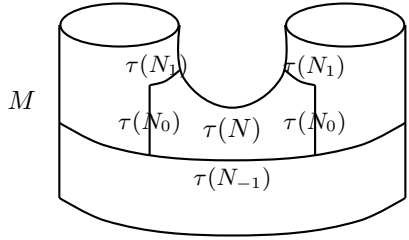
\begin{figure}[h]
\begin{pspicture}(3.7,-0.5)(12.7,2.6)
\psellipse(10,2.3)(0.8,0.3)
\psellipse(13,2.3)(0.8,0.3)
\psline(9.2,2.3)(9.2,0)
\psline(13.8,2.3)(13.8,0)
\pscurve(9.2,0)(9.8,-0.3)(11.5,-0.5)(13.2,-0.3)(13.8,0)
\pscurve(9.2,1)(9.8,0.7)(11.5,0.5)(13.2,0.7)(13.8,1)
\pscurve(10.8,2.3)(10.85,1.6)(11.5,1.2)(12.15,1.6)(12.2,2.3)
\pscurve(10.8,1.7)(10.6,1.55)(10.4,1.5)
\psline(10.4,1.5)(10.4,0.57)
\pscurve(12.2,1.7)(12.4,1.55)(12.6,1.5)
\psline(12.6,1.5)(12.6,0.57)
\rput(8.7,1.3){$M$}
\rput(11.5,0.9){$\tau(N)$}
\rput(10.5,1.75){\small$\tau(N_1)$}
\rput(12.59,1.74){\small$\tau(N_1)$}
\rput(12.6,1){\small$\tau(N_0)$}
\rput(10.4,1){\small$\tau(N_0)$}
\rput(11.5,0.25){\small$\tau(N_{-1})$}
\end{pspicture}
\caption{Handle-body}\end{figure}\par
We use handle-bodies.
The idea of application of handle-bodies is the following:
Let $M,\ g$ and $x$ be as in Remark 14.
Given a problem on $M$, we assume that it is already proved for smaller dimensional $M$.
Then we solve it at $\tau(N)$ by shrinking $N$ and outside $\tau(N)$ by reducing the problem to the case of dimension $m-1$ by the triviality in Remark 14.
\vskip2mm\noindent
{\it Proof of the first statement of Theorem 2\,$'$.}
We postpone the proof of the uniqueness.
We proceed as follows:
First we clarify the problem and state statement 1.
Next we reduce statement 1 to statement 2.
Thirdly, we prove that statement 2 follows from statement 3.
Fourthly, we simplify successively statement 3 to statement 3$'$, statement 3$''$ and statement 4.
Finally, we prove statement 4.\par
Let $M$ be a definable $C^r$ manifold over $R$, $0<r<\infty$.
We can assume $r=1$ because if $M$ is definably $C^1$ diffeomorphic to another definable $C^r$ manifold then $M$ is $C^r$ definably diffeomorphic to the manifold (Theorem II.5.2 in \cite{S2}).
If the o-minimal structure is that of semialgebraic sets then $M$ is semialgebraically $C^1$ diffeomorphic to a Nash manifold (Theorem III.1.3 in \cite{S1}).
Moreover, by Corollary 3.9 in \cite{C-S1}, any Nash manifold is Nash diffeomorphic to some nonsingular algebraic variety defined by polynomials with coefficients in $\R_{\rm alg}$.
Hence it suffices to prove the following statement:\vskip1mm\noindent
{\it Statement 1. The $M^m$ is definably $C^1$ diffeomorphic to a semialgebraic $C^1$ manifold.}\par
We prove statement 1 under statement 2 below by double induction, first on $\dim M$, using a handlebody decomposition.
We assume that $M$ is bounded in some $R^n$.
Let $g$ be a positive definable $C^0$ function on $M$ such that $g(x)\to0$ as a sequence of points $x$ in $M$ converges to a point of $\overline M-M$, e.g.~the function measuring the distance from $\overline M-M$.
Then by Lemma 13 we can assume that $g$ is a definable Morse $C^1$ function such that $g(x_1)\not=g(x_2)$ for any distinct critical points $x_1$ and $x_2$ of $g$.
For simplicity of notation, let $x_1,\ldots,x_p$ be the critical points of $g$ such that $g(x_i)=2i$ for $i=1,\ldots,p$.
Then $p\ge1$ because $x_p$ is the point where $g$ takes the maximum value.
Set $M_i=g^{-1}((0,2i-1]),\ i=0,\ldots,p+1$.
Then $M_0=\emptyset,\ M_{p+1}=M$ and $M_1,\ldots,M_p$ are definable $C^1$ manifolds with boundary.
Let $i\in\N$ with $i\le p$.
As the second induction hypothesis we assume that there is a semialgebraic $C^1$ manifold $S_i$ possibly with boundary and a definable $C^1$ diffeomorphism $\pi_i:S_i\to M_i$.
Then we need to define a semialgebraic $C^1$ manifold $S_{i+1}$ with boundary and a definable $C^1$ diffeomorphism $\pi_{i+1}:S_{i+1}\to M_{i+1}$.\par
Case $i=0$.
Assume that $M$ is not compact.
Then there is no critical points of $g$ in $M_1$ and $g|_{M_1}:M_1\to(0,\,1]$ is proper.
Hence by Lemma 12 there is a definable $C^1$ diffeomorphism from $g^{-1}(1)\times(0,\,1]$ to $M_1$.
By the first induction hypothesis, $g^{-1}(1)$ is definably $C^1$ diffeomorphic to a semialgebraic $C^1$ manifold because $g^{-1}(1)$ is a definable $C^1$ manifold of dimension $=m-1$.
Therefore, $M_1$ is definably $C^1$ diffeomorphic to a semialgebraic $C^1$ manifold with boundary.
In the compact case, $M_1=\emptyset$ and hence there is nothing to prove.\par
Case $0<i\le p$.
Let $x_i$ be of type $k$.
Then $0\le k\le m$.
There are three cases to consider: $k=0,\ k=m$ or $0<k<m$.\par
If $k=0$, we obtain $M_{i+1}$ from $M_i$ attaching an $m$-handle $H$, and $H\cap M_i$ is a definably connected component of $\partial M_i$.
Hence $M_{i+1}$ is definably $C^1$ diffeomorphic to a definable $C^1$ manifold $S_i\cup_\tau N$ with boundary where $N$ is a semialgebraic $C^1$ $m$-handle $D^m$ and attached to $S_i$ by a definable $C^1$ diffeomorphism $\tau$ from $N_{-1}=\partial D^m$, onto a semialgebraically connected component of $\partial S_i$.
If $\tau$ is semialgebraic, so is $S_i\cup_\tau N$.
Hence we need to modify $\tau$ to be semialgebraic.
Consider the double of $S_i$.
Then we have a semialgebraic $C^1$ manifold $S_i\cup{\tau_0}(\partial N\times[0,\,1])$ with boundary which is semialgebraically $C^1$ diffeomorphic to $S_i$, where $\tau_0$ is a semialgebraic $C^1$ diffeomorphism from $\partial N\times\{0\}$ onto the semialgebraically connected component of $\partial S_i$.
By a semialgebraic $C^1$ collar of $\partial N$ in $N$, we regard $N$ as $N'\cup_{\tau_1}(\partial N\times[0,\,1])$, where $N'$ is the closure of the complement of the collar in $N$ and $\tau_1$ is a semialgebraic $C^1$ diffeomorphism from $\partial N\times\{1\}$ onto $\partial N'$.
Then $S_i\cup_\tau N=(S_i\cup_{\tau_0}(\partial N\times[0,\,1]))\cup_{\tau^{-1}_1} N'$.
Hence we modify $\tau^{-1}_1$ to be semialgebraic.
By statement 2 below there exists a definable $C^1$ isotopy $\pi_t:\partial N'\to\partial N\times\{1\},\ 0\le t\le1$, such that $\pi_0=\tau^{-1}_1$ and $\pi_1$ is semialgebraic.
Here $\pi_t$ for each $t$ is a diffeomorphism because it is an embedding, $\pi_0$ is a diffeomorphism and $\partial N'$ is a compact definable $C^0$ manifold.
Set $S_{i+1}=(S_i\cup_{\tau_0}(\partial N\times[0,\,1]))\cup_{\pi_1}N'$.
Then $S_{i+1}$ is a semialgebraic $C^1$ manifold with boundary and there is clearly a definable $C^1$ diffeomorphism from $S_{i+1}$ onto $M_{i+1}$.\par
Assume $0<k<m$.
By Remark 14 and the second induction hypothesis, we obtain $M_{i+1}$ from $S_i$ attaching a semialgebraic $C^1$ $(m-k)$-handle $N$ and smoothing the corners.
Then $N$ is attached to $S_i$ by a definable $C^1$ embedding $\tau:D^k\times\partial D^{m-k}\to\partial S_i$.
Here we can assume that $\tau$ is semialgebraic for the same reason as in the case $0<i\le p$ and $k=0$.
Next we smooth the ``corners" of the definable $C^1$ manifold reflex angles on the boundary and of the semialgebraic $C^1$ manifold reflex angles on the boundary.
By the definition of smoothing, which works also in the $R$ case, we can smooth the manifolds reflex angles on the boundary so that the definable $C^1$ diffeomorphism between the manifolds reflex angles on the boundary induces a definable $C^1$ diffeomorphism between the smoothed ones.\par
If $k=m$, $M_{i+1}$ is regarded as the disjoint union of $M_i$ and a semialgebraic $C^1$ $0$-handle.
Hence $M_{i+1}$ is definably $C^1$ diffeomorphic to a semialgebraic $C^1$ manifold with boundary for the same reason.
Thus statement 1 proved, provided that statement 2 holds.\vskip1mm\noindent
{\it Statement 2.
Let $N_{-1}$ denote $D^k\times\partial D^{m-k}$ for $0\le k<m$, and $\tau:N_{-1}\to\partial S_i$ be a definable $C^1$ embedding.
Then there exists a definable $C^1$ isotopy $\tau_t:N_{-1}\to\partial S_i,\ 0\le t\le1$, such that $\tau_0=\tau$ and $\tau_1$ is semialgebraic.}\par
We will see that statement 2 follows from the following:\vskip1mm\noindent
{\it Statement 3.
Let $N$ be a compact semialgebraic $C^1$ manifold possibly with corners of dimension $m>0$ in $R^{m'}$, and let $\{G\}$ be a family of proper faces of $N$.
Let $U$ be a small open semialgebraic neighborhood of a closed semialgebraic subset $X$ of $N$ such that $\dim X<m$ and $\dim X\cap\partial N<m-1$.
Let $\tau:N\to R^m$ be a definable $C^1$ embedding such that $\tau$ is semialgebraic outside $U$.
Then there exists a definable $C^1$ isotopy $\tau_t:N\to R^m$, $0\le t\le1$, such that $\tau_0=\tau$, $\tau_t=\tau$ on an arbitrarily small semialgebraic neighborhood of $\cup G$ and $\tau_1$ is semialgebraic outside an arbitrarily small semialgebraic neighborhood of $\overline U\cap\cup G$.}\par
Before proving statement 3 and that statement 3 implies statement 2, we define a semialgebraic $C^1$ division of $N_{-1}$ by induction on $\dim N_{-1}$.
It is a finite family $\{M_l\}$ of compact semialgebraic $C^1$ manifolds with corners of dimension $m-1$ in $N_{-1}$ such that $\cup M_l=N_{-1}$, $M_l\cap M_{l'}\subset\partial M_l$ for $l\not=l'$, and for a proper face $F$ of $M_l$, the family $\{F\cap M_{l'}:l'\not=l\}$ is a semialgebraic $C^1$ division of $F$.
It follows that $\{\Int M_l\}$ is a semialgebraic $C^1$ stratification of $N_{-1}$ with the frontier condition.\par
We will use a special semialgebraic $C^1$ division of $N_{-1}$.
For this we give a proof to the $C^1$ triangulation theorem of a semialgebraic $C^1$ manifold of Cairns-Whitehead, which states that any semialgebraic $C^1$ manifold possibly with corners is semialgebraically $C^1$ triangulable, and whose well-known proof is false in the $R$-case.
We need only its correct proof but not the theorem itself.
The proof of Triangulation theorem of definable sets (see the proof of Theorem II.2.1 in \cite{S2}) says that there are a simplicial complex $K$ and a semialgebraic homeomorphism $\tau:|K|\to N_{-1}$ such that $\tau|_{\Int\sigma}$ but not yet $\tau|_\sigma$ is a $C^1$ embedding for each $\sigma\in K$ and $\{\tau(\Int\sigma):\sigma\in K\}$ is a refinement of any given finite semialgebraic open covering of $N_{-1}$.
Note that the $C^1$ triangulation theorem of a semialgebraic $C^1$ manifold says that $\tau|_\sigma$ can be a $C^1$ embedding.
For each pair $\sigma'\subset\sigma$ in $K$ with $\dim\sigma=m-1$ and $\dim\sigma'=m-2$, remove from $\tau(\Int\sigma')$ the points where $(\tau(\Int\sigma),\tau(\Int\sigma'))$ does not satisfy the Whitney condition.
(See p.\,4 in \cite{S2} for the definition the Whitney condition.
Lemma I.2.2 and its proof in \cite{S2} says that the set of the points is semialgebraic and of dimension $<m-2$.)
Choose a semialgebraic triangulation of $\cup_{\sigma''\in K,\dim\sigma''\le m-2}\tau(\sigma'')$ compatible with the removed sets, and repeat the same argument.
Then we obtain a Whitney semialgebraic $C^1$ stratification $\{X_j\}$ of $N_{-1}$ and for each $X_j$, a semialgebraic homeomorphism $\tau_j$ from a simplex $\sigma_j$ to $\overline{X_j}$ such that for each $j$, $\{\tau^{-1}_j(X_{j'}):j'\not=j\}$ is a Whitney semialgebraic $C^1$ stratification of $\sigma_j$.
Moreover, we can assume that $\{\graph\tau_j|_{\overline{X_j}\cap X_{j'}}:j'\}$ is a Whitney semialgebraic $C^1$ stratification for the same reason.
Then the stratified map $\tau_j:\{\tau^{-1}_j(X_{j'}):j'\}\to\{X_{j'}:X_{j'}\subset\overline{X_j}\}$ is a strong isomorphism in the sense explained on p.\,26 in \cite{S2}.
We do not give its definition because the property, which we will use, is only Proposition I.1.13 in \cite{S2}.
For each $j$, let $\rho_j$ denote the function on $N_{-1}$ measuring the distance from $X_j$, let $\epsilon_j\in R$ such that $0<\epsilon_j\ll\epsilon_{j'}\ll a$ for some small $a>0\in R$ and for $j'$ with $X_{j'}\subset\overline{X_j}-X_j$, and set $\epsilon=\{\epsilon_j\}$.
For $j$ and $j'$ with $X_{j'}\subset\overline{X_j}-X_j$, set
$$X^\epsilon_{j,j'}=\{x\in\overline{X_j}:\dis(x,X_{j''})\ge\epsilon_{j''}\ \text{for}\ j''\ \text{with }X_{j''}\subset\overline{X_{j'}}\},\quad X^\epsilon_j=\cap_{j'}X^\epsilon_{j,j'}.$$
Then $X^\epsilon_{j,j'}$ and $X^\epsilon_j$ are compact semialgebraic $C^1$ cells for the following reason:\par
Proposition I.1.13 says that $X^\epsilon_{j,j'}$ and $X^\epsilon_j$ are semialgebraically $C^1$ diffeomorphic to $\sigma^\epsilon_{j,j'}$ and $\sigma^\epsilon_j$, which are defined in the same way for $\sigma_j,\ \tau^{-1}_j(\sigma_{j'})$ and the semialgebraic $C^0$ function $\rho_{j,j'}$ on $\sigma_j$.
Hence we can clarify what we prove as follows:\vskip1mm\noindent
{\it Let $\sigma$ denote one of the above $\sigma_j$'s, $\{Y_j\}$ the Whitney semialgebraic $C^1$ stratification $\{\tau^{-1}_j(X_{j'}):j'\}$ of $\sigma$.
Define $\sigma^\epsilon_j$ and $\sigma^\epsilon$ by the functions on $\sigma$ measuring the distances from $Y_j$ as above.
Then they are semialgebraic $C^1$ cells.}\par
We prove this statement by double induction, first on $\dim\sigma$, and secondly on $\dim Y_j$.
If $\dim Y_j=0$ then this is obvious.
Let $Y_j$ be of dimension 1.
Then $Y_j\cap\cap_{\dim Y_{j'}=0,Y_{j'}\subset\overline{Y_j}}\sigma^\epsilon_{j'}$ is semialgebraically $C^1$ diffeomorphic to a segment.
Hence we easy see that $(\cap_{\dim Y_{j'}=0,Y_{j'}\subset\overline{Y_j}}\sigma^\epsilon_{j'},Y_j\cap\cap_{\dim Y_{j'}=0,Y_{j'}\subset\overline{Y_j}}\sigma^\epsilon_{j'})$ is semialgebraically $C^1$ diffeomorphic to the pair of a cell and a segment contained in some proper face $F$ of the cell such that the segment intersects transversally with any other proper face of $F$.
Therefore, $\sigma^\epsilon_j$ is a semialgebraic $C^1$ cell.
For general $j$, $Y_j\cap\cap_{\dim Y_{j'}<\dim Y_j,Y_{j'}\subset\overline{Y_j}}\sigma^\epsilon_{j'}$ is a semialgebraic $C^1$ cell by the first induction hypothesis.
Then in the same way we show that all $\sigma^\epsilon_j$ and $\sigma^\epsilon$ are semialgebraic $C^1$ cells by the second induction hypothesis.\par
Here we interrupt the proof of the $C^1$ triangulation theorem of a semialgebraic $C^1$ manifold, although we complete it very soon, and we come back to the main stream of the proof of statement 3 and that statement 3 implies statement 2.\par
Let $\{X_j\}$ be such that each $X_j$ is the union of strata of dimension $j$ in the above Whitney semialgebraic $C^1$ stratification of $N_{-1}$, which is a Whitney semialgebraic $C^1$ stratification of $N_{-1}$ but each $X_j$ is not necessarily semialgebraically connected or the frontier condition is not necessarily satisfied.
Let $\epsilon_0,\ldots,\epsilon_{m-2}\in R$ such that $0<\epsilon_{m-2}\ll\cdots\ll\epsilon_0\ll a$ for some small $a>0\in R$.
Set $\epsilon=(\epsilon_0,\ldots,\epsilon_{m-2})$ and
$$M^\epsilon_j=\{x\in N_{-1}:\dis(x,X_j)\le\epsilon_j,\,\dis(x,X_{j'})\ge\epsilon_{j'},\,j'=0,\ldots,j-1\}.$$
Then we have the following properties:\par
The $M^\epsilon_j$ is a compact semialgebraic $C^1$ manifold possibly with corners.
Its semialgebraically connected component is semialgebraically $C^1$ diffeomorphic to the product of a semialgebraic $C^1$ cell and a set $\{x\in R\times\cdots\times R\times[0,\,\infty)\times\cdots\times[0,\,\infty):|x|\le1\}$ by Theorem II.6.5 and Corollary II.6.6 in \cite{S2}.
Here the semialgebraic $C^1$ cell and the set are corresponding to a semialgebraically connected component of $\{x\in X_j:\dis(x,X_{j'})\ge\epsilon_{j'},\,j'=0,\ldots,j-1\}$ and $\{x\in N_{-1}\cap T:|x-x_0|\le\epsilon_j\}$ by the semialgebraic $C^1$ diffeomorphism for some $x_0\in X_j$ and for a linear space $T$ transversal to $X_j$ in $R^m$ and with $T\cap X_j\ni x_0$; the diffeomorphism carries the intersection of $\partial N_{-1}$ and the semialgebraically connected component of $M^\epsilon_j$ to the product of the semialgebraic $C^1$ cell and the set $\{x\in R\times\cdots\times R\times\cdots\times\partial([0,\,\infty)\times\cdots\times[0,\,\infty)):|x|\le1\}$; if the semialgebraically connected component of $X_j$ containing $x_0$ is not contained in $\partial N_{-1}$ then the semialgebraically connected component of $M^\epsilon_j$ is carried to the product of the cell and the unit disk $\{x\in R^{m-1-j}:|x|\le1\}$.
Hence each semialgebraically connected component of $M^\epsilon_j$ is a semialgebraic $C^1$ cell and $\cup_{j_1\not=j}(M^\epsilon_j\cap M^\epsilon_{j_1})=\overline{\partial M^\epsilon_j-\partial N_{-1}}$.
Moreover, $M^\epsilon_{j_1}\cap\cdots\cap M^\epsilon_{j_l}$ and its intersection with $\partial N_{-1}$ have the same properties for any $j_1,\ldots,j_l$.
Therefore, the set of semialgebraically connected components of all $M^\epsilon_{j_1}\cap\cdots\cap M^\epsilon_{j_l}$ and $M^\epsilon_{j_1}\cap\cdots\cap M^\epsilon_{j_l}\cap\partial N_{-1}$ is a semialgebraic $C^1$ subdivision of $N_{-1}$.
Then we can regard it as a cell complex.
Hence, subdividing the cell complex by downward induction on $j$ in the same way as a stellar subdivision, we prove the $C^1$ triangulation theorem of a semialgebraic $C^1$ manifold (see p.\,15 in \cite{R-S} for a stellar subdivision).\par
The set $(M^{2\epsilon}_j\cap\cup_{j_1}\partial M^\epsilon_{j_1}-\partial N_{-1}\overline)$ is of dimension $m-2$.
Each semialgebraically connected component of this set is carried by the semialgebraic $C^1$ diffeomorphism, defined for $2\epsilon$, to the product of the semialgebraic $C^1$ cell and the following subset $Y$ of $\{x\in N_{-1}\cap T:\epsilon_j\le|x-x_0|\le2\epsilon_j\}$:
$$Y=\{x\in N_{-1}\cap T:|x-x_0|=\epsilon_j\text{ or }\epsilon_j<|x-x_0|\le2\epsilon_j\text{ and }x\in\cup_{j_1>j}\partial M^\epsilon_{j_1}\}.$$
Then $Y\cap\{x\in N_{-1}\cap T:\epsilon_j<|x-x_0|\le2\epsilon_j\}$ is semialgebraically $C^1$ diffeomorphic to the complement of the cone with vertex $x_0$ and base $\{x\in N_{-1}\cap T:|x-x_0|=\epsilon_j,\,x\in\cup_{j_1>j}\partial M^\epsilon_{j_1}\}$ in the cone with the same vertex and the base defined by $2\epsilon_j$ in place of $\epsilon_j$.
\vskip1mm\noindent
{\it Proof that statement 3 implies statement 2.}
Assume that statement 3 holds.
To apply statement 3 we need to reduce statement 2 to the case where $\tau$ is a definable $C^1$ embedding into $R^{m-1}$.
Let $R^n$ be the ambient Euclidean space of $\partial S_i$.
Then there are a finite number of linear maps $q_j:R^n\to R^{m-1}$ such that for each $x\in\partial S_i$, the germ of $q_j|_{\partial S_i}$ at $x$ is a $C^1$ embedding germ for some $j$.
Hence there is a finite semialgebraic open covering $\{O_l\}$ of $\partial S_i$ such that for each $l$, $q_j|_{O_l}:O_l\to R^{m-1}$ is a $C^1$ embedding for some $j$.
Here we refine $\{O_l\}$ so that each $O_l$ is definably $C^1$ diffeomorphic to $R^{m-1}$, and we modify $q_j$ so that $q_j|_{O_l}$ is a semialgebraic $C^1$ diffeomorphism from $O_l$ to $R^{m-1}$.
Consider a definable open covering $\{\tau^{-1}(O_l)\}$ of $N_{-1}$.
Then there is the above Whitney semialgebraic $C^1$ stratification $\{X_j\}$ such that the family of the semialgebraically connected components of all $X_j$ is a refinement of $\{\tau^{-1}(O_l)\}$.
Let the symbol $\{M^\epsilon_j\}$ now mean the family of semialgebraically connected components of all the above $M^\epsilon_j$ defined by $\{X_j\}$.
Then each $M^\epsilon_j$ is contained in $\tau^{-1}(O_l)$ for some $l$, say $l_j$.
We write $q_{l_j}\circ\tau|_{M^\epsilon_j}:M^\epsilon_j\to R^{m-1}$ simply as $\tau_j:M^\epsilon_j\to R^{m-1}$.
Thus we replace $\tau:N_{-1}\to\partial S_i$ by the family of all $\tau_j:M^\epsilon_j\to R^{m-1}$.\par
For each $j$, apply statement 3 to $\tau_j:M^\epsilon_j\to R^{m-1}$, $X=U=\emptyset$ and the family $\{G\}$ of the proper faces of $M^\epsilon_j$ which are not contained in $\partial N_{-1}$.
Then there exists a definable $C^1$ isotopy $\tau_{j,t}:M^\epsilon_j\to R^{m-1},\ 0\le t\le1$, satisfying the conditions in statement 3.
Note that $q^{-1}_{l_j}\circ\tau_{j,t},\ 0\le t\le1$, is a definable $C^1$ isotopy from $M^\epsilon_j$ into $\partial S_i$ and extended naturally to a definable $C^1$ isotopy from $N_{-1}$ into $\partial S_i$, which is defined to be $\tau$ outside $M^\epsilon_j$.
Hence by setting $\tau_t=q^{-1}_{l_j}\circ\tau_{j,t}$ on each $M^\epsilon_j$, we obtain a definable $C^1$ isotopy $\tau_t:N_{-1}\to\partial S_i,\ 0\le t\le1$, such that $\tau_0=\tau$ and $\tau_1$ is semialgebraic outside a small semialgebraic neighborhood of $\overline{N_{-1}-M^\epsilon_j}$ in $N_{-1}$.
Repeat the same argument for any $j$.
Then we can assume from the beginning that $\tau$ is semialgebraic outside a small semialgebraic neighborhood of the closed semialgebraic subset $X=\cap_j\overline{\partial M^\epsilon_j-\partial N_{-1}}$, of $N_{-1}$ of dimension $<m-1$.
Note $\dim X\cap\partial N_{-1}<m-2$.
Let $U$ denote the neighborhood.\par
We decrease the dimension of $X$ because the case $\dim X=-1$ is statement 2.
Consider $U,\ X$, $\{M^{2\epsilon}_j\}$ and $\{G\}$, defined by $\{M^{2\epsilon}_j\}$, in place of $\{M^\epsilon_j\}$.
Then by the same argument as above we have a definable $C^1$ isotopy $\tau_t:N_{-1}\to\partial S_i,\,0\le t\le1$, such that $\tau_0=\tau$ and $\tau_1$ is semialgebraic outside a small semialgebraic neighborhood of $X\cap\cap_j\overline{M^{2\epsilon}_j-\partial N_{-1}}$.
This set $X'$ is of dimension $<m-2$ and $\dim X'\cap\partial N_{-1}<m-3$.
Hence repeating the same argument for $3\epsilon,\ldots,m\epsilon$, we arrive at the the case $\dim X=-1$, and statement 2 is proved.
Thus statement 3 implies statement 2.\vskip1mm\noindent
{\it The induction hypothesis in the proof of statement 3.}
We prove statement 3 by induction on $m\ (=\dim N$).
Hence we assume that statement 3 holds for any proper face of $N$.\par
Consider $M^\epsilon_j$ for all $j$ in place of $N$ in statement 3.
Then the above proof also says that we can simplify statement 3 as follows:\vskip1mm\noindent
{\it Statement 3\,$'$.
Let $N^m$ be the product of a semialgebraic $C^1$ cell $C$ and the set $B=\{x\in R\times\cdots\times R\times[0,\,\infty)\times\cdots\times[0,\,\infty):|x|\le1\}$.
Let $\{G\}$ consists of $C\times\{x\in B:|x|=1\}$ and the products of the proper faces of $C$ and $B$.
Let $U$ be an open semialgebraic set in $N$ such that (i) $U=N$ or (ii) $U\subset C\times\{x\in B:1/3<|x|\}$ and $U\cap C\times\{x\in B:1/2\le|x|\}$ is of the form $C\times(U_0*0)\cap C\times\{x\in B:1/2\le|x|\}$ for some semialgebraic open set $U_0$ in $\{x\in B:|x|=1\}$.
Then for the $\tau$ in statement 3, we have the same conclusion as in statement 3.}\par
We further simplify statement 3$'$ to statement 3$''$ below in sequence.
Embed $C$ into $R^{m_1}$ for $m_1=\dim C$, and assume that $C$ is a cell and $D^{m_1}\ (=\{x\in R^{m_1}:|x|\le1\})\subset C\subset2D^{m_1}$, where $2D^{m_1}=\{2x:x\in D^{m_1}\}$.
First we can suppose that $\tau$ is the identity map on a semialgebraic neighborhood of 0 in $N$ for the following reason:\par
Let $\phi$ be a semialgebraic $C^1$ function on $R$ such that $0\le\phi\le1$, $\phi=1$ on $[0,\,1]$ and $\phi=0$ on $[2,\,\infty)$.
Let $\epsilon>0\in R$, and set
$$\tau'(x)=(d\tau)^{-1}_0\big(1-t\phi(\epsilon|x|))\tau(x)+t\phi(\epsilon|x|)(d\tau)_0x\big)\ \ \text{for }x\in N.$$
Choose sufficiently small $\epsilon$.
Then $\tau':N\to R^m,\ 0\le t\le1$, is a definable $C^1$ embedding such that $\tau'$ is the identity map on a small neighborhood of 0 in $N$.
If statement 3$'$ is proved for $\tau'$ and if $\tau'_t,\ 0\le t\le1$, is the resulting definable $C^1$ isotopy then the definable $C^1$ isotopy $\tau_t=(d\tau)_0\circ\tau'_t,\ 0\le t\le1$, satisfies the conditions in statement 3$'$ for $\tau$.
Hence replacing $\tau$ with $\tau'$, we assume that $\tau$ is the identity map and hence semialgebraic on a semialgebraic neighborhood of 0 in $N$.\par
Secondly, we treat only the case $B=D^{m-m_1}$, i.e., $X_j$, appeared in the proof that statement 3 implies statement 2, is contained in $\Int N_{-1}$ but not in $\partial N_{-1}$.
The reason is the following:
There is a semialgebraic $C^1$ isotopy $\tau'_t:N_{-1}\to N_{-1},\ 0\le t\le1$, such that $\tau'_0=\id$ and $\Ima\tau'_1\subset\Int N_{-1}$.
(The $\tau'_t,\ 0\le t\le1$, is not an isotopy of $N_{-1}$.)
Hence by considering the definable $C^1$ isotopy $\tau\circ\tau'_t:N_{-1}\to\partial S_i,\ 0\le t\le1$, we can replace the condition in statement 2 that $\tau_1$ is semialgebraic on the overall $N_{-1}$ by that $\tau_1$ is semialgebraic on a compact neighborhood of $\{0\}\times\partial D^{m-k}$ in $\Int N_{-1}$.
Then we consider only $X_j$ in $\Int N_{-1}$.
Note that if $B=D^{m-m_1}$ then $\cup G$ is the boundary of $N$, $N=(\cup G)*0$ and $\tau_t(N)=\tau(N)$ for the resultant $\tau_t,\ 0\le t\le1$, in statement 3$'$.\par
Thirdly, we can replace $N$ and $\cup G$ by $D^m$ and the sphere $S^{m-1}=\{x\in R^m:|x|=1\}$ for the following reason:
We have $S^{m-1}*0\subset N$, and each segment with ends 0 and a point of $\cup G$ intersects transversally $S^{m-1}$ once and only once.
Then keeping these properties we can modify $S^{m-1}$ through a semialgebraic $C^1$ isotopy of $N$ so that $N-D^m$ is a small neighborhood of $\cup G$ (see smoothing of corners in the proof of Theorem VI.2.1 in \cite{S1}).
Hence, as the required $\tau_t$ satisfies the condition $\tau_t=\tau$ on a small neighborhood of $\cup G$, it suffices to consider $D^m$ in place of $N$.\par
Fourthly, we can assume
$$|d\tau_x v-d\tau_{x'}v|\le|d\tau_x v|/2\quad\text{for }x,x'\in N\text{ and }v\in R^m.\leqno(*)$$
Indeed, there is a finite definable open covering $\{Q_i\}$ of $N$ such that if $x$ and $x'$ are contained in one $Q_i$ then $(*)$ is satisfied.
Hence, repeating the same argument as above, we can assume $(*)$ for any $x$ and $x'$ in $N$.
From $(*)$ it follows that the angle of $d\tau_x v$ and $d\tau_{x'}v$ is not larger than $\pi/6$.\par
By the first and fourth simplification, we have
$$1/2\le\big|\frac{d\tau\,}{d x_i}\big|\le3/2\ \ \text{and }\ |x|/2\le|\tau(x)-x|\le3|x|/2\ \ \text{for }x=(x_1,\ldots,x_m)\in N.\leqno(**)$$
Under this condition we will prove a sort of the smooth Schoenflies problem even if $m=4$, which is the key lemma of statement 3$''$ below.\par
In $(*)$ and $(**)$, the numbers $1/2$ and $3/2$ are not strict.
We will adequately change them in the forthcoming argument without mentioning it.
For example, the first inequality in $(**)$ is replaced by one $c<|d\tau/d x_i|<1/c$ for some $c>0\in\R$.
Note that if $\tau$ is semialgebraic and defined by polynomial functions with coefficients in $\R$ then $(**)$ is satisfied.
We add the following condition to the conclusion in statements 3 and 3$'$ so that the resultant of the induction hypothesis satisfies the condition:\vskip1mm\noindent
{\it Added condition.} The $\tau_t$ for any $t$ is sufficiently close to $\tau$ in the $C^1\ \R$-topology.\par
Consequently, we have simplified statement 3$'$ as follows:\vskip1mm\noindent
{\it Statement 3\,$''$.
Let $N$ be the unit disk $D^m$, and $\cup G$ the unit sphere $S^{m-1}$.
Let $U$ be such that $U\subset D^m-D^m/3$ and $U-D^m/2$ is of the form $U_0*0-D^m/2$ for some semialgebraic open set $U_0$ in $S^{m-1}$, where $D^m/k=\{x/k:x\in D^m\}$.
Let the $\tau$ in statement 3 be the identity map on $D^m/3$ and satisfy $(*)$ and $(**)$.
Then we have the same conclusion as in statement 3.}\vskip1mm\noindent
{\it Proof of statement 3\,$''$.} Set $M=\Ima\tau$, which is a compact definable $C^1$ manifold with boundary.
Note that $M=0*\partial M$ and, moreover, any line in $R^m$ passing through 0 is transversal to $\partial M$ by $(*)$ and $(**)$.
First we show the following statement:\vskip1mm\noindent
{\it Set $\tau_S=\tau|_{S^{m-1}}$.
There exist a definable $C^1$ isotopy $\tau_{S,t}:S^{m-1}\to R^m$ and a definable $C^1$ homotopy $\theta_t:S^{m-1}\to P^{m-1}(R),\ 0\le t\le1$, such that $\tau_{S,0}=\tau_S$, $\theta_t(x)$ is transversal to $\tau_{S,t}(S^{m-1})$ at $\tau_{S,t}(x)$ for each $x\in S^{m-1}$, $\tau_{S,1}$ and $\theta_1$ are semialgebraic, $\tau_{S,t}$ and $\theta_t$ are sufficiently close to $\tau_S$ and $\theta_0$, respectively, in the $C^1\ \R$-topology.
Set
$$W=\{\tau_{S,1}(x)+y:x\in S^{m-1},\,y\in|\theta_1(x)|,\,|y|\le\epsilon\}\ \text{ for small }\epsilon>0\in\R,$$
where $|\theta_1(x)|$ means the line in $R^m$.
Let a map $p:W\to\tau_{S,1}(S^{m-1})$ be defined by $p(\tau_{S,1}(x)+y)=\tau_{S,1}(x)$.
Then $p:W\to\tau_{S,1}(S^{m-1})$ is a semialgebraic $C^1$ tubular neighborhood of $\tau_{S,1}(S^{m-1})$ in $R^m$, and we can choose $\tau_{S,t}$ and $\theta_t$ so that $\Ima\tau_{S,t}\subset W$ and $\tau_{S,t}\circ\tau^{-1}_{S,1}\circ p:W\to\tau_{S,t}(S^{m-1})$ is a definable $C^1$ tubular neighborhood of $\tau_{S,t}(S^{m-1})$ in $R^m$.
We call $\tau_S\circ\tau^{-1}_{S,1}\circ p:W\to M$ a {\rm definable $C^1$ tubular neighborhood with semialgebraic direction}.}\par
We prove this statement.
We naturally regard $P^{m-1}(R)$ as a Nash manifold over $R$.
Define $\theta_0$ of class definable $C^0$ so that $\theta_0(x)$ is orthogonal to $M$ at $\tau_S(x)$ for each $x\in S^{m-1}$, approximate it by a definable $C^1$ map in the $C^0$ topology by Theorem II.5.2 in \cite{S2}, and use the same notation $\theta_0$.
Moreover, by the polynomial approximation theorem we can approximate $\theta_0$ by semialgebraic $C^1$ $\theta_1$ in the $C^1\ \R$-topology.
Then we have a definable $C^1$ homotopy $\theta_t:S^{m-1}\to P^{m-1}(R)$ such that $\theta_t$ is close to $\theta_0$ in the $C^1\ \R$-topology, and by $(**)$, $\theta_t(x)$ is transversal to $M$ at $\tau_S(x)$ for each $x\in S^{m-1}$.
For the same reason as above we have a finite semialgebraic open covering $\{Q_k\}$ of $S^{m-1}$ and for each $k$ there is a linear map $q_k:R^m\to R^{m-1}$ such that $q_k\circ\tau_S|_{\overline{Q_k}}$ is a definable $C^1$ embedding and $q^{-1}_k(y)$ as an element of $P^{m-1}(R)$ is so close to $\theta_1(x)$ in the $C^1\ \R$-topology for $x\in\overline{Q_k}$ and $y\in q_k\circ\tau_S(\overline{Q_k})$ that the angle of $q^{-1}_k(y)$ and $|\theta_1(x)|$ is smaller than $c\pi$ for some constant $c<1/2\in\R$.
Then by the above argument, it suffices to consider $\tau_{S,t}$ only on $Q_k$ and we can replace $\tau_S$ with $q_k\circ\tau_S|_{Q_k}$.
Hence we obtain the required $\tau_{S,t}$ by the induction hypothesis of statement 3.
Then existence of $W$, i.e., of $\epsilon$ and the latter half of the statement are obvious.
Thus the statement is proved.\par
For the time being, we assume that statement 3$''$ holds in the case $\partial M=S^{m-1}$.
We will reduce $\partial M$ to $S^{m-1}$ so that statement 3$''$ in the general case follows from this assumption.
For this it suffices to find a definable $C^1$ diffeomorphism $\tau':D^m\to M$ such that $\tau'$ is semialgebraic outside a small semialgebraic neighborhood of $\overline{U_0}$ in $D^m$, $\tau'=\id$ on $D^m/3$, $\tau'=\tau$ on a small neighborhood of $S^{m-1}$ in $D^m$ and $\tau'$ satisfies $(*)$ and $(**)$.
Indeed, assume that such a $\tau'$ exists.
Then the map $\tau^{\prime-1}\circ\tau:D^m\to D^m$ satisfies the conditions on $\tau$.
Therefore, statement 3$''$ holds for the $\tau^{\prime-1}\circ\tau:D^m\to D^m$, and there exists a definable $C^1$ isotopy $\tau''_t:D^m\to D^m,\ 0\le t\le1$, such that $\tau''_0=\tau^{\prime-1}\circ\tau$, $\tau''_t=\tau^{\prime-1}\circ\tau$ outside a small semialgebraic neighborhood of $S^{m-1}$, $\tau''_1$ is semialgebraic outside a smaller semialgebraic neighborhood of $\overline{U_0}$ and $\tau''_t$ for any $t$ is sufficiently close to $\tau^{\prime-1}\circ\tau$ in the $C^1$ $\R$-topology.
Hence the definable $C^1$ isotopy $\tau'\circ\tau''_t:D^m\to R^m,\ 0\le t\le1$, is the required one in statement 3$''$.\par
Setting construction of $\tau'$ aside, meanwhile, we suppose that $M$ and hence $\tau(S^{m-1})$ are semialgebraic.
Then in the same way as in the proof of the above statement we see that $\tau(S^{m-1})$ is semialgebraically $C^1$ diffeomorphic to $S^{m-1}$ through a semialgebraic $C^1$ diffeomorphism satisfying $(**)$.
Under this condition we will prove that $M$ is semialgebraically $C^1$ diffeomorphic to $D^m$ through a semialgebraic $C^1$ diffeomorphism satisfying $(**)$, which is a semialgebraic version of the smooth Schoenflies problem over $R$.
This is the present aim to obtain the note below.
We use $(**)$.
Let $R\supset\R$ for simplicity of notation, and assume $R\not=\R$ because statement 3$''$ is obvious in the Archimedean case.
The $M$ is a semialgebraic $C^1$ manifold with boundary.
Hence, by Theorem III.1.3 in \cite{S1} and Corollary 3.9 in \cite{C-S1}, $M$ is semialgebraically $C^1$ diffeomorphic to a Nash manifold with boundary defined by polynomials with coefficients in $\R$.
Moreover, by their proofs using the handle-body theory, we obtain a semialgebraic $C^1$ embedding $M\to R^m$ satisfying $(*)$ and $(**)$ such that $M$ is carried to such a Nash manifold with boundary.
Hence we assume from the beginning that $M$ is such a Nash manifold with boundary.
Set $\tau(S^{m-1})_{\R}=\tau(S^{m-1})\cap\R^m$ and define $M_{\R}$ and $D^m_{\R}$ in the same way.
Then it suffices to prove that $M_{\R}$ is semialgebraically $C^1$ diffeomorphic to $D^m_{\R}$ through a semialgebraic $C^1$ diffeomorphism $M_{\R}\to D^m_{\R}$ satisfying $(**)$ because such a diffeomorphism is extended to a semialgebraic $C^1$ diffeomorphism from $M$ to $D^m$ satisfying $(**)$.\par
We define the diffeomorphism.
Let a semialgebraic $C^1$ vector field $v$ on $R^n-\Int D^m/3$ be defined by
$$v_x=\sum_{i=1}^m x_i\frac{\partial\ }{\partial x_i}/|x|\ \ \text{for }x=(x_1,\ldots,x_m)\in R^m-\Int D^m/3.$$
Set $v_{\R}=v|_{\R^m-\Int D^m/3}$.
Then each integral curve of $v_{\R}$ intersects with each of $\partial D^m_{\R}/3$ and $\tau(S^{m-1})_{\R}$ once and only once.
Define a semialgebraic $C^1$ diffeomorphism $\alpha:M_{\R}\to D^m_{\R}$ so that the image under $\alpha$ of the integral curve in $M_{\R}$ passing through $y\in\partial D^m_{\R}/3$ is the integral curve in $D^m_{\R}$ passing through the same point and $\alpha$ carries the integral curve linearly.
Thus $M$ is semialgebraically $C^1$ diffeomorphic to $D^m$ through a semialgebraic $C^1$ diffeomorphism satisfying $(**)$.\vskip1mm\noindent
{\it Note. If $M$ is semialgebraic then statement 3\,$''$ holds.}\par
The reason is the following:
We have seen that there is a semialgebraic $C^1$ diffeomorphism $\eta:M\to D^m$ such that $\eta=\id$ on $D^m/3$ and $(**)$ is satisfied.
Consider $\eta\circ\tau:D^m\to D^m$ in place of $\tau:D^m\to M$.
Then by the hypothesis that statement 3$''$ holds if $\partial M=S^{m-1}$ there exists a definable $C^1$ isotopy $\hat\tau_t:D^m\to D^m,\ 0\le t\le1$, close to $\eta\circ\tau$ in the $C^1\ \R$-topology such that $\hat\tau_0=\eta\circ\tau$, $\hat\tau_t=\id$ on $D^m/3$, $\hat\tau=\eta\circ\tau$ on a small semialgebraic neighborhood of $S^{m-1}$ and $\hat\tau_1$ is semialgebraic outside a small definable neighborhood of $U_0$.
Hence $\eta^{-1}\circ\hat\tau_t:D^m\to M,\ 0\le t\le1$, is the required definable $C^1$ isotopy.\par
We prepare for construction of $\tau'$.
Let $0<\epsilon'\ll\epsilon\ll\epsilon''\in\R$ be sufficiently small.
Then we can assume $\tau(s x)=s\tau(x)$ for $(x,s)\in S^{m-1}\times[1-6\epsilon,\,1]$ for the following reason:\par
First we reduce the problem to the case that the map $S^{m-1}\times[1-7\epsilon,\,1]\ni(x,s)\to\tau(s x)-s\tau(x)\in R^m$ is close to 0 in the $C^1\ \R$-topology.
Let $\beta$ be a semialgebraic $C^1$ diffeomorphism of $[0,\,1]$ defined by polynomial functions with coefficients in $\R$ such that $\beta=\id$ on $[0,\,1/3]$, $\beta(1-6\epsilon)=1-\epsilon'$ and $\beta|_{[1-\epsilon'',\,1]}$ is linear.
Set $\check\tau(x)=\tau(\beta(|x|)x)/(\beta\circ\tau(\beta(|x|)x))$ for $x\in D^m$.
Then $\check\tau$ is a definable $C^1$ diffeomorphism from $D^m$ to $M$ satisfying the conditions on $\tau$ and such that the map $S^{m-1}\times[1-7\epsilon,\,1]\ni(x,s)\to\check\tau(s x)-s\check\tau(x)\in R^m$ is close to 0 in the $C^1\ \R$-topology.
If we consider $\check\tau$ in place of $\tau$ and obtain $\check\tau':D^m\to M$ as required then the map $\tau':D^m\to M$, defined by $\tau'(x)=(\beta\circ\check\tau'(x/\beta(|x|)))\check\tau'(x/\beta(|x|))$ for $x\in D^m$, satisfies the conditions.
Hence we can suppose that the map $S^{m-1}\times[1-7\epsilon,\,1]\ni(x,s)\to\tau(s x)-s\tau(x)\in R^m$ is close to 0 in the $C^1\ \R$-topology.\par
We further modify $\tau$.
Define a definable $C^1$ diffeomorphism $\delta:D^m-\Int(1-7\epsilon)D^m\to M-\Int(0*(1-7\epsilon)\partial M)$ by $\delta(s x)=s\tau(x)$ for $(x,s)\in S^{m-1}\times[1-7\epsilon,\,1]$, which is close the map $\tau|_{D^m-\Int(1-7\epsilon)D^m}:D^m-\Int(1-7\epsilon)D^m\to M$ in the $C^1\ \R$-topology.
Using $\delta$ and a semialgebraic $C^1$ partition of unity we easily modify $\tau$ so that $\tau(s x)=s\tau(x)$ on $S^{m-1}\times[1-6\epsilon,\,1]$.\par
Now we construct the $\tau':D^m\to M$.
We change the circumstance so that we can apply the above note.
Let $p:W\to\partial M$ be a definable $C^1$ tubular neighborhood of $\partial M$ in $R^m$ with semialgebraic direction.
Here we can choose $p$ so that $p^{-1}(y)$ is contained in the line  in $R^m$ passing through $y$ and 0 for each $y\in\partial M$ by the proof of the statement at the beginning of the proof of statement 3$''$.
Define a definable $C^1$ embedding $\rho:\partial M\times[-1/2,\,1]\to R^m$ so that $\Ima\rho=W$, $\rho^{-1}(M)=\partial M\times[-1/2,\,0]$, $p\circ\rho:\partial M\times[-1/2,\,1]\to\partial M$ is the projection and
$$\rho(\tau(x),s-1)=\tau(s x)\ \ \text{for }(x,s)\in S^{m-1}\times[1-6\epsilon,\,1].$$
We regard $\tau|_{\tau^{-1}(W)}:\tau^{-1}(W)\to M\cap W$ as the map $S^{m-1}\times[-1/2,\,0]\ni(x,s)\to\rho^{-1}\circ\tau((1-s)x)\in\partial M\times[-1/2,\,0]$ and translate the problem of construction of $\tau':D^m\to M$ to that on $S^{m-1}\times[-1/2,\,0]\to\partial M\times[-1/2,\,0]$.\par
By the statement at the beginning of the proof of statement 3$''$ we have a definable $C^1$ isotopy $\partial\tau_t:S^{m-1}\to W,\ 0\le t\le1$, such that $\partial\tau_0=\tau|_{S^{m-1}}$, $\Ima\partial\tau_t\subset\rho(\partial M\times[-\epsilon',\,\epsilon'])$, $\partial\tau_1$ is semialgebraic, $p\circ\partial\tau_t:S^{m-1}\to\partial M,\ 0\le t\le1$, is a definable $C^1$ isotopy and $\partial\tau_t$ for any $t$ is close to $\tau|_{S^{m-1}}$ in the $C^1$ $\R$-topology.
We write $\rho^{-1}\circ\partial\tau_t(x)=(p\circ\partial\tau_t(x),\partial\tau_{2,t}(x))$ for $x\in(1-\epsilon)S^{m-1}$.
Then $\partial\tau_{2,t}:S^{m-1}\to[-\epsilon',\,\epsilon'],\ 0\le t\le1$, is a definable $C^1$ homotopy close to the zero function in the $C^1$ $\R$-topology such that $\partial\tau_{2,0}=0$.
Using this we define a definable homeomorphism $\omega:S^{m-1}\times[-1/2,\,0]\to\partial M\times[-1/2,\,0]$ by
$$\omega(x,s)=\left \{
\begin{array}{l}
(\tau(x),s)\qquad\qquad\quad\quad\ \ \,\quad\text{for }(x,s)\in S^{m-1}\times[-\epsilon,\,0]\\
(p\circ\partial\tau_{(-\epsilon-s)/\epsilon}(x),s)\qquad\ \ \text{for }(x,s)\in S^{m-1}\times[-2\epsilon,\,-\epsilon]\\
(p\circ\partial\tau_1(x),s+\partial\tau_{2,1}(x))\quad\text{for }(x,s)\in S^{m-1}\times[-4\epsilon,\,-3\epsilon]\\
(p\circ\partial\tau_{(6\epsilon+s)/\epsilon}(x),s)\quad\qquad\text{for }(x,s)\in S^{m-1}\times[-6\epsilon,\,-5\epsilon]\\
(p\circ\tau((1-s)x),s)\quad\quad\quad\ \ \text{for }(x,s)\in S^{m-1}\times[-1/2,\,-6\epsilon]
\end{array}
\right.$$
and so that $\omega$ is linear on each of segments $\{x\}\times[-3\epsilon,\,-2\epsilon]$ and $\{x\}\times[-5\epsilon,\,-4\epsilon]$.
Here $\omega$ is piecewise smooth but not smooth at $S^{m-1}\times\{-k\epsilon:k=1,\ldots,6\}$.
However, we can smooth there and semialgebraically smooth outside the product of a small semialgebraic neighborhood $U_1$ of $\overline{U_0}$ in $S^{m-1}$ and $\{-k\epsilon:k=1,\ldots,6\}$ by the usual method fixing on $S^{m-1}\times([-1/2,\,-6\epsilon]\cup[-4\epsilon,\,-3\epsilon]\cup[-\epsilon,\,0])$.
Thus we obtain a definable $C^1$ diffeomorphism $\omega:S^{m-1}\times[-1/2,\,0]\to\partial M\times[-1/2,\,0]$ close to the map $S^{m-1}\times[-1/2,\,0]\ni(x,s)\to(p\circ\tau((1-s)x),s)\in\partial M\times[-1/2,\,0]$ in the $C^1$ $\R$-topology.\par
The $\rho\circ\omega$ is semialgebraic on $S^{m-1}\times[-4\epsilon,\,-3\epsilon]$ since $\tau(s x)=s\tau(x)$ for $(x,s)\in S^{m-1}\times[1-6\epsilon,\,1]$.
However, $\omega$ is not yet a solution of the replaced problem of construction of $\tau'$ since $\rho\circ\omega$ is not necessarily semialgebraic on $S^{m-1}\times[-1/2,\,-4\epsilon]$.
To solve this problem we return to the original problem of existence of $\tau'$.
Define a definable $C^1$ diffeomorphism $\tilde\tau:D^m-\Int D^m/2\to\rho(\partial M\times[-1/2,\,0])\ (=\tau(D^m-\Int D^m/2))$ by
$$\tilde\tau((1+s)x)=\rho\circ\omega(x,s)\ \ \ \text{for }(x,s)\in S^{m-1}\times[-1/2,\,0].$$
Then $\tilde\tau=\tau$ on $(D^m-\Int(1-\epsilon)D^m)\cup((1-6\epsilon-\epsilon')D^m-\Int D^m/2)$, $\tilde\tau$ is semialgebraic on $\big((S^{m-1}-U_1)*0-\Int D^m/2\big)\cup\big((1-3\epsilon)D^m-\Int(1-4\epsilon)D^m\big)$, and $\tilde\tau$ is close to $\tau$ in the $C^1$ $\R$-topology.
Hence we obtain a definable $C^1$ diffeomorphism $\tilde\tau:D^m\to M$ by setting $\tilde\tau=\tau$ on $D^m/2$.
Then we can assume from the beginning that $\tau$ is semialgebraic on $\big((S^{m-1}-U_1)*0-\Int D^m/2\big)\cup\big((1-3\epsilon)D^m-\Int(1-4\epsilon)D^m\big)\cup D^m/3$.\par
We modify $\tau$ to the $\tau'$, i.e., modify it semialgebraic outside $U_1*0-\Int(1-4\epsilon)D^m$.
The $\tau(D^m)\ (=M)$ is not necessarily semialgebraic.
However, $\tau((1-3\epsilon)D^m)$ is semialgebraic.
Hence by the above note, statement 3$''$ holds for $\tau|_{(1-3\epsilon)D^m}$.
Let $\tau_t:(1-3\epsilon)D^m\to\tau((1-3\epsilon)D^m),\ 0\le t\le1$, be a resulting definable $C^1$ isotopy.
We can extend it to $\tau_t:D^m\to M$ since $\tau_t=\tau$ on a small neighborhood of $(1-3\epsilon)S^{m-1}$ in $(1-3\epsilon)D^m$.
Then $\tau_1:D^m\to M$ is a definable $C^1$ diffeomorphism and semialgebraic outside a small semialgebraic neighborhood of $\overline{U_0}$ in $D^m$, $\tau_1=\tau$ on a small semialgebraic neighborhood of $S^{m-1}$ in $D^m$, and $\tau_1$ is close to $\tau$ in the $C^1$ $\R$-topology.
Thus $\tau_1$ satisfies the conditions on $\tau'$, and statement 3$''$ holds under the hypothesis that it holds in the case $\partial M=S^{m-1}$.\par
It remains to prove statement 3$''$ under the condition $\partial M=S^{m-1}$.
Assume that $\partial M=S^{m-1}$.
Then we can replace the problem to the case where $\tau$ is close to the identity map in the $C^1\ \R$-topology for the following reason:\par
By the induction hypothesis of statement 3 we can approximate $\tau|_{S^{m-1}}$ by a semialgebraic $C^2$ diffeomorphism $\partial\tau$ of $S^{m-1}$ in the $C^1\ \R$-topology.
Let $\tau_1:D^m\to R^m$ be a semialgebraic $C^2$ extension of $\partial\tau$, which is not necessarily an inclusion.
Let $\theta$ be the Nash function on $R^m$ defined by $\theta(x)=1-|x|^2$ for $x\in R^m$.
Then $\theta^{-1}(0)=S^{m-1}$, $\theta$ is non-critical at $S^{m-1}$, and the map $(\tau-\tau_1)/\theta:\Int D^m\to R^m$ is extended to a definable $C^1$ map $\tau_2$ from $D^m$ to $R^m$.
Let $\tau_3:D^m\to R^m$ be a polynomial approximation of $\tau_2$ in the $C^1\ \R$-topology.
Then $\tau_1+\theta\tau_3$ is a semialgebraic $C^1$ approximation of $\tau$ in the $D^1\ \R$-topology, and clearly it is a diffeomorphism from $D^m$ to itself.
We can consider $(\tau_1+\theta\tau_3)^{-1}\circ\tau$ in place of $\tau$.
Hence we assume that $\tau$ is close to the identity map in the $C^1\ \R$-topology.\par
Under this assumption, statement 3$''$ is obvious.
Indeed, by using the polynomial approximation theorem and a semialgebraic $C^1$ partition of unity, we can find a definable isotopy $\tau_t,\ 0\le t\le1$, so that $\tau_t=\id$ on $D^m-(1-\epsilon)D^m$ and $\tau_1$ is semialgebraic outside a small neighborhood of $\overline{U_0}$.
Thus the proofs of statement 3$''$ and hence of the existence in the first statement of Theorem 2$'$ are complete.\vskip1mm
We prove the uniqueness in the first statement of Theorem 2$'$.
What we prove is that if two semialgebraic $C^1$ manifolds are definably $C^1$ diffeomorphic then they are semialgebraically $C^1$ diffeomorphic.
We have already proved this in the above argument.
Hence the first statement of Theorem 2$'$ is proved.\qed\vskip2mm\noindent
{\it Proof of the second statement of Theorem 2\,$'$.}
Theorem VI.2.1 in \cite{S1} says that any noncompact semialgebraic $C^1$ manifold is semialgebraically $C^1$ diffeomorphic to the interior of a compact Nash manifold with boundary.
Hence, by the first statement, a noncompact definable $C^1$ manifold is definably $C^1$ diffeomorphic to the interior of some compact Nash manifold with boundary.
Its uniqueness is shown in Theorem VI.2.2 in \cite{S1}.\qed\vskip2mm\noindent
{\it Proof of Theorem 1\,$'$.}
We can prove Theorem 1$'$ by Theorem 2$'$ likewise Theorem 1 by Theorem 2.
We omit the details.\qed\vskip2mm\noindent
{\it Proof of the corollary.} (i) The statement (i) in the case $R=R'$ quickly follows from Uniqueness theorem of definable triangulation and Theorem 4.\par
Consider the case $R\not=R'$ and the ``only if\," part in (i).
Let $X\in S_m,\ Y\in S_n$ and $f:X\to Y$ a $\{S_n\}$-definable homeomorphism.
We will reduce the problem to the case where $m=n+1$, $f$ is the restriction to $X$ of the projection $R^{n+1}\to R^n$ forgetting the last factor and $X$ is the graph of some $\{S_n\}$-definable $C^0$ function on $Y$.
We set $Z=\graph f$ and naturally define $\{S_n\}$-definable $C^0$ homeomorphisms $g_X:Z\to X$ and $g_Y:Z\to Y$.
Then we have $\{S'_n\}$-definable $C^0$ maps $g^{R'}_X:Z^{R'}\to R^{\prime m}$ and $g^{R'}_Y:Z^{R'}\to R^{\prime n}$ by the conditions on $X^{R'}$, which are the respective restrictions of the projections $R^{\prime m}\times R^{\prime n}\to R^{\prime m}$ and $R^{\prime m}\times R^{\prime n}\to R^{\prime n}$.
We need to see that $g^{R'}_X$ and $g^{R'}_Y$ are homeomorphisms onto $X^{R'}$ and $Y^{R'}$ respectively.
Here ignoring $g^{R'}_X$ we can assume that $X\subset Y\times R^m\subset R^n\times R^m$, $f:X\to Y$ is the restriction to $X$ of the projection $Y\times R^m\to Y$ and $X$ is the graph of some $\{S_n\}$-definable $C^0$ map from $Y$ into $R^m$.
Moreover, we can reduce the problem to the case $m=1$ by induction on $m$.
Thus the problem to solve is the case where $X$ is the graph of some $\{S_n\}$-definable $C^0$ function.
Then it suffices to show that $X^{R'}$ is the graph of some $\{S'_n\}$-definable $C^0$ function.
This coincides with one of the conditions on $X^{R'}$.
Hence the ``\,only if\," part is proved.\par
Next we prove the ``\,if\," part in (i) in the case $R\not=R'$.
Let $X$ and $Y$ be $\{S_n\}$-definable sets such that $X^{R'}$ and $Y^{R'}$ are $\{S'_n\}$-definably homeomorphic.
These $X$ and $Y$ are $\{S_n\}$-definably homeomorphic to standard or compact semilinear sets $X_1$ and $Y_1$, respectively, by Theorem 4 and Uniqueness theorem of definable triangulation.
The above proof says that $X^{R'}$ and $Y^{R'}$ are $\{S'_n\}$-definably homeomorphic to $X^{R'}_1$ and $Y^{R'}_1$ respectively.
Hence it suffices to see that $X_1$ and $Y_1$ are $\{S_n\}$-definably homeomorphic, i.e., we can replace $X$ and $Y$ with $X_1$ and $Y_1$ respectively.
The $X^{R'}_1$ and $Y^{R'}_1$ are standard or compact semilinear sets and semilinearly homeomorphic by Theorem 4 and Uniqueness theorem of definable triangulation.
Moreover, $X^{R'}_1$ and $Y^{R'}_1$ are replaced by finite families of compact polyhedra by definition of a standard semilinear set.
In such a case, the property of being semilinearly homeomorphic does not depend on an ordered field because the problem is combinatorial.
Hence $X_1$ and $Y_1$ are semilinearly homeomorphic, which proves the ``\,if\," part.\par
Note.
Let $A\subset X$ be a $\{S_n\}$-definable sets in $R^m$.
Then $A$ is $\{S_n\}$-definable and closed in $X$ if and only if $A^{R'}$ is $\{S'_n\}$-definable and closed in $X^{R'}$.
The reason is the same as above.\par
(ii) Consider the first statement in (ii).
The case $r=0$ is obvious because as in the proof of (i) we can regard any $\{S_n\}$-definable set as a standard or compact semilinear set.
Suppose that $0<r\le\infty$ and $X$ is a $\{S_n\}$-definable $C^0$ manifold of dimension $m$ contained in $R^n$.
We will understand that the case $r=1$ is sufficient to prove.
Hence we assume $r=1$.
As the problem is local, for the proof of the ``\,only if\," part in (ii), we only need to show the following:
Let $U$ be a $\{S_n\}$-definable open subset of $R^m$ and $f$ a $\{S_n\}$-definable $C^1$ function on $U$.
Then $f^{R'}:U^{R'}\to R'$ is of class $C^1$.\par
Let $i=1,\ldots,m$.
Define a $\{S_n\}$-definable $C^0$ function $F_i$ on $U\times R$ by $F_i=0$ on $U\times\{0\}$ and
$$F_i(y,t)=\frac{f(y+(0,\ldots,0,\overset{i}{t},0,\ldots,0))-f(y)}{t}-\frac{\partial f}{\partial y_i}(y)\ \text{ on }U\times(R-\{0\}).$$
Then $F_i^{R'},\ f^{R'}$ and $(\partial f/\partial y_i)^{R'}$ are $\{S'_n\}$-definable $C^0$ functions, $F_i^{R'}=0$ on $U^{R'}\times\{0\}$ and
$$F_i^{R'}(y,t)=\frac{f^{R'}(y+(\ldots,0,t,0,\ldots))-f^{R'}(y)}{t}-(\frac{\partial f}{\partial y_i})^{R'}(y)\text{ on }U^{R'}\times(R'-\{0\}).$$
Hence $\partial(f^{R'})/\partial y_i$ exists, equals $(\partial f/\partial y_i)^{R'}$ and is continuous.
Thus $f^{R'}$ is of class $C^1$, and the ``\,only if\," part is proved.\par
We prove the ``\,if\," part.
Assume that $X^{R'}$ is of class $C^1$.
Let $x_0\in X$.
Then we will prove that $X$ is $C^1$ smooth at $x_0$.
For this we reduce the problem to the case where $X$ is the graph of some $\{S_n\}$-definable $C^0$ function defined on some $\{S_n\}$-definable open subset $R^m$.
Since $X^{R'}$ is smooth at $x_0$ there is a projection $p:R^n\ni(x_1,\ldots,x_n)\to (x_{i_1},\ldots,x_{i_m})\in R^m,\ \{i_1,\ldots,i_m\}\subset\{1,\ldots,m\}$, such that $X^{R'}$ around $x_0$ is the graph of some $\{S'_n\}$-definable $C^1$ map $V\to R^{\prime m-n}$ where $V$ is a $\{S'_n\}$-definable open subset of $R^{\prime m}$.
Hence the restriction to $p^{R'}$ of some $\{S_n\}$-definable neighborhood of $x_0$ in $X^{R'}$ is a homeomorphism onto a $\{S'_n\}$-definable neighborhood of $p^{R'}(x_0)$ in $R^{\prime m}$.
By the above argument, this property holds for $p|_X$, i.e., the restriction to $p$ of some $\{S_n\}$-definable neighborhood of $x_0$ in $X$ is a homeomorphism onto a $\{S_n\}$-definable neighborhood of $p(x_0)$ in $R^m$.
Thus we can assume from the beginning that $X$ is the graph of some $\{S_n\}$-definable $C^0$ map from a $\{S_n\}$-definable open subset $U$ of $R^m$ to $R^{n-m}$ and, moreover, of some $\{S_n\}$-definable $C^0$ function $f$ as usual.\par
We need to prove that $f$ is of class $C^1$ under the condition that $f^{R'}$ is so.
Define a $\{S_n\}$-definable $C^0$ function $G_i$ on $U\times(R-\{0\})$ by
$$G_i(y,t)=\big(f(y+(0,\ldots,0,\overset{i}{t},0,\ldots,0))-f(y)\big)/t.$$
Then $G_i^{R'}$ is continuously extended to $U^{R'}\times R'$.
We will show that $G_i$ is also continuously extended to $U\times R$ by reductio ad absurdum.
Assume that $G_i$ is not continuously extendable.
Then there are two cases to consider: (iv) $G_i(y,t)$ diverges to $\pm\infty$ as $(y,t)$ converges to some $(y_0,0)$ in some way, or (v) $G_i(y,t)$ converges to distinct numbers, say, $a$ and $b$, as $(y,t)$ converges to some $(y_0,0)$ in some two ways.
In the case (iv) we choose $\{S_n\}$-definable closed subset $Y$ of $X\times R$ such that $Y\cap X\times\{0\}=\{(y_0,0)\}$, $(y_0,0)\in\overline{Y-X\times\{0\}}$, $G_i$ does not vanish on $Y-\{(y_0,0)\}$ and $G_i(y,t)$ always diverges to $\pm\infty$ as $(y,t)$ in $Y-\{(y_0,0)\}$ converges to $(y_0,0)$.
Then, considering the function $1/G_i$ on $Y-\{(y_0,0)\}$ as above we see that $G^{R'}_i$ diverges to $\pm\infty$ as $(y,t)$ in $U^{R'}\times(R'-\{0\})$ converges to $(y_0,0)$ in some way, which is a contradiction.
In the case (v) also, replacing $G_i$ with $G_i-a$ we arrive at a contradiction as in the case (iv).
Thus $f$ is of class $C^1$, and the ``\,if\," part and the first statement are proved.\par
The second statement of (ii) is obvious because by definition, a Nash manifold and a Nash map are a semialgebraic $C^\infty$ manifold and a semialgebraic $C^\infty$ respectively.\par
(iii) By Theorem 2$'$ we can assume that $X_1$ and $X_2$ are compact Nash manifolds possibly with boundary of dimension $m$ defined by polynomial functions with coefficients in $\R_{\rm alg}$.
Moreover, if $\Int X_1$ and $\Int X_2$ are $\{S_n\}$-definably $C^1$ diffeomorphic then they are Nash diffeomorphic.
Hence we suppose from the beginning that $\{S_n\}$ and $\{S'_n\}$ are the semialgebraic structure and further $R=\R$ because the case of $R'=\R$ and $R=\R_{\rm alg}$ is easily proved.
Then what we prove is that $X_1$ and $X_2$ are semialgebraically $C^2$ diffeomorphic if and only if so are $X^{R'}_1$ and $X^{R'}_2$.
The ``\,only if\," part is obvious by the proof of (ii).\par
We prove the ``\,if\," part.
Let $\pi:X^{R'}_1\to X^{R'}_2$ be a semialgebraic $C^2$ diffeomorphism, and let $\phi$ be a nonnegative semialgebraic Morse $C^2$ function on $X_2$ such that $\phi^{-1}(0)=\partial X_2$, the critical points are all contained in $\Int X_2$ and $\phi(x)\not=\phi(x')$ for distinct critical points $x$ and $x'$.
Let $a_1<\cdots<a_k$ be the critical points of $\phi$.
Given a function $f$ on $\R^m$ of the form $f(x)=\sum_{i=1}^l x_i^2-\sum_{i=l+1}^m x_i^2+c$ for a constant $c$, we denote $N_f$ the set $\{x\in\R^m:|f(x)-c|\le\epsilon^3,|x|\le\epsilon\}$ for small $\epsilon>0\in\R$.
Then $N_f$ is a compact Nash manifold with corners, and $f$ is constant or non-critical on each proper face of $N_f$.
In the case of a Morse function $f$ on $X_2$ with one critical point also we define $N_f$ in the same way using a semialgebraic local coordinate neighborhood of the critical point.
Moreover, choosing a semialgebraic local coordinate neighborhoods of each $a_j$ we obtain a semialgebraic compact $C^2$ manifold $N_\phi$ with corners contained in $X_2$.
For the first $f$ we have a semialgebraic nonsingular $C^2$ vector field $v_f$ on $\R^n-\Int N_f$ such that $v_f f>0$ and for each proper face $F$ of $N_f$, $v_f|_F$ is a vector field on $F$ or $v_f$ is transversal to $F$, e.g.
$$(v_f)_{(x_1,\ldots,x_n)}=\sum_{i=1}^l x_i\frac{\partial\ \ }{\partial x_i}-\sum_{i=l+1}^m x_i\frac{\partial\ \ }{\partial x_i}.$$
In the same way we define a semialgebraic nonsingular $C^1$ vector field $v_\phi$ on $X_2-\Int N_\phi$ having the same properties and such that $v_\phi$ is transversal to $\partial X_2$.
Here we can choose $\phi$, $N_\phi$ and $v_\phi$ so that $|v_\phi|=1$ everywhere, $\phi$ is constant on each proper face $F$ of $N_\phi$ where $v_\phi$ is transversal and the value of each connected component of $N_\phi$ under $\phi$ is $[a_i-\epsilon,\,a_i],\, [a_i-\epsilon,\,a_i+\epsilon]$ or $[a_i,\,a_i+\epsilon]$.
Set $a_0=-1,\ X_{2,i}=\phi^{-1}([a_{i-1}+\epsilon,\,a_i-\epsilon])$ and $Y_{2,i}=\big(\phi^{-1}([a_i-\epsilon,\,a_i+\epsilon])-N_\phi\overline{\big)}$ (see Figure 5).
Then $X_{2,i}$ and $Y_{2,i}$ are compact $C^2$ manifolds with corners, their union is $X_2-\Int D_\phi$, and the vector fields $v_\phi|_{X_{2,i}}$ and $v_\phi|_{Y_{2,i}}$ give $C^1$ triviality of he functions $\phi|_{X_{2,i}}$ and $\phi|_{Y_{2,i}}$ respectively.\par
Set $N_2^{R'}=(N_\phi)^{R'},\ N_1^{R'}=\pi^{-1}(N_2)^{R'}$ and $v^{R'}_2=(v_\phi)^{R'}$, and define a semialgebraic $C^1$ vector field $v^{R'}_1$ on $X^{R'}_1-\Int N^{R'}_1$ by $d\pi(v^{R'}_{1x})=v^{R'}_{2\pi(x)}$ and replace it by $v^{R'}_{1x}/|v^{R'}_{1x}|$.
Then $X^{R'}_i,\ N^{R'}_i$ and $v^{R'}_i,\ i=1,2$, have the same properties as $X_2,\ N_\phi$ and $v_\phi$.
Moreover, the restriction of $v^{R'}_2$ to $X_2-\Int N_\phi$ is $v_\phi$.
If the restriction of $v^{R'}_1$ to $X_1-\Int N^{R'}_1$ is a vector field on $X_1-\Int N^{R'}_1$ and if $X_1\cap N^{R'}_1$ and the vector field have the same properties as $N_\phi$ and $v_\phi$ then we can see that $(X_1,N_1)$ is $C^1$ diffeomorphic and hence semialgebraically $C^2$ diffeomorphic to $(X_2,N_2)$.
Hence we modify $\pi$ and $v^{R'}_1$ so that this is the case.
First we find a semialgebraic $C^1$ isotopy $\rho_t,\ 0\le t\le1$, of $X^{R'}_1$ so that $(\pi\circ\rho_1)^{-1}(a_j)\in X_1$ for any $a_j$.
Hence we can assume $\pi^{-1}(a_j)\in X_1$.
Secondly, we choose a semialgebraic $C^1$ isotopy $\rho'_t,\ 0\le t\le1$, of $X^{R'}_1$ so that $\rho'_0=\id$, $\rho'_t=\id$ on each $\pi^{-1}(a_j)$ and $\phi\circ\pi\circ\rho'_1|_{X_1\cap N^{R'}_1}$ is a real valued Morse function.
Thus we can suppose that $X_1\cap N^{R'}_1\ (=N_1)$ has the same properties as $N_\phi$.
Here we choose $\rho_t$ and $\rho'_t$ so that the new $v^{R'}_1$ keeps the property $|v^{R'}_1|=1$.
Thirdly, keeping the  properties of $v^{R'}_1$ we replace $v^{R'}_1$ with a semialgebraic $C^1$ vector field $v_1$ close to $v^{R'}_1$ in the $C^1$ $\R$-topology whose restriction to $X_1-\Int N_1$ is a nonsingular $C^1$ vector field on $X_1-\Int N_1$, and we modify $v_1$ so that its absolute value is constant 1 as above.
Thus we obtain $N_1$ and $v_1$ which have the same properties as $N_2$ and $v_2$.\par
It remains to prove that $(X_1,N_1)$ and $(X_2,N_2)$ are $C^1$ diffeomorphic.
Let $\xi:Z\to X_1-\Int N_1,\ Z\subset(X_1-\Int N_1)\times\R$, be the union of the local one-parameter groups of local transformations generated by $v_1$.
Let $N_{1,i}$ denote the connected component of $N_1$ containing $\pi^{-1}(a_i)$.
Then modifying $v_1$ outside a small neighborhood of $\partial N_1$ in $X_1-\Int N_1$ we easily obtain compact $C^1$ manifolds $X_{1,i}$ and $Y_{1,i}$, $i=1,\ldots,k$, with corners such that the following conditions are satisfied:
$X_1-\Int N_1$ is their union; $X_{1,i}\cap X_{1,i'}=\emptyset$ and $Y_{1,i}\cap Y_{1,i'}=\emptyset$ for $i\not=i'$; $X_{1,i}\cap Y_{1,i}$ and $X_{1,i+1}\cap Y_{1,i}$ are unions of common proper faces of $X_{1,i}$, $Y_{1,i}$ and $X_{1,i+1}$; $\xi|_{Z\cap\partial X_1\times\R},\ \xi|_{Z\cap(X_{1,i}\cap Y_{1,i})\times\R}$ and $\xi|_{Z\cap(X_{1,i+1}\cap Y_{1,i})\times\R}$ are $C^1$ diffeomorphisms onto $X_{1,1},\ Y_{1,i}$ and $X_{1,i+1}$ respectively.
Hence $v_1$ and $\xi$ give $C^1$ diffeomorphisms $\pi_{X,i}:X_{1,i}\to X_{2,i}$ and $\pi_{Y,i}:Y_{1,i}\to Y_{2,i},\ i=1,\ldots,k$, such that $\pi_{X,i+1}(X_{1,i+1}\cap Y_{1,i})=X_{2,i+1}\cap Y_{2,i}$, $\pi_{X,i}(X_{1,i}\cap Y_{1,i})=X_{2,i}\cap Y_{2,i},\ \pi_{Y,i}(X_{1,i+1}\cap Y_{1,i})=X_{2,i+1}\cap Y_{2,i}$, $\pi_{Y,i}(X_{1,i}\cap Y_{1,i})=X_{2,i}\cap Y_{2,i}$.
Thus it suffices to modify $\pi_{X,i}$ and $\pi_{Y,i}$ so that $\pi_{X,i}=\pi_{Y,i}$ on $X_{1,i}\cap Y_{1,i}$ and $\pi_{X,i+1}=\pi_{Y,i}$ on $X_{1,i+1}\cap Y_{1,i}$ keeping the properties of $\pi_{X,i}$ and $\pi_{Y,i}$.
Indeed, such $\pi_{X,i}$ and $\pi_{Y,i}$ give a $C^1$ diffeomorphism from $X_1-\Int D_1$ onto $X_2-\Int D_2$, which is extendable to one from $X_1$ onto $X_2$.
In turn, we modify $\pi_{X,1},\pi_{Y,1},\pi_{X,2}$, and so on.
Then the problem to solve is only the following statement:\vskip1mm\noindent
{\it There exists a $C^1$ isotopy $\chi_{i,t},\ 0\le t\le1$, of $X_{1,i}$ such that $\chi_{i,t}=\id$ on $X_{1,i}\cap Y_{1,i-1}$ and $\pi_{X,i}\circ\chi_{i,1}(X_{1,i}\cap N_{1,i})=X_{2,i}\cap N_{2,i}$, where $N_{2,i}$ is naturally defined.}\par
\begin{figure}[h]
\begin{pspicture}(0,-0.3)(12.7,3)
\pscurve(2.6,0)(2.7,-0.15)(3.3,-0.3)(3.9,-0.15)(4,0)
\pscurve(5.3,0)(5.4,-0.15)(6,-0.3)(6.6,-0.15)(6.7,0)
\pscurve(2.6,0)(2.7,1)(4.7,2.6)(6.6,1)(6.7,0)
\pscurve(4,0)(4.1,0.4)(4.7,0.8)(5.2,0.4)(5.3,0)

\pscustom[linewidth=0.2pt]{\pscurve(2.6,0.5)(2.7,0.35)(3.3,0.25)(3.9,0.35)(4.13,0.45)}
\pscustom[linewidth=0.2pt]{\pscurve(5.17,0.45)(5.4,0.35)(6,0.25)(6.6,0.35)(6.7,0.5)}
\pscustom[linewidth=0.2pt]{\pscurve(2.85,1.3)(2.95,1.2)(4.7,1.1)(6.3,1.2)(6.4,1.3)}
\pscustom[linewidth=0.2pt]{\pscurve(3.4,2)(3.5,1.9)(4.7,1.8)(5.8,1.9)(5.9,2)}
\pscustom[linewidth=0.2pt]{\pscurve(3.8,0.33)(3.9,0.8)(4.1,1.08)}
\pscustom[linewidth=0.2pt]{\pscurve(5.5,0.33)(5.4,0.8)(5.2,1.08)}
\pscircle(4.6,2.6){0.05}
\pscircle(4.6,0.8){0.05}
\rput(2.3,1.3){$X_1$}
\rput(4.6,2.15){$N_{1,2}$}
\rput(4.6,2.85){\small$\pi^{-1}(a_2)$}
\rput(4.6,1.4){$X_{1,2}$}
\rput(3.4,0){$X_{1,1}$}
\rput(6,0){$X_{1,1}$}
\rput(3.4,0.65){$Y_{1,1}$}
\rput(6,0.65){$Y_{1,1}$}
\rput(4.95,0.9){\small$N_{1,1}$}
\rput(4.7,0.1){\small$\pi^{-1}(a_1)$}
\pscustom[linewidth=0.2pt]{\psline{->}(4.7,0.3)(4.63,0.72)}
\pscurve(8.6,0)(8.7,-0.15)(9.3,-0.3)(9.9,-0.15)(10,0)
\pscurve(11.3,0)(11.4,-0.15)(12,-0.3)(12.6,-0.15)(12.7,0)
\pscurve(8.6,0)(8.7,1)(10.7,2.6)(12.6,1)(12.7,0)
\pscurve(10,0)(10.1,0.4)(10.7,0.8)(11.2,0.4)(11.3,0)
\pscustom[linewidth=0.2pt]{\pscurve(8.6,0.5)(8.7,0.35)(9.3,0.25)(9.9,0.35)(10.13,0.45)}
\pscustom[linewidth=0.2pt]{\pscurve(11.17,0.45)(11.4,0.35)(12,0.25)(12.6,0.35)(12.7,0.5)}
\pscustom[linewidth=0.2pt]{\pscurve(8.85,1.3)(8.95,1.2)(10.7,1.1)(12.3,1.2)(12.4,1.3)}
\pscustom[linewidth=0.2pt]{\pscurve(9.4,2)(9.5,1.9)(10.7,1.8)(11.8,1.9)(11.9,2)}
\pscustom[linewidth=0.2pt]{\pscurve(9.8,0.33)(9.9,0.8)(10.1,1.08)}
\pscustom[linewidth=0.2pt]{\pscurve(11.5,0.33)(11.4,0.8)(11.2,1.08)}
\pscircle(10.6,2.6){0.05}
\pscircle(10.6,0.8){0.05}
\rput(13,1.3){$X_2$}
\rput(10.6,2.15){$N_{2,2}$}
\rput(10.6,2.8){$a_2$}
\rput(10.6,1.4){$X_{2,2}$}
\rput(9.4,0){$X_{2,1}$}
\rput(12,0){$X_{2,1}$}
\rput(9.4,0.65){$Y_{2,1}$}
\rput(12,0.65){$Y_{2,1}$}
\rput(10.95,0.9){\small$N_{2,1}$}
\rput(10.7,0.6){\small$a_1$}
\psline{->}(7.2,0)(8.2,0)
\psline{->}(7.2,0.65)(8.2,0.65)
\psline{->}(7.2,1.4)(8.2,1.4)
\rput(7.7,0.2){$\pi_{X,1}$}
\rput(7.7,0.85){$\pi_{Y,1}$}
\rput(7.7,1.6){$\pi_{X,2}$}
\end{pspicture}
\caption{Diffeomorphism from $X_1$ onto $X_2$}\end{figure}
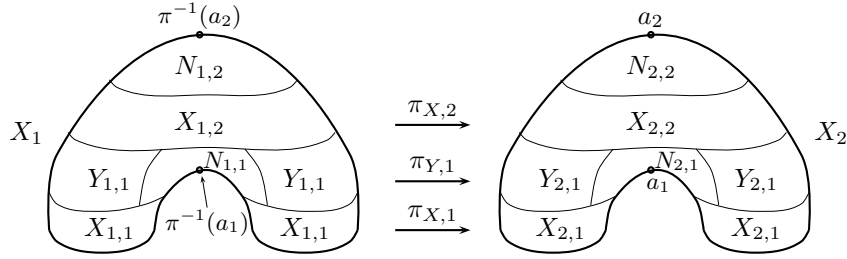\par
We can prove this statement in the same way as in the proof of Theorem 2$'$ but more easily.
We omit the details.
Thus we complete the proofs of the ``\,if\," part and the corollary.\qed
{\small
\qquad Graduate school of mathematics, Nagoya University, Chikusa, Nagoya, 464-8602, Japan\\
\qquad{\it E-mail address}: shiota@math.nagoya-u.ac.jp}
\enddocument